\pdfoutput=1

\documentclass[preprint,3p]{elsarticle}

\usepackage{lineno,hyperref}
\usepackage{bm}
\usepackage{amsmath}
\usepackage{color}
\usepackage{algorithm}
\usepackage{algorithmic}
\usepackage{multirow}

\usepackage{amssymb}

\usepackage{enumitem}

\usepackage{subfigure}
\usepackage{array}
\newcolumntype{C}[1]{>{\centering\arraybackslash}p{#1}}

\usepackage{xcolor}

\newcommand{\Div}[1]{\nabla \cdot {#1}}
\newcommand{\Grad}[1]{\nabla {#1}}

\newcommand{\boundary}[1]{\Gamma^{\mathrm{#1}}}
\newcommand{\hboundary}[1]{\Gamma_{h}^{\mathrm{#1}}}

\usepackage{stmaryrd} 

\newcommand{\avg}[1]{\{\!\{#1\}\!\}}

\newcommand{\jump}[1]{\llbracket {#1} \rrbracket }

\newcommand{\intele}[2]{ \left( {#1},{#2} \right)_{\Omega_{e}} }
\newcommand{\inteleface}[2]{ \left( {#1},{#2} \right)_{\partial\Omega_{e}} }
\newcommand{\intelefaceInterior}[2]{ \left( {#1},{#2} \right)_{\partial\Omega_{e}\setminus\Gamma_h }}

\newcommand{\intelefaceDirichlet}[2]
{ {\left( {#1},{#2} \right)}_{\partial\Omega_{e}\cap \Gamma^{\rm{D}}_{h}}}
\newcommand{\intelefaceNeumann}[2]{ {\left( {#1},{#2} \right)}_{\partial\Omega_{e}\cap \Gamma^{\rm{N}}_{h}}}

\newenvironment{remark}[1][Remark]{\begin{trivlist}
\item[\hskip \labelsep {\bfseries #1}]}{\end{trivlist}}

\journal{}

\begin{document}

\begin{frontmatter}

\title{On the stability of projection methods\\ for the incompressible Navier--Stokes equations\\ based on high-order discontinuous Galerkin discretizations}

\author{Niklas Fehn}
\ead{fehn@lnm.mw.tum.de}
\author{Wolfgang A. Wall}
\ead{wall@lnm.mw.tum.de}
\author{Martin Kronbichler\corref{correspondingauthor1}}
\cortext[correspondingauthor1]{Corresponding author at: Institute for Computational Mechanics, Technical University of Munich, Boltzmannstr. 15, 85748 Garching, Germany. Tel.: +49 89 28915300; fax: +49 89 28915301}
\ead{kronbichler@lnm.mw.tum.de}
\address{Institute for Computational Mechanics, Technical University of Munich,\\ Boltzmannstr. 15, 85748 Garching, Germany}

\begin{abstract}
The present paper deals with the numerical solution of the incompressible Navier--Stokes equations using high-order discontinuous Galerkin (DG) methods for discretization in space. For DG methods applied to the dual splitting projection method, instabilities have recently been reported that occur for coarse spatial resolutions and small time step sizes. By means of numerical investigation we give evidence that these instabilities are related to the discontinuous Galerkin formulation of the velocity divergence term and the pressure gradient term that couple velocity and pressure. Integration by parts of these terms with a suitable definition of boundary conditions is required in order to obtain a stable and robust method. Since the intermediate velocity field does not fulfill the boundary conditions prescribed for the velocity, a consistent boundary condition is derived from the convective step of the dual splitting scheme to ensure high-order accuracy with respect to the temporal discretization. This new formulation is stable in the limit of small time steps for both equal-order and mixed-order polynomial approximations. Although the dual splitting scheme itself includes inf--sup stabilizing contributions, we demonstrate that spurious pressure oscillations appear for equal-order polynomials and small time steps highlighting the necessity to consider inf--sup stability explicitly.
\end{abstract}


\end{frontmatter}


\section{Introduction}\label{Intro}
The numerical solution of the incompressible Navier--Stokes equations is a key issue in computational fluid dynamics. With respect to discretization in time and space, two aspects are of primary importance regarding the present work. On the one hand, operator splitting techniques are well established solution approaches for the incompressible Navier--Stokes equations that are particularly efficient for high Reynolds number flows~\cite{Karniadakis13}. On the other hand, high-order discontinuous Galerkin methods have gained significance as compared to state-of-the-art discretization methods like Finite Volume Methods and Finite Element Methods. They exhibit favorable properties such as high-order accuracy and hp-adaptivity, stability in convection-dominated flows, and geometric flexibility~\cite{Hesthaven07}. The present work is devoted to the problem of instabilities reported and analyzed in~\cite{Ferrer11,Ferrer14} for the discontinuous Galerkin method proposed in~\cite{Hesthaven07} that is based on the high-order dual splitting scheme developed in~\cite{Karniadakis1991}. 

\subsection{Splitting methods for the incompressible Navier--Stokes equations}
The coupling of velocity and pressure in the momentum equation in combination with the incompressibility constraint poses a major challenge in terms of the numerical solution of the incompressible Navier--Stokes equations~\cite{Guermond06}. For monolithic solution approaches, discretization in space and time leads to a system of equations of indefinite saddle point type for velocity and pressure unknowns. Splitting methods, instead, aim at seperating the computation of velocity and pressure in the solution algorithm in order to obtain a set of equations that can be solved more efficiently from a linear algebra point of view such as a convection--diffusion type problem for the velocity and a Poisson equation for the pressure.

Splitting methods for the incompressible Navier--Stokes equations may be subdivided into the four main groups of pressure-correction schemes, velocity-correction schemes, consistent splitting schemes, and algebraic splitting schemes, see~\cite{Karniadakis13,Guermond06} for a comprehensive overview. In the present work the focus is on velocity-correction methods~\cite{Karniadakis1991,Orszag1986,Guermond2003} and pressure-correction methods~\cite{Chorin68,hirt1972,Goda1979,VanKan1986,Timmermans1996,Guermond2004}. In case of pressure-correction methods the momentum equation is solved in the first step using an extrapolation of the pressure gradient term. The pressure and a divergence-free velocity field are obtained in the second step by projecting the intermediate velocity onto the space of divergence-free vectors. This order is reversed for velocity-correction schemes. While the pressure and a divergence-free velocity field are  calculated in the first step taking into account the convective term, the viscous term is considered in the second step.

When using projection methods, a Neumann boundary condition has to be prescribed for the pressure on Dirichlet boundaries. Inconsistent formulations of this pressure Neumann boundary condition can cause unphysical boundary layers which also limit the temporal accuracy of projection schemes~\cite{Karniadakis13,Guermond06}. To obtain higher order accuracy with respect to the temporal discretization, consistent formulations of the pressure Neumann boundary condition are crucial leading to so-called rotational formulations.

\subsection{Discontinuous Galerkin methods for the incompressible Navier--Stokes equations}
The local discontinuous Galerkin method (LDG) is analyzed in~\cite{cockburn2002local} for the steady Stokes equations, in~\cite{cockburn2004local} for the Oseen equation, and a locally conservative LDG method for the steady incompressible Navier--Stokes equations is proposed in~\cite{cockburn2005locally}. A stable equal-order formulation for the steady Navier--Stokes equations using a pressure-stabilization term is proposed in~\cite{cockburn2009equal} considering the local discontinuous Galerkin method and the interior penalty method for the discretization of the viscous term and an upwind flux formulation for the convective term. Regarding the discontinuous Galerkin methods proposed in~\cite{girault2005discontinuous,girault2005splitting}, both the symmetric and the non-symmetric interior penalty method are considered for the viscous term, Lesaint--Raviart upwinding fluxes for the convective term, and central fluxes for the pressure gradient term and velocity divergence term.
The DG method of~\cite{bassi2006artificial, bassi2007implicit} for unsteady incompressible flow solves a local Riemann problem to compute the inviscid numerical fluxes.

Similar approaches in terms of the DG discretization of the convective term and the viscous term are proposed in~\cite{Hesthaven07} and~\cite{Shahbazi07}. The convective term is written in divergence form to ensure local conservativity and is discretized using the local Lax--Friedrichs flux. The discretization of the viscous term is based on the symmetric interior penalty Galerkin (SIPG) method.  The temporal discretization is based  on the high-order dual splitting scheme in~\cite{Hesthaven07} and on an algebraic splitting scheme in~\cite{Shahbazi07}. This DG discretization of the convective term and viscous term is also applied in~\cite{Klein13,Klein15}. The DG discretization of the convective term and viscous term used in the present work follows the approach of~\cite{Hesthaven07,Shahbazi07}. While central fluxes are used for the velocity divergence term and pressure gradient term in~\cite{Shahbazi07,Klein13,Klein15}, no integration by parts of these terms is considered in~\cite{Hesthaven07}. As discussed below, the DG discretization of the velocity--pressure coupling terms is of particular importance with respect to the stability of the method for small time step sizes and both variants are analyzed in the present work.

DG formulations applied to pressure-correction schemes are proposed in~\cite{Botti11} using mixed discontinuous--continuous approximations for velocity and pressure and in~\cite{Piatkowski16} using a fully discontinuous formulation. A hybridizable discontinuous Galerkin method is proposed in~\cite{nguyen2011implicit} and hybrid discontinuous Galerkin methods are considered in~\cite{Lehrenfeld16} using a standard DG discretization for the convective term and an H(div)-conforming HDG discretization for the velocity occuring in the Stokes operator.

\subsection{Novel contributions of the present paper}
A discontinuous Galerkin formulation for the high-order dual splitting scheme~\cite{Karniadakis1991} has first been proposed in~\cite{Hesthaven07} using equal-order approximations for velocity and pressure. Instabilities of this method have been reported in~\cite{Ferrer11} in the limit of small time steps. These instabilities are analyzed in more detail in~\cite{Ferrer14} where the instabilities occuring for small time step sizes are related to inf--sup instabilities. 
A stabilization is proposed in~\cite{Ferrer14} by scaling the penalty parameter of the interior penalty method used to discretize the pressure Poisson equation by the inverse of the time step size. Instabilities have also been reported in~\cite{Steinmoeller13} for the under-reolved and low viscosity regime where a stabilization is proposed that is based on a postprocessing step projecting the intermediate velocity onto the space of exactly divergence-free vectors. The instabilities considered in~\cite{Ferrer14} and~\cite{Steinmoeller13} might also be related as discussed in~\cite{Krank16b}, where several stabilization approaches are reviewed and where a stabilization similar to the postprocessing in~\cite{Steinmoeller13,Joshi16} is proposed by adding a consistent div--div penalty term to the projection equation. 
The analysis in~\cite{Emamy14} suggests that the instabilities might be related to the discretization of the velocity divergence term and pressure gradient term, but the formulation is unclear regarding the imposition of boundary conditions. The approach proposed recently in~\cite{Emamy17} to overcome these instabilities has already been proposed by~\cite{Leriche2000,Leriche2006} in a different context and has been analyzed in~\cite{Krank16b} in the context of instabilities for small time step sizes and discontinuous Galerkin discretizations. Instabilities for the standard projection are also reported in~\cite{Piatkowski16} in the context of pressure-correction methods, where the div--div penalty based projection is compared to a postprocessing technique using~$H(\mathrm{div})$ reconstruction with Raviart--Thomas spaces. To the best of the authors' knowledge, no clear understanding of the instabilities for small time step sizes is currently available. Our numerical investigations indicate that these instabilities are neither related to inf--sup instabilities nor to inaccuracies of the spatially discretized projection operator resulting in velocity fields that do not exactly fulfill the divergence-free constraint.

For the DG formulation proposed in~\cite{Hesthaven07} and analyzed in~\cite{Ferrer11,Ferrer14}, the velocity divergence term and pressure gradient term are not integrated by parts when deriving the weak formulation. This might be due to the following reasons:
\begin{itemize}
\item For the high-order dual splitting scheme (and more general projection methods), one is not forced to perform integration by parts of these terms as they appear on the right-hand side of the pressure Poisson equation and projection equation. Hence, the resulting systems of equations are still solvable without integration by parts. Note that this is fundamentally different for a monolithic solution approach. In that case, the system of equations is not solvable which becomes obvious when looking at the fact that performing integration by parts and defining numerical fluxes for the pressure gradient term is a necessary prerequisite to enforce continuity of the pressure solution in a weak sense. 

\item Performing integration by parts and defining numerical fluxes also requires a treatment of boundary conditions. In case of the high-order dual splitting scheme, this is not straight-forward because the intermediate velocity does not fulfill the Dirichlet boundary conditions prescribed for the velocity. Integration by parts of the velocity divergence term is mentioned in~\cite{Steinmoeller13}, however, without defining a numerical flux function and specifying boundary conditions. Integration by parts of both terms using central fluxes is considered in~\cite{Krank16b,Emamy14,Emamy17} but uncertainties with respect to the treatment of boundary conditions are avoided by defining exterior values on domain boundaries as a function of interior values only, or by using inconsistent velocity Dirichlet boundary conditions.
\end{itemize}

Using the discontinuous Galerkin formulation proposed in~\cite{Hesthaven07}, we demonstrate that we can reproduce the instabilities analyzed in~\cite{Ferrer11,Ferrer14} in the limit of small time steps and that these instabilities occur similarly for both equal-order and mixed-order approximations.
The discontinuous Galerkin formulation of the velocity--pressure coupling terms occuring on the right-hand side of the pressure Poisson equation and the projection step play a crucial role with respect to the instabilities described above. Integration by parts of these terms with consistent boundary conditions should be performed in order to obtain a stable and robust method. In this respect, we propose a stable and high-order accurate boundary condition for the intermediate velocity field. By means of numerical investigation we demonstrate that this new formulation is stable in the limit of small time steps for both equal-order and mixed-order approximations.
In addition, we show that inf--sup instabilities in form of spurious pressure oscillations are present when using equal-order polynomial approximations. Hence, our results significantly differ from the conclusions drawn in~\cite{Ferrer14} and~\cite{Emamy17}. As a means of verification of our results we compare the results for the high-order dual splitting scheme to alternative solution strategies such as a fully coupled, monolithic solution approach and pressure-correction schemes. 

\subsection{Outline}

The outline of this paper is as follows. Section~\ref{MathematicalModel} describes the mathematical model of the incompressible Navier--Stokes equations. Aspects related to the temporal discretization are discussed in Section~\ref{TemporalDiscretization} where we present different solution strategies such as a coupled solution approach and splitting methods with an emphasis on the description of boundary conditions. Section~\ref{SpatialDiscretization} is devoted to the discontinuous Galerkin discretization of the incompressible Navier--Stokes equations using a general framework for the different solution strategies. In Sections~\ref{TemporalDiscretization} and~\ref{SpatialDiscretization} we intentionally choose a comprehensive discussion of boundary conditions and weak forms for reasons of reproducibility and clarity. Numerical evidence for our hypotheses is given in Section~\ref{NumericalResults}. In Section~\ref{Conclusion} we summarize our results.

\section{Mathematical model}\label{MathematicalModel}
We consider the incompressible Navier--Stokes equations in a domain~$\Omega \subset \mathbb{R}^d$, consisting of the momentum equation
\begin{equation}
\frac{\partial \bm{u}}{\partial t} + \nabla \cdot \bm{F}_{\mathrm{c}}(\bm{u}) - \nabla \cdot \bm{F}_{\mathrm{v}} (\bm{u}) + \nabla p = \bm{f} \;\; \text{in}\; \Omega \times [0, T]\label{MomentumEquation}
\end{equation}
and the continuity equation
\begin{equation}
\nabla \cdot \bm{u} = 0 \;\; \text{in}\; \Omega \times [0, T] \; ,\label{ContinuityEquation}
\end{equation}
where the unknowns are the velocity~$\bm{u}=(u_1,...,u_d)^T$ and the kinematic pressure~$p$. The body force vector is denoted by~$\bm{f} = (f_1,...,f_d)^T$. The convective term is written in conservative (divergence) formulation where the convective flux is~$\bm{F}_{\mathrm{c}}(\bm{u}) = \bm{u}\otimes \bm{u}$. The viscous term is written in Laplace formulation so that the viscous flux is given as~$\bm{F}_{\mathrm{v}}(\bm{u})=\nu \nabla \bm{u}$ with the constant kinematic viscosity~$\nu$.

The incompressible Navier--Stokes equations~\eqref{MomentumEquation} and~\eqref{ContinuityEquation} are subject to the initial condition
\begin{equation}
\bm{u}(\bm{x}, t=0) = \bm{u}_0(\bm{x}) \;\; \text{in} \; \Omega \; ,
\end{equation}
where~$\bm{u}_0(\bm{x})$ has to be divergence-free and fulfills the velocity Dirichlet boundary condition~$\bm{g}_{u}$ described below. On the boundary ~$\Gamma = \partial \Omega$, Dirichlet and Neumann boundary conditions are prescribed
\begin{align}
\bm{u} &= \bm{g}_{u}\;\; \text{on} \; \boundary{D} \times [0, T]\; ,\label{DirichletBC}\\
\left(\bm{F}_{\mathrm{v}} (\bm{u})  - p \bm{I} \right) \cdot \bm{n} &= \bm{h}\;\; \text{on} \; \boundary{N} \times [0, T]\; ,\label{NeumannBC_Coupled}
\end{align}
where the Dirichlet and Neumann part of the boundary are denoted by~$\boundary{D}$ and~$\boundary{N}$, respectively, with~$\Gamma = \boundary{D} \cup \boundary{N}$ and~$\boundary{D} \cap \boundary{N} = \emptyset$. The outward pointing unit normal vector is denoted by~$\bm{n}$. For the projection methods considered in this work a splitting of the Neumann boundary condition according to~$\bm{h} = \bm{h}_u - g_p \bm{n}$ into a viscous part~$ \bm{h}_u$ and a pressure part~$g_p$ is necessary due to the operator splitting. Accordingly, the viscous forces and the pressure have to be prescribed seperately on~$\boundary{N}$
\begin{align}
\bm{F}_{\mathrm{v}} (\bm{u})\cdot \bm{n} &= \bm{h}_u\;\; \text{on} \; \boundary{N} \times [0, T]\; ,\label{NeumannBC_Viscous}\\
p &= g_p\;\; \text{on} \; \boundary{N} \times [0, T]\; .\label{NeumannBC_Pressure}
\end{align}
In case of pure Dirichlet boundary conditions,~$\Gamma = \boundary{D}$, the velocity Dirichlet boundary condition~\eqref{DirichletBC} has to fulfill the constraint~$\int_{\boundary{D}} \bm{g}_{u} \cdot \bm{n}\; \mathrm{d}\Gamma = 0$. In that case, the pressure is only defined up to an additive constant and a unique pressure solution can be obtained by requiring~$\int_{\Omega} p \;\mathrm{d}\Omega = 0$.

The unsteady (generalized) Stokes equations and their temporal and spatial discretization are obtained from the incompressible Navier--Stokes equations by neglecting the convective term in equation~\eqref{MomentumEquation}. In the following sections we discuss the full incompressible Navier--Stokes equations.

\section{Temporal discretization}\label{TemporalDiscretization}
As solution strategies for the incompressible Navier--Stokes equations~\eqref{MomentumEquation} and~\eqref{ContinuityEquation} we consider a fully coupled solution approach as well as projection type solution methods such as pressure-correction schemes and velocity-correction schemes. For the latter, an operator splitting of the incompressible Navier--Stokes equations is performed at the level of differential equations when discretizing the equations in time. As a time integration method we consider BDF (backward differentiation formula) time integration.

\subsection{Coupled solution approach}\label{TemporalDiscretization_CoupledSolution}
To discretize the incompressible Navier--Stokes equations~\eqref{MomentumEquation} and~\eqref{ContinuityEquation} in time, the time interval~$[0,T]$ is divided into~$N$ time steps of uniform size~$\Delta t = T/N$ leading to the grid~${\lbrace t_{i}\rbrace}_{i=0}^{N} = {\lbrace i \Delta t\rbrace}_{i=0}^{N}$. Let~$n=0,...,N-1$ denote the time step number. Then, the equations are advanced from time~$t_{n} = n \Delta t$ to time~$t_{n+1} = (n+1) \Delta t$ in time step~$n$. 
Applying the BDF scheme to the incompressible Navier--Stokes equations yields
\begin{align}
\frac{\gamma_0 \bm{u}^{n+1}-\sum_{i=0}^{J-1}\left(\alpha_i\bm{u}^{n-i}\right)}{\Delta t} + \Div{\bm{F}_{\mathrm{c}}(\bm{u}^{n+1})}
- \Div{\bm{F}_{\mathrm{v}} (\bm{u}^{n+1})} + \Grad{p^{n+1}} &= \bm{f}\left(t_{n+1}\right)\; ,\label{TemporalDiscretization_Coupled_Momentum}\\
\Div{\bm{u}^{n+1}} &= 0 \; ,\label{TemporalDiscretization_Coupled_Continuity}
\end{align}
In the present work we consider BDF schemes of order~$J=1,2$ which are A-stable time integration schemes. The coefficients of the time integration scheme are given as~$\gamma_0=\alpha_0=1$ for~$J=1$ and~$\gamma_0=3/2$,~$\alpha_0 = 2$,~$\alpha_1=-1/2$ for~$J=2$.

\subsection{High-order dual splitting scheme}
The high-order dual splitting scheme~\cite{Karniadakis1991} treats the convective term, the pressure term, and the viscous term seperately in different substeps, where the convective term is formulated explicitly and the viscous term implicitly in time. Discretization in time is based on BDF time integration and an extrapolation scheme is used in order to extrapolate the convective term to time~$t_{n+1}$.

\subsubsection{Convective step} 
In the first substep, the convective term and the body force term are considered. An intermediate velocity field~$\hat{\bm{u}}$ is obtained from the following equation
\begin{equation}
\frac{\gamma_0\hat{\bm{u}}-\sum_{i=0}^{J-1}\left(\alpha_i\bm{u}^{n-i}\right)}{\Delta t} = 
- \sum_{i=0}^{J-1}\left(\beta_i \Div{\bm{F}_{\mathrm{c}}\left(\bm{u}^{n-i}\right)}\right)
+ \bm{f}\left(t_{n+1}\right)\; ,\label{DualSplitting_ConvectiveStep}
\end{equation}
where Dirichlet boundary conditions are imposed for the velocity field on~$\boundary{D}$ at old time instants~$t_{n-i}$
\begin{equation}
\bm{u} = \bm{g}_u \;\; \text{on} \; \boundary{D} \;\; .
\end{equation}
For the convective term an exptrapolation scheme of order~$J$ is used, where the coefficients~$\beta_i$ are~$\beta_0=1$ for~$J=1$ and~$\beta_0=2$,~$\beta_1=-1$ for~$J=2$.

\subsubsection{Pressure step and projection step}
In the second substep, the pressure solution~$p^{n+1}$ at time~$t_{n+1}$ as well as a second intermediate velocity field~$\hat{\hat{\bm{u}}}$ are computed by decomposing the intermediate velocity~$\hat{\bm{u}}$ into an irrotational part~$\Grad{p^{n+1}}$ and a solenoidal part~$\hat{\hat{\bm{u}}}$ (projection method)
\begin{align}
\frac{\gamma_0 }{\Delta t}\hat{\hat{\bm{u}}}+\Grad{p^{n+1}} =\frac{\gamma_0 }{\Delta t}\hat{\bm{u}} \; ,\label{ProjectionMethodEq1}\\
\Div{\hat{\hat{\bm{u}}}} = 0 \label{ProjectionMethodEq2} \; .
\end{align}
To solve these equations a Poisson equation is derived for the pressure by taking the divergence of equation~\eqref{ProjectionMethodEq1} and making use of equation~\eqref{ProjectionMethodEq2}. The second intermediate velocity~$\hat{\hat{\bm{u}}}$ is then obtained from equation~\eqref{ProjectionMethodEq1} by projecting~$\hat{\bm{u}}$ onto the space of divergence-free vectors
\begin{align}
-\nabla^2 p^{n+1}= -\frac{\gamma_0 }{\Delta t}\Div{\hat{\bm{u}}} \; ,\label{DualSplitting_PressureStep}\\
\hat{\hat{\bm{u}}} = \hat{\bm{u}} - \frac{\Delta t}{\gamma_0} \Grad{p^{n+1}}\; .\label{DualSplitting_ProjectionStep}
\end{align}
The pressure Poisson equation~\eqref{DualSplitting_PressureStep} is subject to the boundary conditions
\begin{align}
\Grad{p^{n+1}}\cdot\bm{n} &= h_p\left(t_{n+1}\right)\;\; \text{on} \; \boundary{D}\; ,\label{DualSplitting_PressureBC_GammaD}\\
p^{n+1} &= g_p \left(t_{n+1}\right)\;\; \text{on} \; \boundary{N}\; .\label{DualSplitting_PressureBC_GammaN}
\end{align}
The consistent Neumann boundary condition~$h_p$ is derived by multiplying the momentum equation of the incompressible Navier--Stokes equations by the normal vector~$\bm{n}$ and solving for the pressure term~\cite{Karniadakis13,Karniadakis1991}
\begin{equation}
h_p\left(t_{n+1}\right) =
- \left[\frac{\partial \bm{g}_u\left(t_{n+1}\right)}{\partial t} 
+\sum_{i=0}^{J_p-1} \beta_i\left(\Div{\bm{F}_{\mathrm{c}}\left(\bm{u}^{n-i}\right)}
+\nu\nabla\times\boldsymbol{\omega}^{n-i}\right)
 -\bm{f}\left(t_{n+1}\right)\right]\cdot \bm{n}\; ,\label{DualSplitting_PressureBC_GammaD_2}
\end{equation}
where we added the time derivative term and the body force term compared to the formulation in~\cite{Karniadakis1991} in order to extend the formulation to the more general case of time dependent boundary conditions and right-hand side vectors~$\bm{f}\neq \bm{0}$. The time derivative term is calculated using the given boundary values~$\bm{g}_u$ on~$\boundary{D}$. The convective term and the viscous term are formulated explicitly using an extrapolation scheme of order~$J_p$ and the known velocity solution at previous time instants. Note that the viscous term is written in rotational form~$\nabla \times \bm{\omega}$, where~$\bm{\omega}=\nabla \times \bm{u}$ denotes the vorticity. This formulation is obtained by applying the vector identity~$\nabla^2 \bm{u} = \Grad{\left(\Div{\bm{u}}\right)}-\nabla \times \left( \nabla \times \bm{u} \right)=-\nabla \times \left( \nabla \times \bm{u} \right)$ and making use of the incompressibility constraint~$\Div{\bm{u}}=0$. The rotational formulation has first been proposed and analyzed in~\cite{Orszag1986,Karniadakis1991}. It is well known that the rotational formulation effectively reduces boundary divergence errors as compared to the Laplace formulation and is essential in obtaining high-order accuracy in time, see also~\cite{Guermond06,Karniadakis13}. An alternative point of view is provided in~\cite{Leriche2000}, where it is shown that the ellipticity of the Stokes operator is lost with the viscous term written in Laplace formulation.

\subsubsection{Viscous step}
In the final step of the dual splitting scheme the viscous term is considered, leading to the following Helmholtz-like equation
\begin{equation}
\frac{\gamma_0 }{\Delta t} \bm{u}^{n+1}  -  \Div{\bm{F}_{\mathrm{v}}\left(\bm{u}^{n+1}\right)}=
\frac{\gamma_0 }{\Delta t}\hat{\hat{\bm{u}}} \; ,\label{DualSplitting_ViscousStep}
\end{equation} 
where the velocity~$\bm{u}^{n+1}$ has to fulfill the boundary conditions
\begin{align}
\bm{u}^{n+1} &= \bm{g}_{u}\left(t_{n+1}\right)\;\; \text{on} \; \boundary{D}\; ,\label{DualSplitting_ViscousStep_DBC}\\
\bm{F}_{\mathrm{v}} (\bm{u}^{n+1})\cdot \bm{n} &= \bm{h}_u\left(t_{n+1}\right)\;\; \text{on} \; \boundary{N} \; .\label{DualSplitting_ViscousStep_NBC}
\end{align} 

\begin{remark}
For velocity-correction schemes theoretical rates of convergence are available for the case of pure Dirichlet boundary conditions. As shown in~\cite{Guermond2003}, the high-order dual splitting scheme with~$J=2$ and~$J_p=1$ is formally equivalent to the rotational velocity-correction scheme proposed by Guermond and Shen~\cite{Guermond2003} who proved stability and theoretical rates of convergence of order~$\Delta t^2$ in the~$L^2$-norm of the velocity and~$\Delta t^{3/2}$ in the~$L^2$-norm of the pressure for the velocity-correction scheme in rotational form. Numerical investigations in~\cite{Leriche2006} show that the rate of convergence of~$\Delta t^{3/2}$ for the pressure is related to the first order extrapolation~$J_p=1$ in the pressure Neumann boundary condition and that~$J_p=2$ has to be used to obtain optimal rates of convergence (of order~$\Delta t^2$) also for the pressure.~\footnote{While Leriche et al.~\cite{Leriche2006} report even rates of convergence of order~$\Delta t^{5/2}$ for the pressure using the scheme with~$J=J_p=2$, this in not in agreement with the present results shown in Section~\ref{NumericalResults}.} Moreover, an eigenvalue analysis~\cite{Leriche2006} reveals that the high-order dual splitting scheme is only conditionally stable for~$J_p>2$, while it is unconditionally stable for~$J_p\leq 2$ independent of the order~$1\leq J \leq 4$ of the BDF scheme. According to that analysis, among the schemes that provide unconditional stability, the method with~$J=3$ and~$J_p=2$ achieves the highest rates of convergence of order~$\Delta t^3$ for the velocity and~$\Delta t^{5/2}$ for the pressure. In the present paper, the analysis is restricted to the choice~$J_p=J\leq 2$ ensuring that instabilities considered in the present work are not related to the temporal discretization scheme.
\end{remark}

\subsection{Pressure-correction scheme}
This section describes the operator splitting and temporal discretization of pressure-correction schemes and presents different formulations of pressure-correction schemes that are summarized in~\cite{Guermond06}. For more detailed information as well as theoretical aspects the reader is referred to~\cite{Guermond06} and references mentioned therein. With respect to the convective term an implicit formulation is considered in the present work.

\subsubsection{Momentum step}
The momentum equation is solved in the first substep, where the pressure gradient term is either neglected (non-incremental formulation) or an extrapolation of the pressure gradient term based on the pressure solution at previous instants of time (incremental formulation) is used. An intermediate velocity field~$\hat{\bm{u}}$ is calculated in the momentum step by solving the equation
\begin{equation}
\frac{\gamma_0\hat{\bm{u}}-\sum_{i=0}^{J-1}\left(\alpha_i\bm{u}^{n-i}\right)}{\Delta t} 
+ \Div{\bm{F}_{\mathrm{c}}\left(\hat{\bm{u}}\right)}
-\Div{\bm{F}_{\mathrm{v}}\left(\hat{\bm{u}}\right)} = 
- \sum_{i=0}^{J_p-1}\left(\beta_i \Grad{p^{n-i}}\right)
+ \bm{f}\left(t_{n+1}\right)\; ,\label{PressureConvection_MomentumStep_Impl}
\end{equation}
where the boundary conditions for the intermediate velocity field~$\hat{\bm{u}}$ are
\begin{align}
\hat{\bm{u}} &= \bm{g}_{u}\left(t_{n+1}\right)\;\; \text{on} \; \boundary{D} \; ,\label{PressureCorrection_DirichletBC_IntermediateVelocity}\\
\bm{F}_{\mathrm{v}} (\hat{\bm{u}})\cdot \bm{n} &= \bm{h}_u\left(t_{n+1}\right)\;\; \text{on} \; \boundary{N} \; .\label{PressureCorrection_NeumannBC_IntermediateVelocity}
\end{align}
The order of the extrapolation of the pressure gradient term is denoted by~$J_p$ where schemes with~$J_p=0$ are called non-incremental and schemes with~$J_p \geq 1$ incremental pressure-correction schemes~\cite{Guermond06}.

\subsubsection{Pressure step and projection step}
In the second substep, the velocity~$\bm{u}^{n+1}$ and the pressure~$p^{n+1}$ at time~$t_{n+1}$ are obtained as the solution of the following projection method
\begin{align}
\frac{\gamma_0 }{\Delta t}\bm{u}^{n+1} +\Grad{\phi^{n+1}} =\frac{\gamma_0 }{\Delta t}\hat{\bm{u}} \; ,\label{PressureCorrection_ProjectionMethodEq1}\\
\Div{\bm{u}^{n+1}} = 0 \; ,\label{PressureCorrection_ProjectionMethodEq2}\\
\phi^{n+1} = p^{n+1} - \sum_{i=0}^{J_p-1}\left(\beta_i p^{n-i}\right) + \chi \nu \Div{\hat{\bm{u}}}\; .\label{PressureCorrection_Phi}
\end{align}
where the formulation is called standard for~$\chi=0$ and rotational for~$\chi=1$.~\footnote{The standard/rotational terminology is explained in~\cite{Guermond2004} and originates from the following consideration. By inserting equations~\eqref{PressureCorrection_ProjectionMethodEq1} and~\eqref{PressureCorrection_Phi} into equation~\eqref{PressureConvection_MomentumStep_Impl} a Neumann boundary condition can be derived for the pressure. If the divergence correction term is used in equation~\eqref{PressureCorrection_Phi},~$\chi=1$, a boundary condition similar to equation~\eqref{DualSplitting_PressureBC_GammaD_2} is obtained with the viscous term written in rotational formulation, while for the standard form,~$\chi=0$, the viscous term is written in Laplace formulation.} Again, a Poisson equation for the pressure increment can be derived by taking the divergence of equation~\eqref{PressureCorrection_ProjectionMethodEq1} and making use of the divergence-free constraint~\eqref{PressureCorrection_ProjectionMethodEq2}. Subsequently, the pressure solution~$p^{n+1}$ and velocity solution~$\bm{u}^{n+1}$ are given by reformulating equations~\eqref{PressureCorrection_Phi} and by projecting~$\hat{\bm{u}}$ onto the space of divergence-free vectors according to equation~\eqref{PressureCorrection_ProjectionMethodEq1}
\begin{align}
-\nabla^2 \phi^{n+1}= -\frac{\gamma_0 }{\Delta t}\Div{\hat{\bm{u}}} \; ,\label{PressureCorrection_PressurePoissonEquation}\\
p^{n+1} = \phi^{n+1} + \sum_{i=0}^{J_p-1}\left(\beta_i p^{n-i}\right) - \chi \nu \Div{\hat{\bm{u}}}\; ,\label{PressureCorrection_PressureUpdate}\\
\bm{u}^{n+1} = \hat{\bm{u}} - \frac{\Delta t}{\gamma_0} \Grad{\phi^{n+1}}\; .\label{PressureCorrection_Projection}
\end{align}
The pressure Poisson equation~\eqref{PressureCorrection_PressurePoissonEquation} is subject to the boundary conditions
\begin{align}
\Grad{\phi^{n+1}}\cdot\bm{n} =h_{\phi}(t_{n+1}) = 0\;\; \text{on} \; \boundary{D}\; ,\label{PressureCorrection_PressureBC_GammaD}
\\
\phi^{n+1} = g_{\phi}(t_{n+1}) = g_p \left(t_{n+1}\right) -\sum_{i=0}^{J_p-1}\left(\beta_i g_{p}\left(t_{n-i}\right)\right) \;\; \text{on} \; \boundary{N}\; .\label{PressureCorrection_PressureBC_GammaN}
\end{align}
The above boundary conditions~\eqref{PressureCorrection_PressureBC_GammaD} and~\eqref{PressureCorrection_PressureBC_GammaN} are in line with the boundary conditions (10.3) in~\cite{Guermond06}, except that we extend the formulation towards the more general case of inhomogeneous and time-dependent pressure boundary conditions on the Neumann part~$\boundary{N}$ of the boundary.

\begin{remark}
Theoretical rates of convergence of pressure-correction schemes are derived in~\cite{Guermond2004} and summarized in~\cite{Guermond06}. The non-incremental pressure-correction scheme with~$J=1$,~$J_p=0$ in standard form is~$\Delta t$ accurate in the~$L^2$-norm of the velocity and~$\Delta t^{1/2}$ accurate in the~$L^2$-norm of the pressure. The incremental pressure-correction scheme with~$J=2$,~$J_p=1$ is~$\Delta t^2$ accurate in the~$L^2$-norm of the velocity for both the standard formulation and the rotational formulation. While the standard formulation achieves an accuracy of order~$\Delta t$ in the~$L^2$-norm of the pressure, the rotational form is~$\Delta t^{3/2}$ accurate in the~$L^2$-norm of the pressure. As reported in~\cite{Guermond06}, numerical results give evidence that pressure-correction schemes are only conditionally stable for~$J_p \geq 2$. For this reason, we only consider schemes with~$J_p=J-1$ and~$J=1,2$ in the present work which are unconditionally stable.
\end{remark}

\section{Spatial discretization}\label{SpatialDiscretization}
In this section, the discontinuous Galerkin spatial discretization is derived for the different solution strategies discussed in Section~\ref{TemporalDiscretization}. The local Lax--Friedrichs flux is used to discretize the convective term and the symmetric interior penalty Galerkin (SIPG) method to discretize both the viscous term and the Laplace operator in the pressure Poisson equation when dealing with projection methods. The velocity divergence term and pressure gradient term are integrated by parts using a central flux formulation. As demonstrated in what follows, the DG discretization of these terms plays a decisive role in terms of stability and is a central aspect of the present work.

\subsection{Notation}
The physical domain~$\Omega$ is approximated by the computational domain~$\Omega_h \in \mathbb{R}^d$ with boundary~$\Gamma_h = \partial \Omega_h$, where~$\Gamma_h = \hboundary{D} \cup \hboundary{N}$ and~$\hboundary{D}\cap\hboundary{N} = \emptyset$. The computational domain~$\Omega_h = \bigcup_{e=1}^{N_{\text{el}}} \Omega_{e}$ consists of~$N_{\text{el}}$ non-overlapping finite elements~$\Omega_{e}$, where we consider quadrilateral/hexahedral element geometries in this work. The velocity~$\bm{u}(\bm{x},t)$ and pressure~$p(\bm{x},t)$ are approximated by functions ${\bm{u}_h(\bm{x},t)\in\mathcal{V}^{u}_h}$ and~$p_h(\bm{x},t)\in \mathcal{V}^{p}_h$. In the context of discontinuous Galerkin finite element methods, the solution is polynomial inside elements but discontinuous between elements. The spaces of test and trial functions for velocity and pressure are defined as
\begin{align}
\mathcal{V}^{u}_{h} &= \left\lbrace\bm{u}_h\in \left[L_2(\Omega_h)\right]^d\; : \; \bm{u}_h\left(\bm{x}(\boldsymbol{\xi})\right)\vert_{\Omega_{e}}= \tilde{\bm{u}}_h^e(\boldsymbol{\xi})\vert_{\tilde{\Omega}_{e}}\in \mathcal{V}^{u}_{h,e}=[\mathcal{P}_{k_u}(\tilde{\Omega}_{e})]^d\; ,\;\; \forall e=1,\ldots,N_{\text{el}} \right\rbrace\;\; ,\\
\mathcal{V}^{p}_{h} &= \left\lbrace p_h\in L_2(\Omega_h)\; : \; p_h\left(\bm{x}(\boldsymbol{\xi})\right)\vert_{\Omega_{e}} = \tilde{p}_h^e(\boldsymbol{\xi})\vert_{\tilde{\Omega}_{e}}\in \mathcal{V}^{p}_{h,e}=\mathcal{P}_{k_p}(\tilde{\Omega}_{e})\; ,\;\; \forall e=1,\ldots,N_{\text{el}} \right\rbrace\; ,
\end{align}
respectively, where~$\mathcal{P}_{k}(\tilde{\Omega}_{e})$ denotes the space of polynomials of tensor degree~$\leq k$ on the reference element~$\tilde{\Omega}_e=[0,1]^d$ with reference coordinates~$\boldsymbol{\xi}=(\xi_1,...,\xi_d)^T$. In the above equations,~$\bm{x}(\boldsymbol{\xi}) : \tilde{\Omega}_e \rightarrow \Omega_e$ denotes the mapping from reference space to physical space. We use a nodal approach so that the approximate solutions of velocity and pressure on element~$e$ can be written as
\begin{align}
\tilde{\bm{u}}_h^e(\boldsymbol{\xi},t) = \sum_{i_1,...,i_d=0}^{k_u} N_{i_1...i_d}^{k_u}(\boldsymbol{\xi})\bm{u}_{i_1...i_d}^e(t)\;\; , \;\; \tilde{p}_h^e(\boldsymbol{\xi},t) = \sum_{i_1,...,i_d=0}^{k_p} N_{i_1...i_d}^{k_p}(\boldsymbol{\xi})p_{i_1...i_d}^e(t)\; ,
\end{align}
where~$ \bm{u}_{i_1...i_d}^e$ and~$p_{i_1...i_d}^e$ denote the nodal degrees of freedom of the velocity and pressure solution on element~$e$, respectively. The multidimensional shape functions~$N_{i_1...i_d}^{k}$ are given as the tensor product of one-dimensional shape functions,~$N_{i_1...i_d}^{k}(\boldsymbol{\xi})=\prod_{n=1}^{d} l_{i_n}^k(\xi_n)$, where~$l_i^k(\xi)$ are the Lagrange polynomials of degree~$k$ based on the Legendre--Gauss--Lobatto nodes. For geometries with curved boundaries a polynomial mapping~$\bm{x}(\boldsymbol{\xi})$ of degree~$k_u$ is used for a high order accurate interpolation of the geometry.
In this work, both equal-order polynomials for velocity and pressure ($k_p=k_u$) and mixed-order polynomials ($k_p=k_u-1$) are analyzed. 

The interface of two adjacent elements~$\Omega_{e^-}$ and~$\Omega_{e^+}$ is denoted by~$f_{e^-/e^+}=\partial \Omega_{e^-} \cap \partial \Omega_{e^+}$ where the outward pointing normal vectors on~$f_{e^-/e^+}$ are denoted by~$\bm{n}^{-}$ for~$\Omega_{e^-}$ and~$\bm{n}^{+}$ for~$\Omega_{e^+}$. Furthermore, let~$u^{-}_h$ and~$u^{+}_h$ denote the solution~$u_h$ on~$f_{e^-,e^+}$ evaluated from the interior of element~$e^-$ and element~$e^+$, respectively. Following~\cite{bassi2006artificial,bassi2007implicit}, the average operator~$\avg{\cdot}$ and jump operator~$\jump{\cdot}$ are defined as~$\avg{u} = (u^- + u^+)/2$ and~$ \jump{u} = u^- \otimes \bm{n}^- + u^+ \otimes \bm{n}^+$, respectively. We note that both operators can be applied to a scalar, vectorial or tensorial quantity~$u$ and that defining numerical fluxes in terms of these operators guarantees conservativity of the numerical flux.

We use an element-by-element formulation when deriving the weak formulation, i.e., we use a notation where volume integrals are performed over the current element~$\Omega_e$ and face integrals over the boundary~$\partial \Omega_e$ of element~$e$. Integrals over~$\Omega_e$ and~$\partial\Omega_e$ are abbreviated by using the shorthand notation~$\intele{v}{u} = \int_{\Omega_e} v \odot u \; \mathrm{d}\Omega$ and~$\inteleface{v}{u} = \int_{\partial \Omega_e} v \odot u \; \mathrm{d} \Gamma$, where the operator~$\odot$ symbolizes inner products, i.e.,~$v u$ for rank-0 tensors,~$\bm{v}\cdot\bm{u} = v_i u_i$ for rank-1 tensors, and~$\bm{v} : \bm{u} = v_{ij} u_{ij}$ for rank-2 tensors.\\
Moreover, we introduce the convention that interior information on the current element~$\Omega_e$ is denoted by the superscript~$(\cdot)^-$ and exterior information from neighboring elements by the superscript~$(\cdot)^+$. Accordingly, the normal vector~$\bm{n}$ of the current element~$\Omega_e$ is equal to~$\bm{n}^-$, while~$\bm{n}^+=-\bm{n}^-=-\bm{n}$.

\subsection{General form of spatial discretization of coupled solution approach}
In the following, the weak discontinuous Galerkin formulation of the incompressible Navier--Stokes equations using a coupled solution approach described in Section~\ref{TemporalDiscretization_CoupledSolution} is derived in two steps: 

\begin{enumerate}[label=(\roman*), ref=(\roman*)]
\item\label{DGDiscretizationStep1} 
 by requiring the time discrete residuals of the momentum equation~\eqref{TemporalDiscretization_Coupled_Momentum} and the continuity equation~\eqref{TemporalDiscretization_Coupled_Continuity} to be orthogonal to all test functions~$\bm{v}_h\in \mathcal{V}^{u}_h$ and~$q_h \in \mathcal{V}^{p}_h$, respectively, and 
\item\label{DGDiscretizationStep2}  by performing integration by parts including the definition of numerical fluxes, which represents the core of discontinuous Galerkin methods. 
\end{enumerate}
Multiplying the residual of the momentum equation by test functions~$\bm{v}_h$ and the residual of the continuity equation by test functions~$q_h$ as well as integration over~$\Omega_h$ yields the following set of equations
\begin{align}
\begin{split}
\intele{\bm{v}_h}{\frac{\gamma_0 \bm{u}^{n+1}_h-\sum_{i=0}^{J-1}\left(\alpha_i\bm{u}^{n-i}_h\right)}{\Delta t}}
+ \intele{\bm{v}_h}{\nabla \cdot \bm{F}_{\mathrm{c}}(\bm{u}^{n+1}_h)}  & \\
- \intele{\bm{v}_h}{\nabla \cdot \bm{F}_{\mathrm{v}} (\bm{u}^{n+1}_h)}
+ \intele{\bm{v}_h}{\nabla p^{n+1}_h} - \intele{\bm{v_h}}{\bm{f}(t_{n+1})}& = 0 \;\;
\forall \bm{v}_h \in \mathcal{V}^{u}_{h,e} \; ,
\end{split}\\
\intele{q_h}{-\nabla \cdot \bm{u}^{n+1}_h} & = 0  \;\; \forall q_h \in \mathcal{V}^{p}_{h,e}\label{WeakForm_CoupledSolution_Continuity_Step1} \; ,
\end{align}
for all elements~$e=1,...N_{\text{el}}$. Terms involving spatial derivative operators are then integrated by parts. In this step, Gauss' divergence theorem is applied in order to transform volume integrals into surface integrals. Subsequently, physical fluxes are replaced by numerical fluxes in order to enforce continuity in a weak sense. Numerical fluxes are defined as a function of the approximate solution on both elements adjacent to an interior face and as a function of the interior solution and prescribed boundary data on boundary faces. By the example of the convective term, this second step can be generically written as
\begin{align}
\intele{\bm{v}_h}{\nabla \cdot \bm{F}_{\mathrm{c}}(\bm{u}_h)}  \hspace{0,5cm}\rightarrow \hspace{0,5cm} c^e_h\left(\bm{v}_h,\bm{u}_h,\bm{g}_u(t_{n+1})\right) \; .
\end{align}
As a result, we obtain the following weak discontinuous Galerkin formulation of the fully discrete, incompressible Navier--Stokes equations: Find~$\bm{u}^{n+1}_h\in\mathcal{V}^u_h$,~$p^{n+1}_h\in \mathcal{V}^{p}_h$ such that 
\begin{align}
\begin{split}
m^e_{h,u}\left(\bm{v}_h,\frac{\gamma_0 \bm{u}^{n+1}_h-\sum_{i=0}^{J-1}\left(\alpha_i\bm{u}^{n-i}_h\right)}{\Delta t} \right)
+ c^e_h\left(\bm{v}_h,\bm{u}^{n+1}_h,\bm{g}_u(t_{n+1})\right)& \\
+ v^e_h\left(\bm{v}_h,\bm{u}^{n+1}_h,\bm{g}_u(t_{n+1}),\bm{h}_u(t_{n+1})\right)
+ g^e_h\left(\bm{v}_h,p^{n+1}_h,g_p(t_{n+1})\right) -  b^e_h\left(\bm{v}_h,\bm{f}(t_{n+1})\right) &= 0 \; ,
\end{split} \label{WeakForm_CoupledSolution_Momentum}\\
-d^e_h(q_h,\bm{u}^{n+1}_h,\bm{g}_u(t_{n+1}))& = 0  \; ,\label{WeakForm_CoupledSolution_Continuity}
\end{align}
for all~$(\bm{v}_h, q_h) \in \mathcal{V}^{u}_{h,e} \times \mathcal{V}^{p}_{h,e}$ and for all elements~$e=1,...,N_{\text{el}}$. The minus sign is inserted in equation~\eqref{WeakForm_CoupledSolution_Continuity_Step1} and equation~\eqref{WeakForm_CoupledSolution_Continuity} to ensure that the matrix representation of the (linearized) system of equations corresponding to the weak formulation~\eqref{WeakForm_CoupledSolution_Momentum} and~\eqref{WeakForm_CoupledSolution_Continuity} is symmetric with respect to the pressure gradient term and the velocity divergence term.

The velocity mass matrix operator is given in elementwise notation as~$m^{e}_{h,u}\left(\bm{v}_h,\bm{u}_h\right) = \intele{\bm{v}_h}{\bm{u}_h}$ and the body force operator as~$ b^{e}_{h}\left(\bm{v}_h,\bm{f}\right) = \intele{\bm{v}_h}{\bm{f}}$. Since these terms do not contain spatial derivative operators, there is no need to perform step~\ref{DGDiscretizationStep2} described above. Consequently, the mass matrix operator is block-diagonal and can be inverted locally (element-by-element). A detailed description of the convective term~$c^e_h$, viscous term~$v^e_h$, pressure gradient term~$g^e_h$, and velocity divergence term~$d^e_h$ is given below. Apart from these operators, we introduce in Section~\ref{DGBasicOperators} the DG formulation of the negative Laplace operator~$l_h^{e}$ required when discretizing the projection type solution methods in space.
 
 \begin{table}[!h]
\caption{Choice of exterior values~$\left(\cdot\right)^+$ on domain boundaries as a function of interior values~$\left(\cdot\right)^-$ and prescribed boundary data for velocity and pressure in order to weakly impose boundary conditions using a mirror principle.}\label{BCsWeak}
\renewcommand{\arraystretch}{1.1}
\begin{center}
\begin{tabular}{ccc}
\hline  
 & $\hboundary{D}$ & $\hboundary{N}$\\ 
\hline
\multirow{2}{*}{velocity} & $\bm{u}_h^{+} = -\bm{u}_h^{-} + 2 \bm{g}_{u}$ & $\bm{u}_h^{+} = \bm{u}_h^{-}$\\
 & $\Grad{\bm{u}_{h}^{+}}\cdot \bm{n} = \Grad{\bm{u}_{h}^{-}}\cdot\bm{n}$ & $\Grad{\bm{u}_{h}^{+}}\cdot \bm{n} =  -\Grad{\bm{u}_{h}^{-}}\cdot\bm{n}+\frac{2\bm{h}_u}{\nu}$\\

\multirow{2}{*}{pressure} & $p^+_h = p^-_h$ & $p^+_h = - p^-_h + 2 g_p$\\
& $\Grad{p_{h}^{+}}\cdot \bm{n} = -\Grad{p_{h}^{-}}\cdot\bm{n} + 2 h_p$ & $\Grad{p_{h}^{+}}\cdot \bm{n} = \Grad{p_{h}^{-}}\cdot\bm{n}$\\
\hline
\end{tabular}
\end{center}
\renewcommand{\arraystretch}{1}
\end{table}

\subsection{DG formulation of basic operators}\label{DGBasicOperators}
\subsubsection{Convective term}
In order to derive the discontinuous Galerkin formulation of the convective term we perform step~\ref{DGDiscretizationStep2}. Integration by parts of the convective term~$\intele{\bm{v}_h}{\nabla \cdot \bm{F}_{\mathrm{c}}(\bm{u}_h)} $ and replacing the physical flux~$\bm{F}_{\mathrm{c}}(\bm{u}_h)$ by the numerical flux~$\bm{F}^{*}_{\mathrm{c}}(\bm{u}_h)$ yields
\begin{equation}
c^e_h\left(\bm{v}_h,\bm{u}_h,\bm{g}_u\right) = -\intele{\Grad{\bm{v}_h}}{\bm{F}_{\mathrm{c}}(\bm{u}_h)} + 
\inteleface{\bm{v}_h}{\bm{F}^{*}_{\mathrm{c}}(\bm{u}_h)\cdot \bm{n}} \; .\label{WeakFormConvectiveTerm}
\end{equation}
The local Lax--Friedrichs flux is defined as~\cite{Hesthaven07,Shahbazi07,Klein13}
\begin{equation}
\bm{F}^{*}_{\mathrm{c}} (\bm{u}_h) = \avg{\bm{F}_{\mathrm{c}}(\bm{u}_h)} + \frac{\Lambda}{2}\jump{\bm{u}_h} \; ,\label{LaxFriedrichsFlux}
\end{equation}
where~$\Lambda = \max \left(\lambda^{-},\lambda^{+}\right)$.
The maximum eigenvalue~$\lambda$  (in terms of absolute values) of the flux Jacobian is
\begin{align}
\lambda^{\pm} &=  \max_{i} \left\vert \lambda_i\left(\left.\frac{\partial \bm{F}(\bm{u})\cdot\bm{n}}{\partial\bm{u}}\right\vert_{\bm{u}_h^{\pm}}\right) \right\vert = 2 \vert \bm{u}_h^{\pm} \cdot \bm{n}\vert \; .
\end{align}
In the above equation,~$\bm{u}_h^{\pm}$ is the local velocity evaluated in each quadrature point. Note that mean values of the velocity are used in~\cite{Hesthaven07,Shahbazi07,Klein13}.
Boundary conditions are imposed by calculating exterior values~$\bm{u}_h^+$ on~$\Gamma_h$ as defined in Table~\ref{BCsWeak}. In order to highlight that the convective term depends on the prescribed boundary data~$\bm{g}_{u}$ on Dirichlet boundaries~$\hboundary{D}$, we use the notation~$c^e_h\left(\bm{v}_h,\bm{u}_h,\bm{g}_u\right)$.

\subsubsection{Velocity divergence term}
The DG formulation of the velocity divergence term is derived by performing integration by parts of~$\intele{q_h}{\Div{\bm{u}_h}}$ and 
replacing the physical flux~$\bm{u}_h$ by the numerical flux~$\bm{u}^*_h$ to obtain
\begin{equation}
d^e_h\left(q_h,\bm{u}_h,\bm{g}_u\right) = -\intele{\Grad{q_h}}{\bm{u}_h}+\inteleface{q_h}{\bm{u}^*_h\cdot\bm{n}}\; .\label{WeakFormVelocityDivergence}
\end{equation}
As numerical flux function we use the central flux
\begin{equation}
\bm{u}^*_h = \avg{\bm{u}_h}\; .\label{NumericalFluxDivergenceTerm}
\end{equation}
Inserting equation~\eqref{NumericalFluxDivergenceTerm} into equation~\eqref{WeakFormVelocityDivergence} and calculating exterior values~$\bm{u}_h^+$ according to Table~\ref{BCsWeak} results in the following DG formulation of the divergence operator
\begin{align}
d^e_h\left(q_h,\bm{u}_h,\bm{g}_u\right) = 
-\intele{\Grad{q_h}}{\bm{u}_h}
+\intelefaceInterior{q_h}{\avg{\bm{u}_h}\cdot\bm{n}}
+\intelefaceNeumann{q_h}{\bm{u}_h\cdot\bm{n}}
+\intelefaceDirichlet{q_h}{\bm{g}_{u}\cdot\bm{n}} \; .
\label{WeakFormVelocityDivergenceBoundaryIntegrals}
\end{align}

\begin{remark}
As a reference formulation we consider a modified formulation of the velocity divergence term used in~\cite{Hesthaven07} in the context of the high-order dual splitting scheme
\begin{equation}
d^e_{h,\text{ref}}\left(q_h,\bm{u}_h\right) = \intele{q_h}{\Div{\bm{u}_h}}\; .\label{VelocityDivergence_ReferenceFormulation}
\end{equation}
This formulation does not perform integration by parts as described in step~\ref{DGDiscretizationStep2} above. Accordingly, this formulation does not depend on boundary conditions prescribed for the velocity.
\end{remark}

\subsubsection{Pressure gradient term}
The procedure detailed above is applied to obtain the DG formulation~$g^e_h\left(\bm{v}_h,p_h\right)$ of the pressure gradient term
\begin{equation}
g^e_h\left(\bm{v}_h,p_h,g_p\right) = -\intele{\Div{\bm{v}_h}}{p_h}+\inteleface{\bm{v}_h}{p^*_h\bm{n}}\; .
\label{WeakFormPressureGradient}
\end{equation}
As for the divergence term, the numerical flux~$p^*_h$ is defined as the central flux
\begin{equation}
p^*_h = \avg{p_h}\; ,\label{NumericalFluxPressureGradient}
\end{equation}
where exterior values~$p_h^+$ on domain boundaries~$\Gamma_h$ are calculated as listed in Table~\ref{BCsWeak}. Inserting equation~\eqref{NumericalFluxPressureGradient} along with the respective boundary conditions into equation~\eqref{WeakFormPressureGradient} yields
\begin{align}
g^e_h\left(\bm{v}_h,p_h,g_p\right) =
-\intele{\Div{\bm{v}_h}}{p_h}
+\intelefaceInterior{\bm{v}_h}{\avg{p_h}\bm{n}}
+\intelefaceDirichlet{\bm{v}_h}{p_h\bm{n}}
+\intelefaceNeumann{\bm{v}_h}{g_p\bm{n}}\; .
\label{WeakFormPressureGradientBoundaryIntegrals}
\end{align}

\begin{remark}
As a reference formulation we consider a modified formulation of the pressure gradient term used in~\cite{Hesthaven07} in the context of the high-order dual splitting scheme
\begin{equation}
g^e_{h,\text{ref}}\left(\bm{v}_h,p_h\right) = \intele{\bm{v}_h}{\Grad{p_h}}\; .\label{PressureGradient_ReferenceFormulation}
\end{equation}
This formulation does not perform integration by parts as described in~\ref{DGDiscretizationStep2} above. Accordingly, this formulation does not depend on boundary conditions prescribed for the pressure.
\end{remark}

\subsubsection{Negative Laplace operator}\label{LaplaceOperator}
In case of the velocity-correction scheme and pressure-correction scheme the pressure solution is obtained by solving a Poisson equation. To derive the DG formulation of the negative Laplace operator we consider the following Poisson-type model problem
\begin{align}
- \nabla^2 p = f\;\; \text{in} \; \Omega\; ,\label{PoissonProblem}
\end{align}
subject to boundary conditions
\begin{align}
p &= g_p\;\; \text{on} \; \boundary{D}_{\mathrm{PPE}}=\boundary{N}\; ,\\
\Grad{p} \cdot \bm{n} &= h_p\;\; \text{on} \; \boundary{N}_{\mathrm{PPE}}=\boundary{D} \; .
\end{align}
Note that the role of Dirichlet and Neumann boundaries is interchanged for the pressure Poisson equation as compared to the definition of Dirichlet and Neumann boundaries for the incompressible Navier--Stokes equations. In the DG context, the weak formulation for this problem including second derivatives is typically derived by rewriting the original equation as a system of first order equations and applying the above procedure separately for each of the first order equations~\cite{Hesthaven07}. These equations can then be recombined to obtain the primal formulation for the unknown solution~$p_h$: Find~$p_h\in \mathcal{V}_h^p$ such that
\begin{align}
l_{h}^{e}\left(q_h,p_h,g_p,h_p\right) = \intele{q_h}{f} \;\; \forall q_h \in \mathcal{V}^{p}_{h,e}\; ,\label{DGWeakFormPoissonEquation}
\end{align}
and for all elements~$e=1,...,N_{\text{el}}$, where~$l_h^e$ is given as
\begin{align}
l_{h}^{e}\left(q_h,p_h,g_p,h_p\right) = \intele{\Grad{q_h}}{\Grad{p_h}}-\inteleface{\Grad{q_h}}{\left(p_h-p_h^*\right)\bm{n}}- \inteleface{q_h}{\bm{\sigma}_h^*\cdot\bm{n}} .\label{DGLaplace_PrimalFormVariant1}
\end{align}
We consider the symmetric interior penalty Galerkin (SIPG) method for which the numerical fluxes are defined as~\cite{arnold2000discontinuous,arnold2002unified}
\begin{align}
p_h^* &= \avg{p_h}\; ,\label{DGLaplace_InteriorPenaltyFlux1}\\
\bm{\sigma}_h^* &= \avg{\Grad{p}_h}-\tau \jump{p_h}\; .\label{DGLaplace_InteriorPenaltyFlux2}
\end{align}
The penalty parameter of the SIPG method is denoted by~$\tau$ and has to be large enough to ensure coercivity of the bilinear form. Essentially, the penalty parameter depends on the polynomial degree~$k$ and a characteristic element length~$h$. An explicit expression for the penalty parameter of the SIPG method is derived in~\cite{Shahbazi05} for triangular/tetrahedral elements and in~\cite{Hillewaert13} for other element geometries. For quadrilateral/hexahedral elements the penalty parameter~$\tau_e$ associated to element~$e$ is defined as~\cite{Hillewaert13}
\begin{align}
\tau_e = (k+1)^2 \frac{A\left(\partial \Omega_e \setminus \Gamma_h\right)/2 + A\left(\partial \Omega_e \cap \Gamma_h\right)}{V\left(\Omega_e\right)}\; ,\label{TauIP_Element}
\end{align}
with the element volume~$V\left(\Omega_e\right) = \int_{\Omega_e}\mathrm{d}\Omega$ and the surface area~$A(f) = \int_{f\subset\partial\Omega_e}\mathrm{d}\Gamma$.
On interior faces, the penalty parameter~$\tau$ is obtained by taking the maximum value of both elements adjacent to face~$f$
\begin{align}
\tau = 
\begin{cases}
\max\left(\tau_{e^-},\tau_{e^+}\right) & \text{if face } f \subseteq \partial \Omega_e \setminus \Gamma_h\; ,\\
\tau_e & \text{if face } f \subseteq \partial \Omega_e \cap \Gamma_h\; .
\end{cases}\label{TauIP}
\end{align}
Boundary conditions are imposed in the weak formulation by defining exterior values for both the pressure~$p_h^+$ and the pressure gradient in normal direction~$\Grad{p_h}^+\cdot\bm{n}$. By inserting the numerical fluxes~\eqref{DGLaplace_InteriorPenaltyFlux1} and~\eqref{DGLaplace_InteriorPenaltyFlux2} as well as the boundary conditions specified in Table~\ref{BCsWeak} into equation~\eqref{DGLaplace_PrimalFormVariant1} the weak formulation of the Laplace operator~$l^{e}_{h}$ can be seperated into a homogeneous part~$l^{e}_{h,\text{hom}}$ and an inhomogeneous part~$l^{e}_{h,\text{inhom}}$
according to
\begin{align}
l_{h}^{e}\left(q_h,p_h,g_p,h_p\right)
 = l^{e}_{h,\text{hom}}(q_h,p_h) +l_{h,\text{inhom}}^{e}(q_h,g_p,h_p)\; ,\label{NegativeLaplace_HomogeneousInhomogeneous}
\end{align}
where~$l^{e}_{h,\text{inhom}}$ contains the inhomgoneous parts of boundary face integrals and is shifted to the right-hand side of equation~\eqref{DGWeakFormPoissonEquation} when solving the linear system of equations.

\subsubsection{Viscous term}
The viscous operator represents a generalization of the Laplace operator to vectorial quantities  with the viscosity~$\nu$ as a scaling factor. In analogy to the Laplace operator considered in Section~\ref{LaplaceOperator}, the primal formulation of the viscous term is given as
\begin{align}
v_{h}^{e}\left(\bm{v}_h,\bm{u}_h,\bm{g}_u,\bm{h}_u\right) =
\intele{\Grad{\bm{v}_h}}{\nu\Grad{\bm{u}_h}}
-\inteleface{\Grad{\bm{v}_h}}{\nu \left(\bm{u}_h
-\bm{u}^*_h\right)\otimes \bm{n}}-\inteleface{\bm{v}_h}{\bm{F}^*_{\mathrm{v},h}\cdot\bm{n}} \; .\label{DGViscous_LaplaceForm_PrimalFormVariant1}
\end{align}
For the symmetric interior penalty Galerkin (SIPG) method the numerical fluxes are defined as
\begin{align}
\bm{u}^*_h &=\avg{\bm{u}_h}\; ,\label{DGViscous_LaplaceForm_NumericalFlux1}\\
\bm{F}^*_{\mathrm{v},h} &= \nu \avg{\Grad{\bm{u}_h}} -\nu\tau\jump{\bm{u}_h}\; ,\label{DGViscous_LaplaceForm_NumericalFlux2}
\end{align}
where the interior penalty parameter is defined as in equation~\eqref{TauIP} and equation~\eqref{TauIP_Element}. Again, boundary conditions are incorporated into the formulation by defining exterior values for the velocity~$\bm{u}_h^+$ and the velocity gradient in normal direction~$\Grad{\bm{u}_h^+}\cdot\bm{n}$. Inserting the numerical fluxes~\eqref{DGViscous_LaplaceForm_NumericalFlux1} and~\eqref{DGViscous_LaplaceForm_NumericalFlux2} into equation~\eqref{DGViscous_LaplaceForm_PrimalFormVariant1} and imposing boundary conditions according to Table~\ref{BCsWeak}, the weak formulation of the viscous operator~$v^{e}_{h}$ can be written as the sum of a homogeneous part~$v^{e}_{h,\text{hom}}$ and an inhomogeneous part~$v^{e}_{h,\text{inhom}}$
\begin{align}
v_{h}^{e}\left(\bm{v}_h,\bm{u}_h,\bm{g}_u,\bm{h}_u\right) = v^{e}_{h,\text{hom}}(\bm{v}_h,\bm{u}_h) + v_{h,\text{inhom}}^{e}(\bm{v}_h,\bm{g}_u,\bm{h}_u)\; .\label{DGViscous_LaplaceForm_HomogeneousInhomogeneous}
\end{align}
\begin{remark}
In the weak formulation of the momentum equation of the coupled solution approach, equation~\eqref{WeakForm_CoupledSolution_Momentum}, and for the derivation of the weak formulation of the viscous term and the pressure gradient term, a splitting of the Neumann boundary condition into a viscous part and a pressure part according to equations~\eqref{NeumannBC_Viscous} and~\eqref{NeumannBC_Pressure} is used, although equation~\eqref{NeumannBC_Coupled} defines the Neumann boundary condition to be prescribed in case of the coupled solution approach. However, the inhomogeneous boundary face integrals of the viscous term
\begin{align}
v_{h,\text{inhom}}^{e}(\bm{v}_h,\bm{g}_u,\bm{h}_u) = 
 \intelefaceDirichlet{\Grad{\bm{v}_h}}{\nu\, \bm{g}_{u} \otimes\bm{n}}
- \intelefaceNeumann{\bm{v}_h}{\bm{h}_u}
- \intelefaceDirichlet{\bm{v}_h}{2\nu\tau\bm{g}_{u}}\; ,
\label{DGViscous_LaplaceForm_Inhomogeneous}
\end{align}
and the pressure gradient term according to equation~\eqref{WeakFormPressureGradientBoundaryIntegrals} are added in equation~\eqref{WeakForm_CoupledSolution_Momentum} so that any decomposition~$\bm{h}=\bm{h}_u-g_p\bm{n}$ ensures a correct imposition of the Neumann boundary condition~\eqref{NeumannBC_Coupled}
\begin{align}
-\intelefaceNeumann{\bm{v}_h}{\bm{h}_u} +\intelefaceNeumann{\bm{v}_h}{g_p\bm{n}}&=-\intelefaceNeumann{\bm{v}_h}{\bm{h}_u-g_p \bm{n}}
=-\intelefaceNeumann{\bm{v}_h}{\bm{h}}\; .
\end{align}
Without loss of generality one can use~$\bm{h}_u=\bm{h}$ and~$g_p=0$ for the coupled solution approach. A decoupled treatment of the Neumann boundary condition according to equations~\eqref{NeumannBC_Viscous} and~\eqref{NeumannBC_Pressure} is, however, necessary for the projection methods discussed below.
\end{remark}

\subsection{High-order dual splitting scheme}\label{WeakFormDualSplitting}
In this section, we briefly summarize the discontinuous Galerkin discretization of the dual splitting scheme using the DG formulation of basic operators derived in Section~\ref{DGBasicOperators}.
\subsubsection{Convective step}
The weak DG formulation of the convective step~\eqref{DualSplitting_ConvectiveStep} is given as follows: Find~$\hat{\bm{u}}_h\in\mathcal{V}^u_h$ such that
\begin{align}
m^{e}_{h,u}\left(\bm{v}_h,\frac{\gamma_0 \hat{\bm{u}}_h-\sum_{i=0}^{J-1}\left(\alpha_i\bm{u}^{n-i}_h\right)}{\Delta t} \right)
= 
- \sum_{i=0}^{J-1} \left(\beta_i c^e_h\left(\bm{v}_h,\bm{u}^{n-i}_h,\bm{g}_u(t_{n-i})\right)\right)
+ b^{e}_{h}\left(\bm{v}_h,\bm{f}(t_{n+1})\right)\; ,
\label{DualSplitting_ConvectiveStep_WeakForm}
\end{align}
for all~$\bm{v}_h \in \mathcal{V}^{u}_{h,e}$ and for all elements~$e=1,...,N_{\text{el}}$.

\subsubsection{Pressure step}
The weak formulation of the pressure Poisson equation~\eqref{DualSplitting_PressureStep} reads: Find~$p_h^{n+1}\in\mathcal{V}^p_h$ such that
\begin{align}
l_{h,\text{hom}}^{e}\left(q_h,p_h^{n+1}\right) = - \frac{\gamma_0}{\Delta t} d_{h}^{e}\left(q_h,\hat{\bm{u}}_h,\bm{g}_{\hat{u}}(t_{n+1})\right)
 - l_{h,\text{inhom}}^{e}\left(q_h,g_{p}(t_{n+1}),h_{p}(t_{n+1})\right)\; ,
\label{DualSplitting_Pressure_WeakForm}
\end{align}
for all~$q_h \in \mathcal{V}^{p}_{h,e}$ and for all elements~$e=1,...,N_{\text{el}}$. 
In order to evaluate the discrete divergence operator on the right-hand side of the pressure Poisson equation, a boundary condition~$\bm{g}_{\hat{u}}(t_{n+1})$ has to be specified on~$\Gamma^{\mathrm{D}}_h$ for the intermediate velocity field~$\hat{\bm{u}}_h$ according to equation~\eqref{WeakFormVelocityDivergenceBoundaryIntegrals}. Note that applying~$\bm{g}_{u}(t_{n+1})$ as boundary condition is inconsistent and, hence, does not yield optimal rates of convergence with respect to the temporal discretization. To obtain a consistent boundary condition, we derive the boundary condition~$\bm{g}_{\hat{u}}(t_{n+1})$ by solving equation~\eqref{DualSplitting_ConvectiveStep} for the intermediate velocity
\begin{equation}
\bm{g}_{\hat{u}}\left(t_{n+1}\right) = \sum_{i=0}^{J-1}\left(\frac{\alpha_i}{\gamma_0}\bm{g}_{u}(t_{n-i})\right)-\frac{\Delta t}{\gamma_0} \sum_{i=0}^{J-1}\left(\beta_i \Div{\bm{F}_{\mathrm{c}}\left(\bm{u}_h^{n-i}\right)}\right)
+ \frac{\Delta t}{\gamma_0} \bm{f}\left(t_{n+1}\right)\; ,\label{DualSplitting_DBC_IntermediateVelocity}
\end{equation} 
where the fact has been used that~$\bm{u}$ satisfies the boundary condition~$\bm{g}_u$ on~$\Gamma^{\mathrm{D}}$ according to equation~\eqref{DualSplitting_ViscousStep_DBC}. We note that this boundary condition is essential in order to obtain a method that is both stable in the limit of small time step sizes and that ensures higher order accuracy of the temporal discretization.

The boundary values~$g_{p}$ and~$h_{p}$ in the above pressure Poisson equation are defined in equation~\eqref{DualSplitting_PressureBC_GammaN} and equation~\eqref{DualSplitting_PressureBC_GammaD}, respectively. To evaluate~$\bm{g}_{\hat{u}}(t_{n+1})$ and~$h_{p}(t_{n+1})$ on the right-hand side of equation~\eqref{DualSplitting_Pressure_WeakForm} according to the boundary conditions~\eqref{DualSplitting_DBC_IntermediateVelocity} and~\eqref{DualSplitting_PressureBC_GammaD_2}, the convective term and the viscous term have to be calculated on~$\partial \Omega_e$ as a function of the approximate velocity solution~$\bm{u}_h$ on element~$e$. In the discrete case, the divergence of the convective term is calculated as
\begin{align}
\Div{\bm{F}_{\mathrm{c}}\left(\bm{u}_h\right)} = \bm{u}_h \left(\Div{\bm{u}_h}\right) + \left(\Grad{\bm{u}_h}\right) \cdot \bm{u}_h\; .
\end{align}
The viscous term in equation~\eqref{DualSplitting_PressureBC_GammaD_2} involves second derivatives and is calculated in two steps so that the computation of second derivatives is replaced by a sequence of first derivatives. The vorticity~$\bm{\omega}_h \in \mathcal{V}^u_h$ in equation~\eqref{DualSplitting_PressureBC_GammaD_2} is calculated by a local~$L^2$-projection
\begin{align}
\intele{\bm{v}_h}{\bm{\omega}_h} = \intele{\bm{v}_h}{\nabla \times \bm{u}_h}\; .
\end{align}
The viscous term is then evaluated by calculating the curl of the vorticity~$\bm{\omega}_h$ on the respective boundary.

\subsubsection{Projection step}
In elementwise notation, the weak form of the projection step~\eqref{DualSplitting_ProjectionStep} is to find~$\hat{\hat{\bm{u}}}_h\in\mathcal{V}^u_h$ such that
\begin{align}
m_{h,u}^{e}(\bm{v}_h,\hat{\hat{\bm{u}}}_h) = m_{h,u}^{e}\left(\bm{v}_h,\hat{\bm{u}}_h\right)-\frac{\Delta t}{\gamma_0}g_h^{e}\left(\bm{v}_h,p_h^{n+1},g_{p}\left(t_{n+1}\right)\right)\; ,\label{DualSplitting_Projection_WeakForm}
\end{align}
for all~$\bm{v}_h \in \mathcal{V}^{u}_{h,e}$ and for all elements~$e=1,...,N_{\text{el}}$. Note that the evaluation of the discrete pressure gradient term as defined in equation~\eqref{WeakFormPressureGradientBoundaryIntegrals} requires a pressure Dirichlet boundary condition to be prescribed on~$\Gamma_h^{\mathrm{N}}$, where the Dirichlet boundary value~$g_p$ in equation~\eqref{DualSplitting_Projection_WeakForm} is the value prescribed in equation~\eqref{DualSplitting_PressureBC_GammaN}. 
\subsubsection{Viscous step}
The viscous step completes time step~$n$ of the dual splitting scheme by solving a Helmholtz-like equation for the velcity~$\bm{u}_h^{n+1}$. As for the pressure Poisson equation, inhomogeneous boundary face integrals are shifted to the right-hand side to obtain the following weak formulation of equation~\eqref{DualSplitting_ViscousStep}: Find~$\bm{u}_h^{n+1}\in\mathcal{V}^u_h$ such that
\begin{align}
m^{e}_{h,u}\left(\bm{v}_h,\frac{\gamma_0}{\Delta t} \bm{u}_h^{n+1} \right) 
+ v^{e}_{h,\text{hom}}\left(\bm{v}_h,\bm{u}_h^{n+1}\right)
= 
m^{e}_{h,u}\left(\bm{v}_h,\frac{\gamma_0}{\Delta t}\hat{\hat{\bm{u}}}_h \right)
& - v_{h,\text{inhom}}^{e}(\bm{v}_h,\bm{g}_u(t_{n+1}),\bm{h}_u(t_{n+1}))
\; ,
\label{DualSplitting_ViscousStep_WeakForm}
\end{align}
for all~$\bm{v}_h \in \mathcal{V}^{u}_{h,e}$ and for all elements~$e=1,...,N_{\text{el}}$.

\begin{remark}
When considering the reference formulation used in~\cite{Hesthaven07} that does not perform integration by parts of the velocity divergence term and the pressure gradient term, the terms~$d^e_h$ in equation~\eqref{DualSplitting_Pressure_WeakForm} and~$g^e_h$ in equation~\eqref{DualSplitting_Projection_WeakForm} have to be replaced by~$d^e_{h,\text{ref}}$ defined in equation~\eqref{VelocityDivergence_ReferenceFormulation} and~$g^e_{h,\text{ref}}$ defined in equation~\eqref{PressureGradient_ReferenceFormulation}, respectively.
\end{remark}

\subsection{Pressure-correction scheme}\label{PressureCorrectionWeakForm}
Based on the DG formulation of basic derivative operators derived in Section~\ref{DGBasicOperators}, we summarize the discontinuous Galerkin discretization of pressure-correction schemes in this section.

\subsubsection{Momentum step}
The discontinuous Galerkin discretization of the time discrete momentum equation~\eqref{PressureConvection_MomentumStep_Impl} to be solved in the first substep of the pressure-correction scheme reads: Find~$\hat{\bm{u}}_h\in\mathcal{V}^u_h$ such that
\begin{align}
\begin{split}
m^{e}_{h,u}\left(\bm{v}_h,\frac{\gamma_0 \hat{\bm{u}}_h-\sum_{i=0}^{J-1}\left(\alpha_i\bm{u}^{n-i}_h\right)}{\Delta t} \right)
+ c^e_h\left(\bm{v}_h,\hat{\bm{u}}_h,\bm{g}_u(t_{n+1})\right)
+ v^e_h\left(\bm{v}_h,\hat{\bm{u}}_h,\bm{g}_u(t_{n+1}),\bm{h}_u(t_{n+1})\right)& \\
+ \sum_{i=0}^{J_p-1} \left(\beta_i g^e_h\left(\bm{v}_h,p^{n-i}_h,g_p(t_{n-i})\right) \right)
- b^e_h\left(\bm{v}_h,\bm{f}(t_{n+1})\right)& = 0 \; ,
\end{split} \label{PressureCorrection_MomentumImplicit_Nonlinear_WeakForm}
\end{align}
for all~$\bm{v}_h \in \mathcal{V}^{u}_{h,e}$ and for all elements~$e=1,...,N_{\text{el}}$. The boundary conditions~$\bm{g}_u$ and~$\bm{h}_u$ prescribed for the intermediate velocity~$\hat{\bm{u}}_h$ are defined in equations~\eqref{PressureCorrection_DirichletBC_IntermediateVelocity} and~\eqref{PressureCorrection_NeumannBC_IntermediateVelocity}, repectively. When solving the incremental formulation of the pressure-correction scheme, a boundary condition~$g_p$ according to equation~\eqref{NeumannBC_Pressure} has to be prescribed for the pressure on~$\Gamma_h^{\mathrm{N}}$ in order to evaluate the discrete pressure gradient operator.

\subsubsection{Pressure step}
The discontinuous Galerkin formulation of the pressure Poisson equation~\eqref{PressureCorrection_PressurePoissonEquation} is given as: Find~$\phi_h\in\mathcal{V}^p_h$ such that
\begin{align}
l_{h,\text{hom}}^{e}\left(q_h,\phi_h^{n+1}\right) =
- \frac{\gamma_0}{\Delta t} d_{h}^{e}\left(q_h,\hat{\bm{u}}_h,\bm{g}_u(t_{n+1})\right)
- l_{h,\text{inhom}}^{e}\left(q_h,g_{\phi}(t_{n+1}),h_{\phi}(t_{n+1})\right)\; ,
\label{PressureCorrection_PressureStep_WeakForm}
\end{align}
for all~$q_h \in \mathcal{V}^{p}_{h,e}$ and for all elements~$e=1,...,N_{\text{el}}$. The boundary values~$g_{\phi}$ and~$h_{\phi}$ are defined in equations~\eqref{PressureCorrection_PressureBC_GammaN} and~\eqref{PressureCorrection_PressureBC_GammaD}, respectively.

The approximate pressure solution~$p_h^{n+1}$ at time~$t_{n+1}$ is obtained from equation~\eqref{PressureCorrection_PressureUpdate}. In elementwise notation, the weak formulation of this pressure update reads: Find~$p_h^{n+1}\in\mathcal{V}^p_h$ such that
\begin{align}
m_{h,p}^{e}\left(q_h,p_h^{n+1}\right) = m_{h,p}^{e}\left(q_h,\phi_h^{n+1} + \sum_{i=0}^{J_p-1}\left(\beta_i p_h^{n-i}\right)\right) 
- \chi \nu\; d_{h}^{e}\left(q_h,\hat{\bm{u}}_h,\bm{g}_u\left(t_{n+1}\right)\right) \; ,
\label{PressureCorrection_PressureUpdate_WeakForm}
\end{align}
for all~$q_h \in \mathcal{V}^{p}_{h,e}$ and for all elements~$e=1,...,N_{\text{el}}$. The pressure mass matrix operator in the above equation is~$m^{e}_{h,p}\left(q_h,p_h\right) = \intele{q_h}{p_h}$. In contrast to the dual splitting scheme, the intermediate velocity field~$\hat{\bm{u}}$ fulfills the velocity Dirichlet boundary condition~$\bm{g}_u$ which can be seen from equation~\eqref{PressureCorrection_DirichletBC_IntermediateVelocity}. Consequently, this boundary condition is used to evaluate the discrete divergence operator applied to the intermediate velocity~$\hat{\bm{u}}_h$ on the right-hand side of equations~\eqref{PressureCorrection_PressureStep_WeakForm} and~\eqref{PressureCorrection_PressureUpdate_WeakForm}.

\subsubsection{Projection step}
The projection defined in equation~\eqref{PressureCorrection_Projection} finalizes time step~$n$ of the pressure-correction scheme. In terms of the weak discontinuous Galerkin formulation this local problem can be stated as: Find~$\bm{u}_h^{n+1}\in\mathcal{V}^u_h$ such that
\begin{align}
m_{h,u}^{e}\left(\bm{v}_h,\bm{u}_h^{n+1}\right) = m_{h,u}^{e}\left(\bm{v}_h,\hat{\bm{u}}_h\right)-\frac{\Delta t}{\gamma_0}g_h^{e}\left(\bm{v}_h,\phi_h^{n+1},g_{\phi}\left(t_{n+1}\right)\right)\; ,\label{PressureCorrection_Projection_WeakForm}
\end{align}
for all~$\bm{v}_h \in \mathcal{V}^{u}_{h,e}$ and for all elements~$e=1,...,N_{\text{el}}$. The pressure boundary condition~$g_{\phi}$ is defined in equation~\eqref{PressureCorrection_PressureBC_GammaN}.

\begin{remark}
When considering the reference formulation that does not perform integration by parts of the velocity divergence term and the pressure gradient term, the terms~$d^e_h$ in equations~\eqref{PressureCorrection_PressureStep_WeakForm} and~\eqref{PressureCorrection_PressureUpdate_WeakForm} and~$g^e_h$ in equations~\eqref{PressureCorrection_MomentumImplicit_Nonlinear_WeakForm} and~\eqref{PressureCorrection_Projection_WeakForm} have to be replaced by~$d^e_{h,\text{ref}}$ defined in equation~\eqref{VelocityDivergence_ReferenceFormulation} and~$g^e_{h,\text{ref}}$ defined in equation~\eqref{PressureGradient_ReferenceFormulation}, respectively.
\end{remark}

\subsection{Numerical integration}\label{NumericalIntegration}
Volume and surface integrals occuring in the weak formulations derived above are calculated using Gaussian quadrature. The number of quadrature points in each coordinate direction is selected in order to ensure exact integration on affine element geometries with constant Jacobian. In detail,~$n_{\rm{q}} = k_u+1$ quadrature points are used to integrate the velocity mass matrix term, the viscous term, the velocity divergence term, the pressure gradient term, and the body force term. Similarly,~$n_{\rm{q}}=k_p+1$ quadrature points are used to integrate the Laplace operator occuring in the pressure Poisson equation and the pressure mass matrix operator. Due to the nonlinearity of the convective term we use~$n_{\rm{q}}= \lfloor \frac{3k_u}{2} \rfloor +1$ quadrature points for the integration of the convective operator as well as for boundary face integrals containing the convective term to avoid aliasing effects.

\section{Numerical results}\label{NumericalResults}
In the following, we present numerical results for three different test cases. For an unsteady Stokes problem, instabilities in the limit of small time steps are analyzed for different formulations of the velocity divergence term and pressure gradient term and for both equal-order and mixed-order polynomials. As a means of verifying the results, the dual splitting scheme is compared to the pressure-correction scheme and coupled solution approach. Additionally, the aspect of inf--sup instabilities is analyzed by performing spatial convergence tests for equal-order and mixed-order polynomials.
Subsequently, using mixed-order polynomials we demonstrate optimal rates of convergence with respect to the temporal discretization and spatial discretization by considering an analytical solution of the full incompressible Navier--Stokes equations with non-trivial and time-dependent Dirichlet and Neumann boundary conditions. Finally, laminar flow around a cylinder with unsteady vortex shedding is considered in order to verify stability and accuracy of the presented approach for more complex flow problems.

\subsection{Implementation}\label{Implementation}
The solution of nonlinear systems of equations is based on a Newton--Krylov type solution approach. Linear(ized) system of equations are solved by using state of the art iterative methods such as the conjugate gradient (CG) method and the generalized mininum residual (GMRES) method. Unless otherwise specified, we use a relative solver tolerance of~$10^{-8}$ and an absolute solver tolerance of~$10^{-12}$ as tolerance criterion for all system of equations to be solved. In case of pure Dirichlet boundary conditions the pressure level is undefined resulting in a system of equations that is singular. In order to obtain a consistent system of equations we apply a transformation based on a subspace projection as described in~\cite{Krank16b}.

For the convergence tests presented in this section we use relative~$L^2$-errors that are defined as
\begin{align}
e_u=\frac{\Vert \bm{u}(\bm{x},t=T)-\bm{u}_h(\bm{x},t=T) \Vert_{L^2(\Omega_h)}}{\Vert \bm{u}(\bm{x},t=T) \Vert_{L^2(\Omega_h)}} \; , \; e_p=\frac{\Vert p(\bm{x},t=T)-p_h(\bm{x},t=T) \Vert_{L^2(\Omega_h)}}{\Vert p(\bm{x},t=T) \Vert_{L^2(\Omega_h)}} \; ,
\end{align}
where Gaussian quadrature is used to calculate the volume integrals in the above expressions. The number of one-dimensional quadrature points is~$k_{u}+3$ for the velocity error and~$k_p+3$ for the pressure error in order to ensure that the calculation of errors is not affected by quadrature errors. Experimental rates of convergence for two meshes with characterisitic element lengths~$h_1$ and~$h_2$ are calculated as~$\log(e_{h_1}/e_{h_2})/\log(h_1/h_2)$.

The code is implemented in \texttt{C++} and makes use of the object-oriented finite element library~\texttt{deal.II}~\cite{dealII85}. The incompressible Navier--Stokes solvers described above are implemented using a high-performance framework for generic finite element operator application developed in~\cite{Kronbichler12, Kormann16} that is based on a matrix-free implementation. 

\subsection{Unsteady Stokes equations}

We consider the unsteady Stokes flow problem analyzed in~\cite{Ferrer11,Ferrer14} in the context of instabilities occuring for small time step sizes. The analytical solution of the two-dimensional unsteady Stokes equations with~$\bm{f}=\bm{0}$ is defined as
\begin{align}
\begin{split}
\bm{u}(\bm{x},t) &=  \begin{pmatrix}
\sin(x_1)\left(a\sin(a x_2)-\cos(a)\sinh(x_2)\right)\\
\cos(x_1)\left(\cos(a x_2)+\cos(a)\cosh(x_2)\right)
\end{pmatrix}
\exp\left(-\lambda t\right)\; ,\\
p(\bm{x},t) &=  \lambda\cos(a)\cos(x_1)\sinh(x_2)\exp\left(-\lambda t\right)\; ,
\end{split}
\end{align}
where the parameters~$\lambda, \nu, a$ are given as~$\lambda = \nu (1+a^2)$ with~$\nu=1$ and~$a=2.883356$. The domain~$\Omega = [-L/2,L/2]^2$ is a square of length~$L=2$ and the time interval is~$[0,T]=[0,0.1]$. On domain boundaries~$\Gamma = \partial \Omega$ Dirichlet boundary conditions are prescribed,~$\Gamma =\Gamma^{\rm{D}}$. The Dirichlet boundary condition~$\bm{g}_u$, the time derivative term~$\partial \bm{g}_u/\partial t$ in equation~\eqref{DualSplitting_PressureBC_GammaD_2}, and initial conditions are deduced from the analytical solution.  The solution at previous instants of time~$t_{-J+1},...,t_{-1}$ required by the BDF scheme for~$J>1$ is obtained by interpolation of the analytical solution. Note that~$\nabla p \cdot \bm{n} \neq 0$ on domain boundaries. Hence, this flow example is a suitable test case to assess the temporal convergence properties of the projection-type solution methods. A uniform Cartesian grid consisting of quadrilateral elements of length~$h = L/2^l$ in~$x_1$ and~$x_2$-direction is used, where~$l$ denotes the level of refinement. To fix the pressure level the mean value of the vector containing the pressure degrees of freedom is set to zero. This is consistent with the exact pressure solution due to the symmetry of the analytical solution and the uniformity of the mesh.

\subsubsection{Instabilities in the limit of small time step sizes}\label{instabilities_for_small_time_steps}

\begin{figure}[t]
 \centering 
 \subfigure[equal-order polynomials]{
	\includegraphics[width=0.8\textwidth]{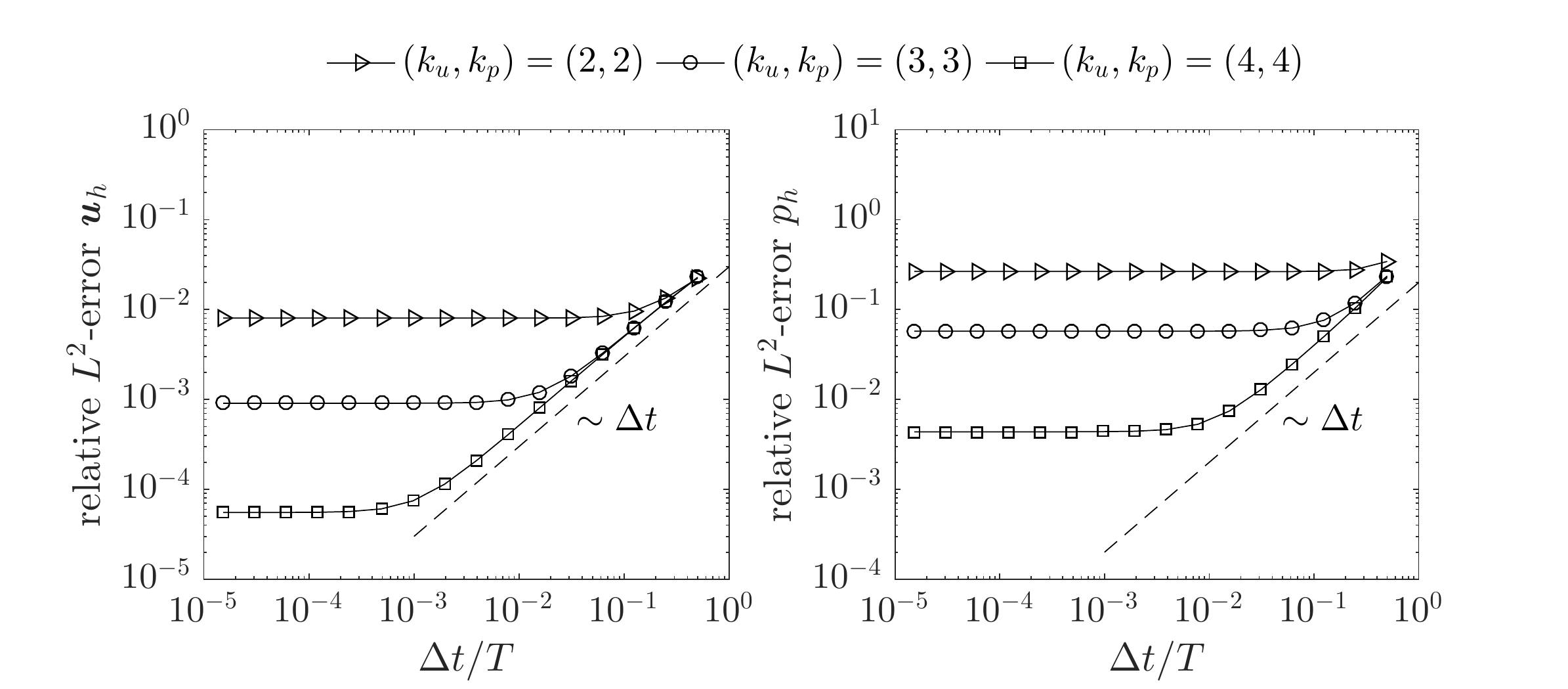}}
 \subfigure[mixed-order polynomials]{
	\includegraphics[width=0.8\textwidth]{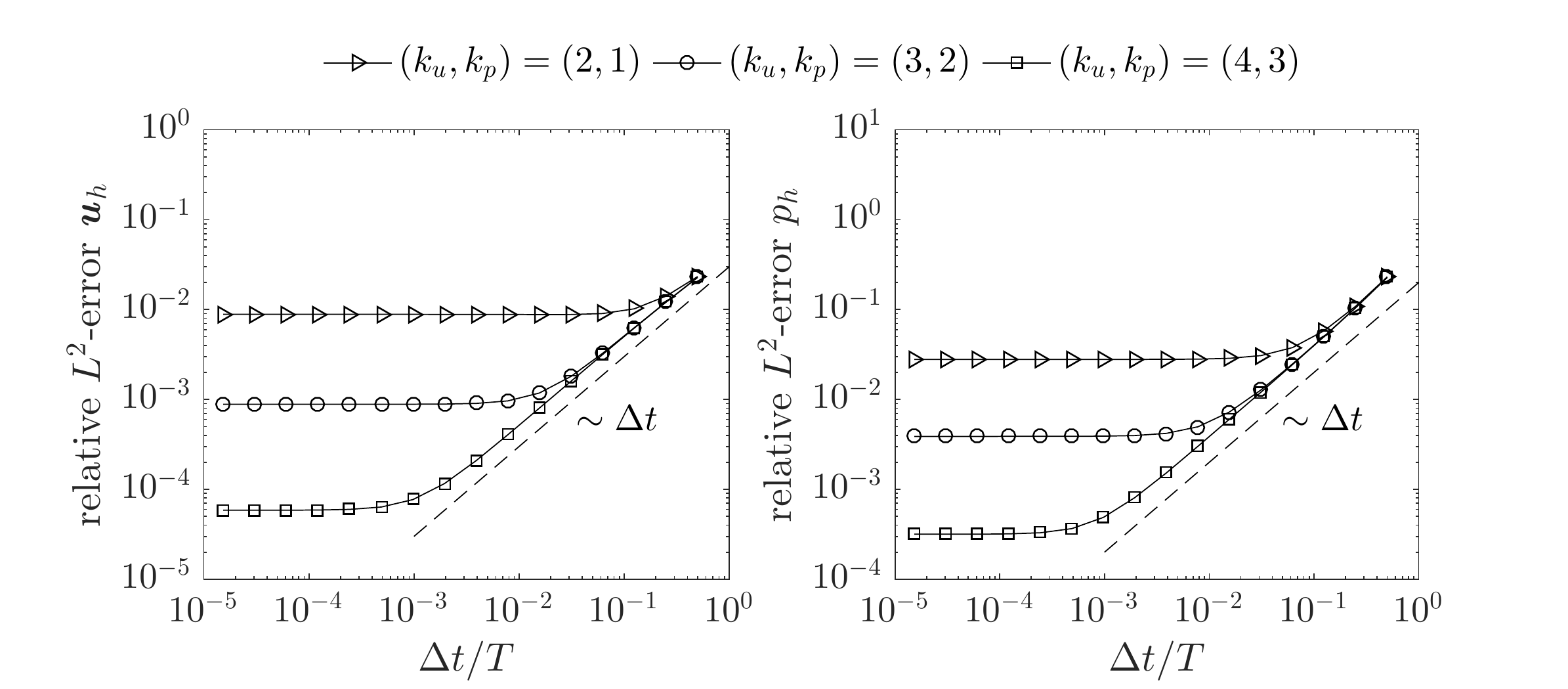}}
\caption{Stability analysis of coupled solution approach in the limit of small time step sizes for BDF1.}
\label{fig:stability_coupled_solver}
\end{figure}

\begin{figure}[t]
 \centering 
 \subfigure[equal-order polynomials]{
	\includegraphics[width=0.8\textwidth]{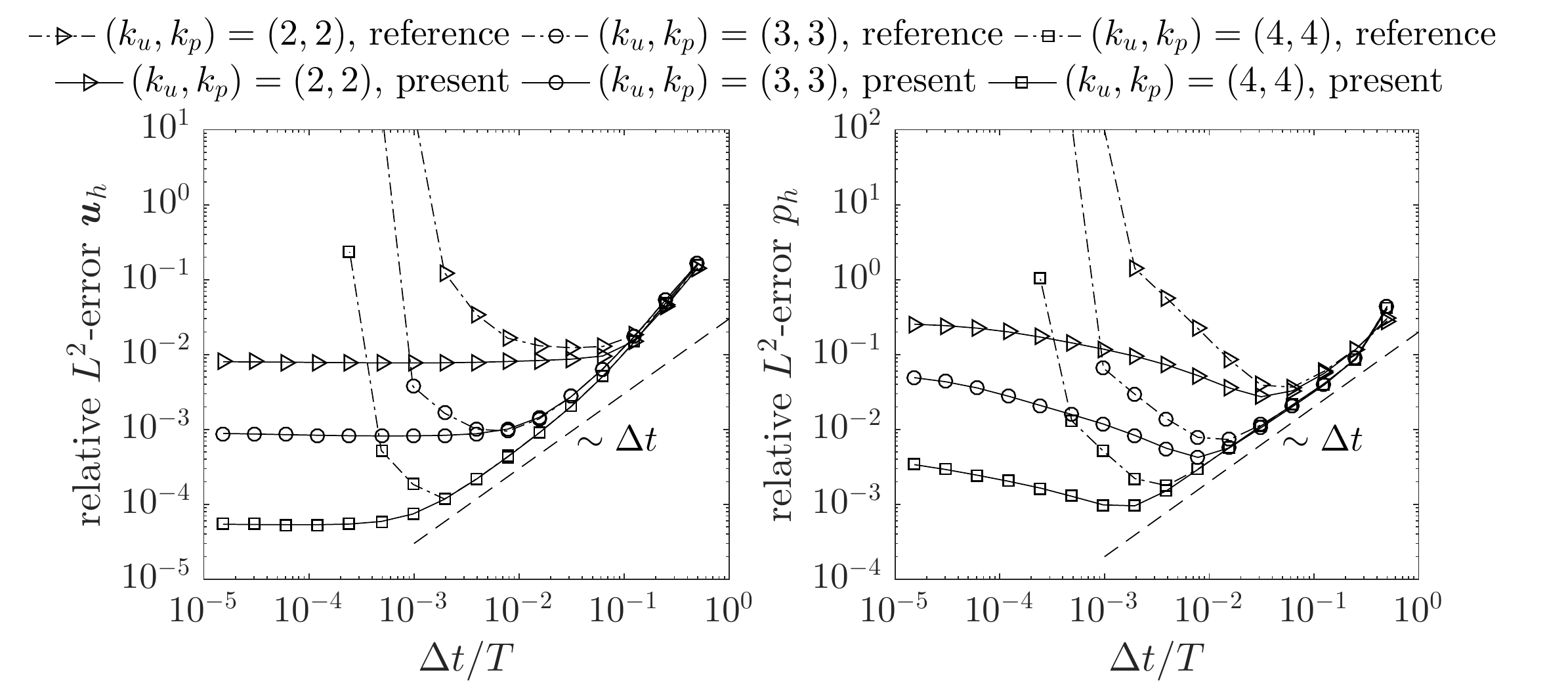}}
 \subfigure[mixed-order polynomials]{
	\includegraphics[width=0.8\textwidth]{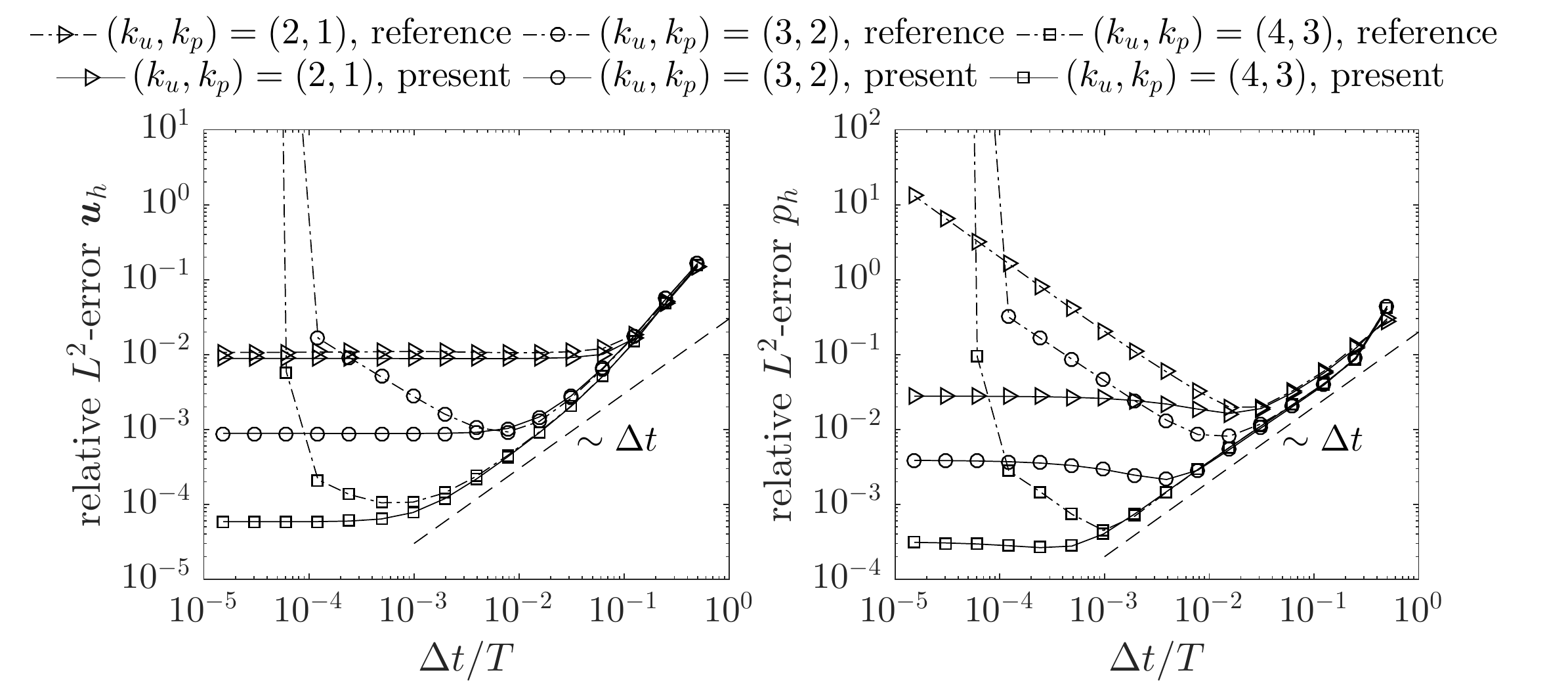}}
\caption{Stability analysis of dual splitting scheme in the limit of small time step sizes for BDF1.}
\label{fig:stability_dual_splitting}
\end{figure}

\begin{figure}[t]
 \centering 
 \subfigure[equal-order polynomials]{
	\includegraphics[width=0.8\textwidth]{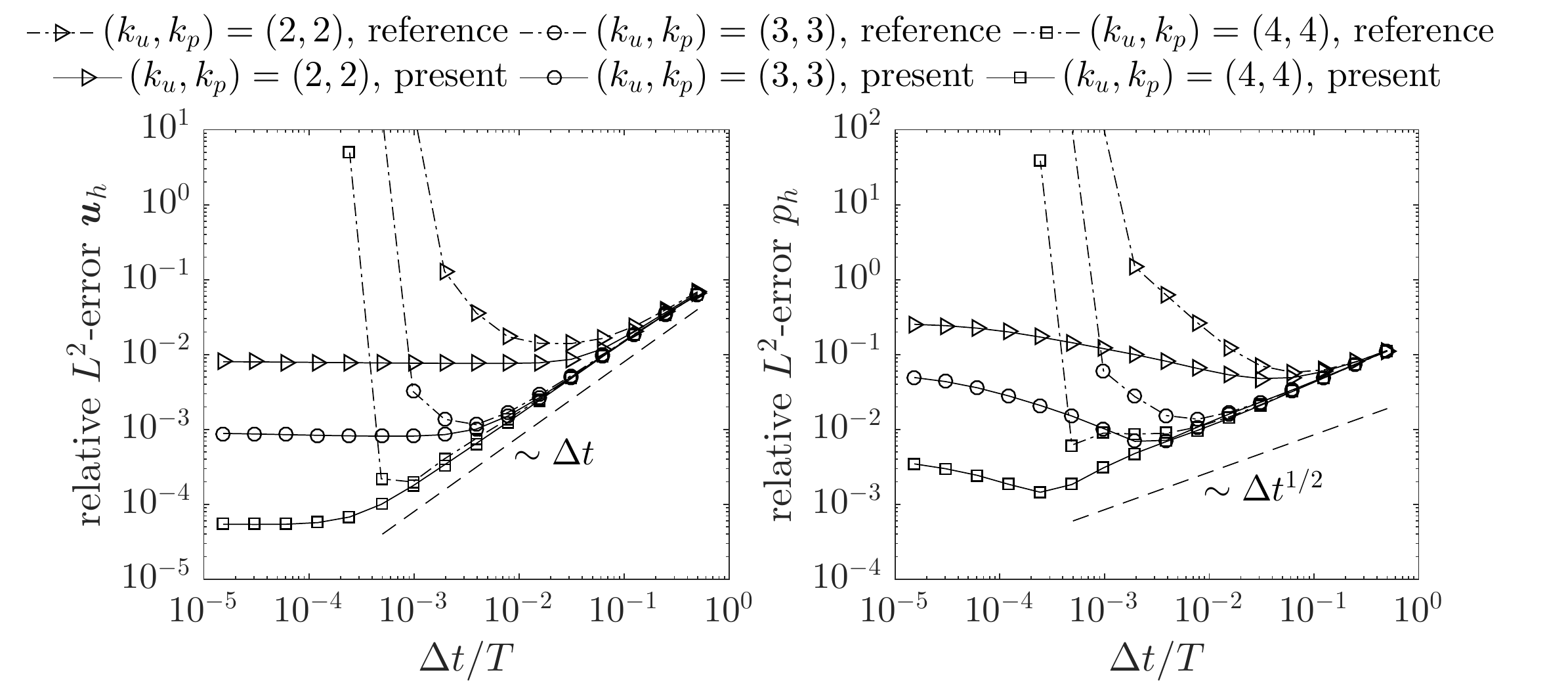}}
 \subfigure[mixed-order polynomials]{
	\includegraphics[width=0.8\textwidth]{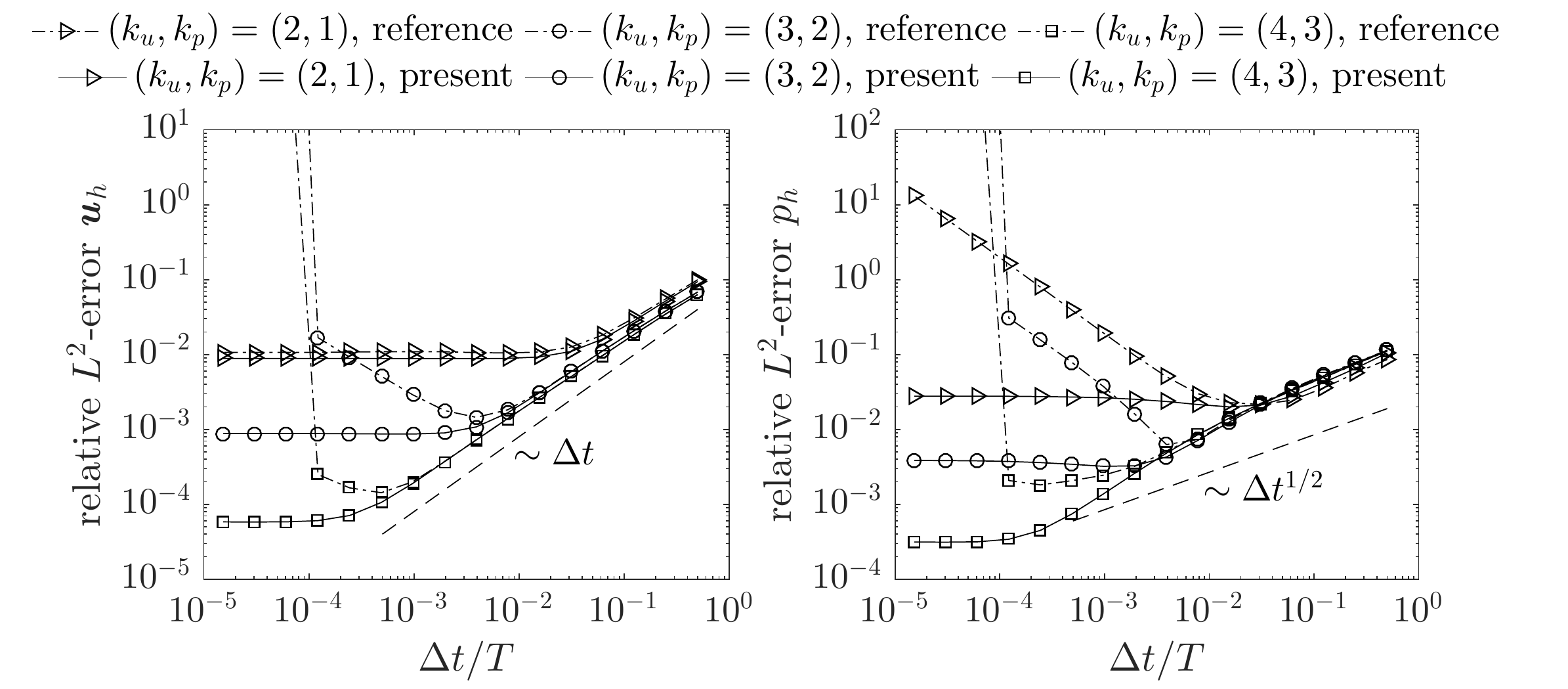}}
\caption{Stability analysis of pressure-correction scheme in the limit of small time step sizes for BDF1.}
\label{fig:stability_pressure_correction}
\end{figure}

To investigate the stability of the different solution approaches in the limit of small time step sizes we perform temporal convergence tests and vary the time step size~$\Delta t/T$ over a wide range. Since the instabilities reported in~\cite{Ferrer11,Ferrer14} occur primarily for coarse spatial resolutions, we select a coarse mesh with refine level~$l=2$. Moreover, the results are compared for both equal-order polynomials and mixed-order polynomials and varying polynomial degree. To show the impact of the temporal discretization error this analysis is performed for first order time integration schemes,~$J=1$. We note, however, that qualitatively similar results in terms of stability are obtained when using second order accurate time integration schemes.

The results for the coupled solution approach presented in Figure~\ref{fig:stability_coupled_solver} show the expected behavior. The error is proportional to~$\Delta t$ for large time steps. For small time steps the temporal error becomes negligible as compared to the spatial discretization error and the overall error approaches a constant value. The pressure error is significantly larger for equal-order polynomials than for mixed-order polynomials while the velocity error is almost the same for both equal-order and mixed-order polynomials. This apsect is analyzed in more detail in Section~\ref{InfSupStability}.

Figure~\ref{fig:stability_dual_splitting} shows results obtained for the dual splitting scheme. As a reference method we consider the DG discretization proposed in~\cite{Hesthaven07} which is based on a modified formulation of the velocity divergence term and pressure gradient term as explained in Section~\ref{WeakFormDualSplitting}. As in~\cite{Ferrer11,Ferrer14}, we observe instabilities in the limit of small time step sizes. We emphasize that these instabilities occur similarly for equal-order and mixed-order polynomials. The situation changes completely, however, when applying the DG formulations~$d^e_h$ and~$g^e_h$ according to equations~\eqref{WeakFormVelocityDivergenceBoundaryIntegrals} and~\eqref{WeakFormPressureGradientBoundaryIntegrals}, repectively, along with the consistent boundary condition~\eqref{DualSplitting_DBC_IntermediateVelocity} for the intermediate velocity field proposed in the present paper. For this formulation, stability is obtained for both equal-order and mixed-order polynomials and the behavior in the limit of small time steps is comparable to the coupled solution approach.

As a further verification of the results, we perform the same simulations for the non-incremental pressure-correction scheme in standard formulation. The results are shown in Figure~\ref{fig:stability_pressure_correction}. Again, we compare the DG formulation used in the present work to a reference method that is based on the modified discretization of the velocity divergence term and pressure gradient term as explained in Section~\ref{PressureCorrectionWeakForm}. The stability behavior is the same as for the dual splitting scheme. We note that for the pressure-correction scheme the reference slope is~$\Delta t^{1/2}$ for the pressure, representing the expected theoretical rate of convergence. 

The above results lead to the following conclusions. Our results are clearly in contradiction to the conlusions drawn in~\cite{Ferrer14} where it is stated that the instabilities in the limit of small time steps are related to the temporal discretization and inf--sup instabilities but not to the spatial discretization. Instead, our results suggest that these instabilities are clearly related to the discontinuous Galerkin formulation of the velocity divergence term and the pressure gradient term. The discretization of these operators is a basic building block of any incompressible Navier--Stokes solver independent of the temporal discretization approach. In fact, the increased pressure error observed for the equal-order formulation is an indication of inf--sup instabilities and is discussed in more detail below. As for the Stokes flow problem analyzed here, we observed a qualitatively similar behavior in terms of instabilities for small time step sizes and the different formulations of the velocity divergence term and pressure gradient term when considering the full Navier--Stokes equations, e.g., for the two-dimensional vortex problem considered in Section~\ref{VortexProblem} or the three-dimensional Beltrami flow problem~\cite{Ethier1994}.

\subsubsection{Temporal convergence test}\label{TemporalConvergenceStokes}
To assess the temporal accuracy of the considered solution approaches for different orders of the time integration scheme, the relative~$L^2$-errors of velocity and pressure are calculated for time step sizes~$\Delta t/T = 1/2^m$ with~$m=1,...,12$ and we use a refine level of~$l=4$ with mixed-order polynomials~$(k_u,k_p)=(6,5)$ so that the spatial discretization error has a negligible impact on the overall accuracy. BDF schemes of order~$J=1$ and~$J=2$ are analyzed for all solution strategies. In case of the pressure-correction scheme the order of extrapolation of the pressure gradient term is~$J_p=J-1$ and both the standard formulation~($\chi=0$) and the rotational formulation~($\chi=1$) are considered, see equation~\eqref{PressureCorrection_Phi}.

\begin{figure}[t]
 \centering 
	\includegraphics[width=1.0\textwidth]{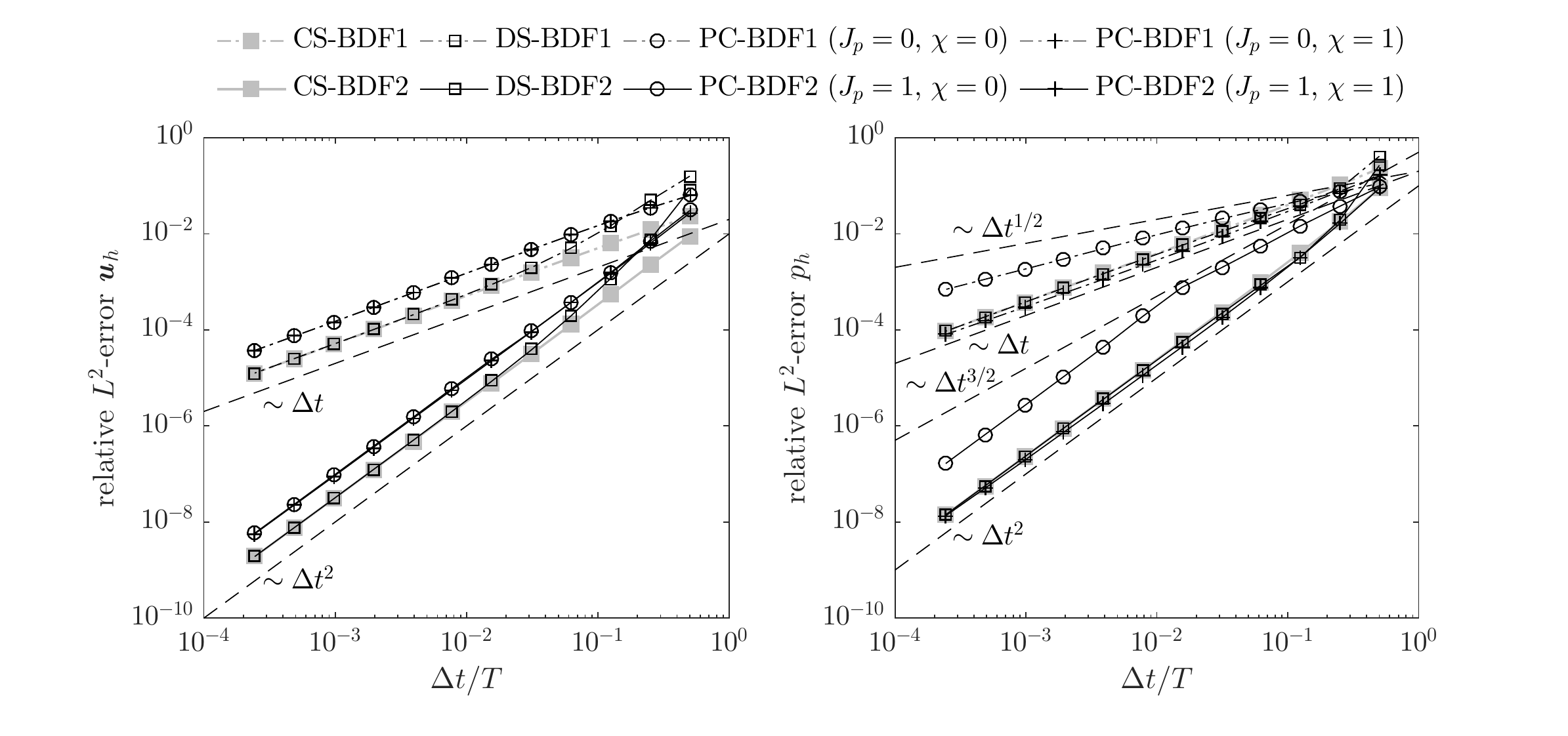}

\caption{Stokes flow problem: temporal convergence tests for coupled solver (CS), dual splitting scheme (DS), and pressure-correction scheme (PC) in standard form~($\chi=0$) and rotational form~($\chi=1$) for BDF schemes of order~$J=1$ and~$J=2$. The spatial resolution is~$l=4$ and~$(k_u,k_p)=(6,5)$.}
\label{fig:temporal_convergence_stokes}
\end{figure}

Results of the temporal convergence test are presented in Figure~\ref{fig:temporal_convergence_stokes}. For the coupled solution approach and the dual splitting scheme, both velocity and pressure converge with optimal rates of convergence of order~$\Delta t^1$ for the BDF1 scheme and~$\Delta t^2$ for the BDF2 scheme. For the pressure-correction scheme, the velocity error converges also with optimal rates of convergence of order~$\Delta t^J$. The velocity error is almost the same for the standard formulation and the rotational formulation, but the absolute error is approximately a factor of~$3$ larger as compared to the coupled solution approach. For the standard formulation of the pressure-correction scheme the pressure converges with a rate between~$\Delta t^{1/2}$ and~$\Delta t^1$ for the BDF1 scheme (experimental rates of convergence are between~$0.6$ and~$0.7$) and a rate of convergence close to~$\Delta t^{3/2}$ for BDF2 (a transition of the convergence rate from~$\Delta t^{3/2}$ to~$\Delta t^{2}$ can be observed where the point of transition depends on the spatial resolution as described in~\cite[Section 7]{Guermond06}). Using the rotational formulation significantly reduces the pressure error and optimal convergence rates of order~$\Delta t^J$ are obtained as for the coupled solution approach and the dual splitting scheme.

\subsubsection{Spatial convergence test and numerical investigation of inf--sup stability}\label{InfSupStability}
In this section, we analyze the inf--sup stability of projection methods numerically by performing simulations for both equal-order polynomials and mixed-order polynomials. As argued in~\cite{Guermond06}, the fact that the pressure Poisson equation and Helmholtz equation for the velocity are solvable independently of the polynomial spaces used to represent the velocity and pressure solutions in case of projection methods does not mean that the inf--sup condition is irrelevant for projection methods. Instead, the corresponding steady-state Stokes problem of projection methods is the decisive metric to evaluate the need of the inf--sup condition.

Following~\cite{Guermond06} and~\cite{Ferrer14}, we briefly derive the steady-state Stokes equations for the dual splitting scheme and the pressure-correction scheme from the equations shown in Section~\ref{TemporalDiscretization} as a basis for the interpretation of numerical results shown below. For the ease of notation, we consider the equations at the level of differential operators, but similar relations can be derived for discretized operators or matrix formulations. For the dual splitting scheme, the following system of equations can be derived for the steady Stokes problem
\begin{align}
\begin{pmatrix}
-\Div{\bm{F}_{\rm{v}}(\bm{u})} & +\Grad{p}\\
-\Div{\bm{u}} & +\frac{\Delta t}{\gamma_0}\nabla^2 p
\end{pmatrix} =
\begin{pmatrix}
\bm{f}\\
\frac{\Delta t}{\gamma_0}\Div{\bm{f}}
\end{pmatrix} \; .
\label{InfSupDualSplitting}
\end{align}
The first equation is obtained by adding equations~\eqref{DualSplitting_ConvectiveStep},~\eqref{ProjectionMethodEq1} and~\eqref{DualSplitting_ViscousStep}, neglecting the convective term in equation~\eqref{DualSplitting_ConvectiveStep}, assuming steady state~$\bm{u}^{n+1}=\bm{u}^{n} = ...=\bm{u}^{n-J+1} = \bm{u}$, and using the fact that the BDF time integration constants fulfill the property~$\gamma_0 = \sum_{i=0}^{J-1}\alpha_i$. The second equation is obtained by taking the divergence of equation~\eqref{ProjectionMethodEq1} and inserting equations~\eqref{ProjectionMethodEq2} and~\eqref{DualSplitting_ConvectiveStep}. Equation~\eqref{InfSupDualSplitting} highlights that the dual splitting scheme introduces an inf--sup stabilization term~$\Delta t /\gamma_0 \nabla^2 p$ that is scaled by the time step size~$\Delta t$. According to this relation, the impact of the stabilization diminishes for small time step sizes~$\Delta t$. 

Similarly, the following system of equations can be derived for the pressure-correction scheme
\begin{align}
\begin{pmatrix}
-\Div{\bm{F}_{\rm{v}}(\hat{\bm{u}})} + \Grad{(\chi \nu \Div{\hat{\bm{u}}})} & +\Grad{p}\\
-\Div{\hat{\bm{u}}} + \frac{\Delta t}{\gamma_0}\Div{\left(\Grad{(\chi \nu\Div{\hat{\bm{u}}})} \right)} & +\frac{\Delta t}{\gamma_0}\left(1-\sum_{i=0}^{J_p-1}\beta_i\right)\nabla^2 p
\end{pmatrix} =
\begin{pmatrix}
\bm{f}\\
0
\end{pmatrix} \; .
\label{InfSupPressureCorrection}
\end{align}
To derive these equations one assumes that the solution reaches a steady state,~$\bm{u}^{n+1}=\bm{u}^{n} = ...=\bm{u}^{n-J+1} = \bm{u}$ and~$p^{n+1}=p^{n}=...=p^{n-J_p+1}=p$, and uses the fact that the time integration constants fulfill~$\gamma_0 = \sum_{i=0}^{J-1}\alpha_i$ and~$\sum_{i=0}^{J_p-1} \beta_i = 1$. The first equation is obtained by adding equations~\eqref{PressureConvection_MomentumStep_Impl} and~\eqref{PressureCorrection_ProjectionMethodEq1}, neglecting the convective term in equation~\eqref{PressureConvection_MomentumStep_Impl}, and replacing~$\phi^{n+1}$ by equation~\eqref{PressureCorrection_Phi}. The second equation is derived by taking the divergence of equation~\eqref{PressureCorrection_ProjectionMethodEq1} and inserting equations~\eqref{PressureCorrection_ProjectionMethodEq2} and~\eqref{PressureCorrection_Phi}. Equation~\eqref{InfSupPressureCorrection} highlights that the pressure-correction scheme introduces an inf--sup stabilization term in case of the nonincremental formulation ($J_p=0$) but not in case of the incremental formulation since~$1-\sum_{i=0}^{J_p-1} \beta_i = 0$ for~$J_p\geq 1$.

\begin{table}[!h]
\caption{Spatial convergence tests for Stokes flow problem for coupled solution approach, dual splitting scheme, and pressure-correction scheme. Relative~$L^2$-errors and experimental rates of convergence are reported in form of tuples composed of two values where the first value corresponds to the equal-order formulation~$(k_u,k_p)=(k,k)$ and the second value to the mixed-order formulation~$(k_u,k_p)=(k,k-1)$.}\label{SpatialConvergenceStokes}
\renewcommand{\arraystretch}{1.06}
\begin{scriptsize}
\begin{center}
\subtable[Coupled solution approach]{
\begin{tabular}{ccccccc}
\hline
& & \multicolumn{2}{c}{velocity (equal-order (k,k) $\vert$ mixed-order (k,k-1))} & & \multicolumn{2}{c}{pressure (equal-order (k,k) $\vert$ mixed-order (k,k-1))} \\ 
\cline{3-4} \cline{6-7} $k$ & $h$ & relative $L^2$-error  & rate of convergence & & relative $L^2$-error & rate of convergence\\ 
\hline 
2 & $L/2$  & 5.00E--002 $\vert$ 5.39E--002 &                   & & 1.30E+000  $\vert$ 8.28E--002 & \\  
  & $L/4$  & 8.06E--003 $\vert$ 8.85E--003 & 2.63 $\vert$ 2.60 & & 2.66E--001 $\vert$ 2.79E--002 & 2.29 $\vert$ 1.57\\  
  & $L/8$  & 1.07E--003 $\vert$ 1.16E--003 & 2.92 $\vert$ 2.94 & & 3.47E--002 $\vert$ 6.69E--003 & 2.94 $\vert$ 2.06\\  
  & $L/16$ & 1.34E--004 $\vert$ 1.44E--004 & 2.99 $\vert$ 3.01 & & 5.96E--003 $\vert$ 1.31E--003 & 2.54 $\vert$ 2.36\\   
\hline 
3 & $L/2$  & 1.45E--002 $\vert$ 1.55E--002 &                   & & 1.59E--001 $\vert$ 5.78E--002 & \\  
  & $L/4$  & 9.08E--004 $\vert$ 8.86E--004 & 4.00 $\vert$ 4.13 & & 5.74E--002 $\vert$ 3.88E--003 & 1.47 $\vert$ 3.90\\  
  & $L/8$  & 7.82E--005 $\vert$ 5.85E--005 & 3.54 $\vert$ 3.92 & & 2.41E--002 $\vert$ 2.72E--004 & 1.25 $\vert$ 3.84\\  
  & $L/16$ & 8.07E--006 $\vert$ 3.79E--006 & 3.28 $\vert$ 3.95 & & 1.09E--002 $\vert$ 2.41E--005 & 1.15 $\vert$ 3.50\\   
\hline
4 & $L/2$  & 1.38E--003 $\vert$ 1.39E--003 &                   & & 8.95E--002 $\vert$ 2.80E--003 & \\  
  & $L/4$  & 5.53E--005 $\vert$ 5.86E--005 & 4.64 $\vert$ 4.57 & & 4.37E--003 $\vert$ 3.17E--004 & 4.36 $\vert$ 3.15\\  
  & $L/8$  & 1.71E--006 $\vert$ 1.85E--006 & 5.02 $\vert$ 4.99 & & 1.42E--004 $\vert$ 1.65E--005 & 4.94 $\vert$ 4.24\\  
  & $L/16$ & 5.33E--008 $\vert$ 5.75E--008 & 5.00 $\vert$ 5.01 & & 6.35E--006 $\vert$ 7.17E--007 & 4.49 $\vert$ 4.54\\   
\hline
5 & $L/2$  & 2.55E--004 $\vert$ 2.62E--004 &                   & & 4.24E--003 $\vert$ 1.76E--003 & \\  
  & $L/4$  & 3.94E--006 $\vert$ 3.78E--006 & 6.02 $\vert$ 6.12 & & 5.51E--004 $\vert$ 2.54E--005 & 2.94 $\vert$ 6.11\\  
  & $L/8$  & 8.36E--008 $\vert$ 6.14E--008 & 5.56 $\vert$ 5.94 & & 5.81E--005 $\vert$ 4.27E--007 & 3.25 $\vert$ 5.89\\  
  & $L/16$ & 2.18E--009 $\vert$ 1.03E--009 & 5.26 $\vert$ 5.90 & & 6.54E--006 $\vert$ 9.13E--009 & 3.15 $\vert$ 5.55\\   
\hline
\end{tabular}}

\subtable[Dual splitting scheme]{
\begin{tabular}{ccccccc}
\hline
& & \multicolumn{2}{c}{velocity (equal-order (k,k) $\vert$ mixed-order (k,k-1))} & & \multicolumn{2}{c}{pressure (equal-order (k,k) $\vert$ mixed-order (k,k-1))} \\ 
\cline{3-4} \cline{6-7} $k$ & $h$ & relative $L^2$-error  & rate of convergence & & relative $L^2$-error & rate of convergence\\ 
\hline 
2 & $L/2$  & 4.92E--002 $\vert$ 5.38E--002 &                   & & 1.27E+000  $\vert$ 8.27E--002 & \\  
  & $L/4$  & 7.91E--003 $\vert$ 8.85E--003 & 2.64 $\vert$ 2.60 & & 2.23E--001 $\vert$ 2.77E--002 & 2.52 $\vert$ 1.58\\  
  & $L/8$  & 1.06E--003 $\vert$ 1.16E--003 & 2.89 $\vert$ 2.94 & & 2.49E--002 $\vert$ 6.48E--003 & 3.16 $\vert$ 2.10\\  
  & $L/16$ & 1.34E--004 $\vert$ 1.44E--004 & 2.98 $\vert$ 3.01 & & 4.39E--003 $\vert$ 1.18E--003 & 2.50 $\vert$ 2.46\\   
\hline 
3 & $L/2$  & 1.42E--002 $\vert$ 1.55E--002 &                   & & 1.56E--001 $\vert$ 5.73E--002 & \\  
  & $L/4$  & 8.50E--004 $\vert$ 8.85E--004 & 4.09 $\vert$ 4.13 & & 3.47E--002 $\vert$ 3.79E--003 & 2.17 $\vert$ 3.92\\  
  & $L/8$  & 5.54E--005 $\vert$ 5.85E--005 & 3.94 $\vert$ 3.92 & & 3.26E--003 $\vert$ 2.58E--004 & 3.41 $\vert$ 3.88\\  
  & $L/16$ & 3.52E--006 $\vert$ 3.79E--006 & 3.98 $\vert$ 3.95 & & 1.72E--004 $\vert$ 2.11E--005 & 4.24 $\vert$ 3.61\\   
\hline
4 & $L/2$  & 1.25E--003 $\vert$ 1.39E--003 &                   & & 7.84E--002 $\vert$ 2.71E--003 & \\  
  & $L/4$  & 5.34E--005 $\vert$ 5.86E--005 & 4.55 $\vert$ 4.57 & & 2.39E--003 $\vert$ 2.93E--004 & 5.03 $\vert$ 3.21\\  
  & $L/8$  & 1.70E--006 $\vert$ 1.85E--006 & 4.97 $\vert$ 4.99 & & 7.91E--005 $\vert$ 1.11E--005 & 4.92 $\vert$ 4.72\\  
  & $L/16$ & 5.38E--008 $\vert$ 5.75E--008 & 4.98 $\vert$ 5.00 & & 2.88E--006 $\vert$ 2.50E--007 & 4.78 $\vert$ 5.48\\   
\hline
5 & $L/2$  & 2.55E--004 $\vert$ 2.62E--004 &                   & & 3.74E--003 $\vert$ 1.65E--003 & \\  
  & $L/4$  & 3.57E--006 $\vert$ 3.78E--006 & 6.16 $\vert$ 6.12 & & 1.72E--004 $\vert$ 2.20E--005 & 4.45 $\vert$ 6.23\\  
  & $L/8$  & 5.77E--008 $\vert$ 6.14E--008 & 5.95 $\vert$ 5.94 & & 2.82E--006 $\vert$ 3.23E--007 & 5.93 $\vert$ 6.09\\  
  & $L/16$ & 9.93E--010 $\vert$ 1.05E--009 & 5.86 $\vert$ 5.87 & & 3.89E--008 $\vert$ 4.91E--009 & 6.18 $\vert$ 6.04\\      
\hline
\end{tabular}  }

\subtable[Pressure-correction scheme]{
\begin{tabular}{ccccccc}
\hline
& & \multicolumn{2}{c}{velocity (equal-order (k,k) $\vert$ mixed-order (k,k-1))} & & \multicolumn{2}{c}{pressure (equal-order (k,k) $\vert$ mixed-order (k,k-1))} \\ 
\cline{3-4} \cline{6-7} $k$ & $h$ & relative $L^2$-error  & rate of convergence & & relative $L^2$-error & rate of convergence\\ 
\hline 
2 & $L/2$  & 5.00E--002 $\vert$ 5.39E--002 &                   & & 1.30E+000  $\vert$ 8.28E--002 & \\  
  & $L/4$  & 8.06E--003 $\vert$ 8.85E--003 & 2.63 $\vert$ 2.60 & & 2.66E--001 $\vert$ 2.79E--002 & 2.29 $\vert$ 1.57\\  
  & $L/8$  & 1.07E--003 $\vert$ 1.16E--003 & 2.92 $\vert$ 2.94 & & 3.47E--002 $\vert$ 6.69E--003 & 2.94 $\vert$ 2.06\\  
  & $L/16$ & 1.34E--004 $\vert$ 1.44E--004 & 2.99 $\vert$ 3.01 & & 5.97E--003 $\vert$ 1.31E--003 & 2.54 $\vert$ 2.36\\   
\hline 
3 & $L/2$  & 1.45E--002 $\vert$ 1.55E--002 &                   & & 1.59E--001 $\vert$ 5.78E--002 & \\  
  & $L/4$  & 9.08E--004 $\vert$ 8.86E--004 & 4.00 $\vert$ 4.13 & & 5.74E--002 $\vert$ 3.88E--003 & 1.47 $\vert$ 3.90\\  
  & $L/8$  & 7.82E--005 $\vert$ 5.85E--005 & 3.54 $\vert$ 3.92 & & 2.42E--002 $\vert$ 2.72E--004 & 1.25 $\vert$ 3.84\\  
  & $L/16$ & 8.11E--006 $\vert$ 3.79E--006 & 3.27 $\vert$ 3.95 & & 1.10E--002 $\vert$ 2.41E--005 & 1.14 $\vert$ 3.50\\   
\hline
4 & $L/2$  & 1.38E--003 $\vert$ 1.39E--003 &                   & & 8.95E--002 $\vert$ 2.80E--003 & \\  
  & $L/4$  & 5.53E--005 $\vert$ 5.86E--005 & 4.64 $\vert$ 4.57 & & 4.37E--003 $\vert$ 3.17E--004 & 4.36 $\vert$ 3.15\\  
  & $L/8$  & 1.71E--006 $\vert$ 1.85E--006 & 5.02 $\vert$ 4.99 & & 1.42E--004 $\vert$ 1.67E--005 & 4.94 $\vert$ 4.24\\  
  & $L/16$ & 5.33E--008 $\vert$ 5.76E--008 & 5.00 $\vert$ 5.00 & & 6.22E--006 $\vert$ 7.16E--007 & 4.51 $\vert$ 4.54\\   
\hline
5 & $L/2$  & 2.55E--004 $\vert$ 2.62E--004 &                   & & 4.24E--003 $\vert$ 1.76E--003 & \\  
  & $L/4$  & 3.94E--006 $\vert$ 3.78E--006 & 6.02 $\vert$ 6.12 & & 5.51E--004 $\vert$ 2.54E--005 & 2.94 $\vert$ 6.11\\  
  & $L/8$  & 8.37E--008 $\vert$ 6.13E--008 & 5.56 $\vert$ 5.94 & & 5.85E--005 $\vert$ 4.27E--007 & 3.24 $\vert$ 5.90\\  
  & $L/16$ & 2.37E--009 $\vert$ 1.36E--009 & 5.14 $\vert$ 5.49 & & 7.00E--006 $\vert$ 8.93E--009 & 3.06 $\vert$ 5.58\\      
\hline
\end{tabular}  }
\end{center}
\end{scriptsize}
\renewcommand{\arraystretch}{1}
\end{table}

\begin{figure}[t]
 \centering 
 \subfigure[$\Delta t/T = 10^{-1}$, equal-order formulation~$(k_u,k_p)=(k,k)$, refine level~$l=4-k$]{
	\includegraphics[width=0.2\textwidth]{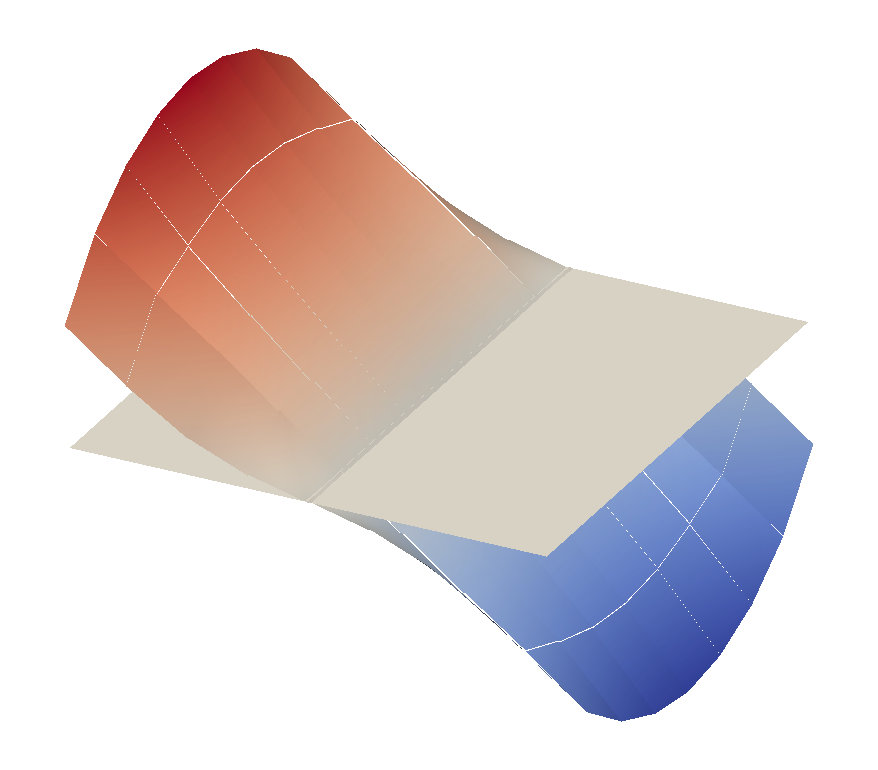}
    \includegraphics[width=0.2\textwidth]{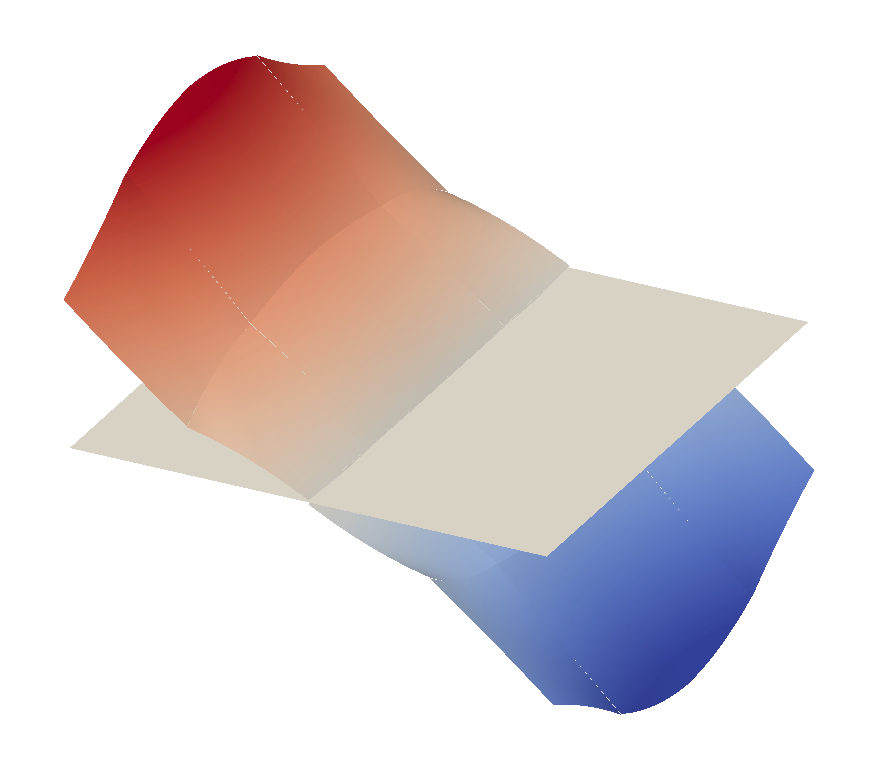}
    \includegraphics[width=0.2\textwidth]{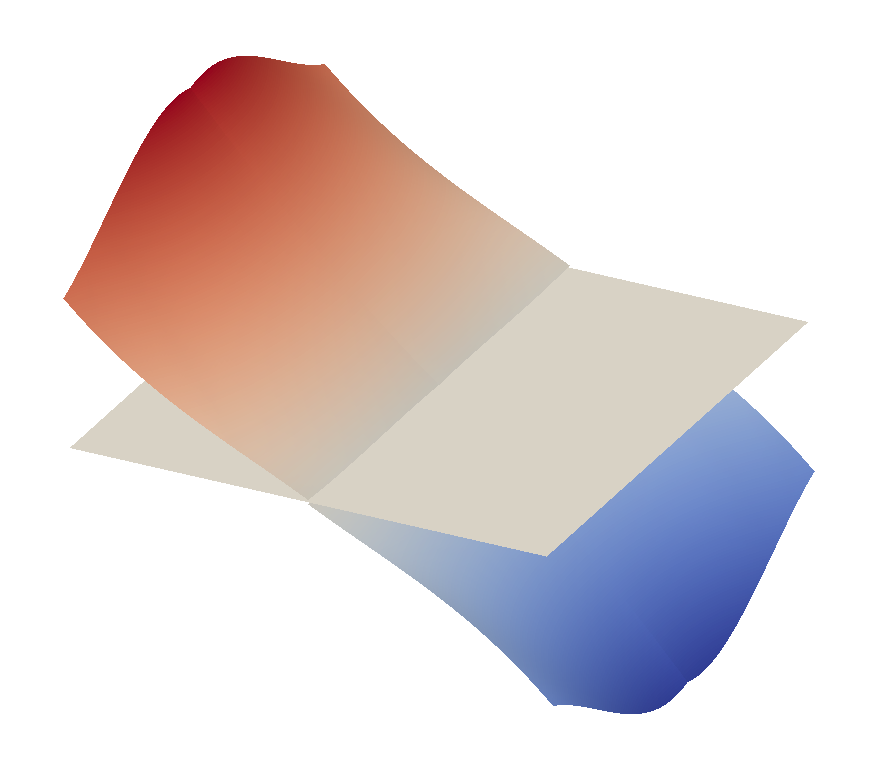}
    \includegraphics[width=0.2\textwidth]{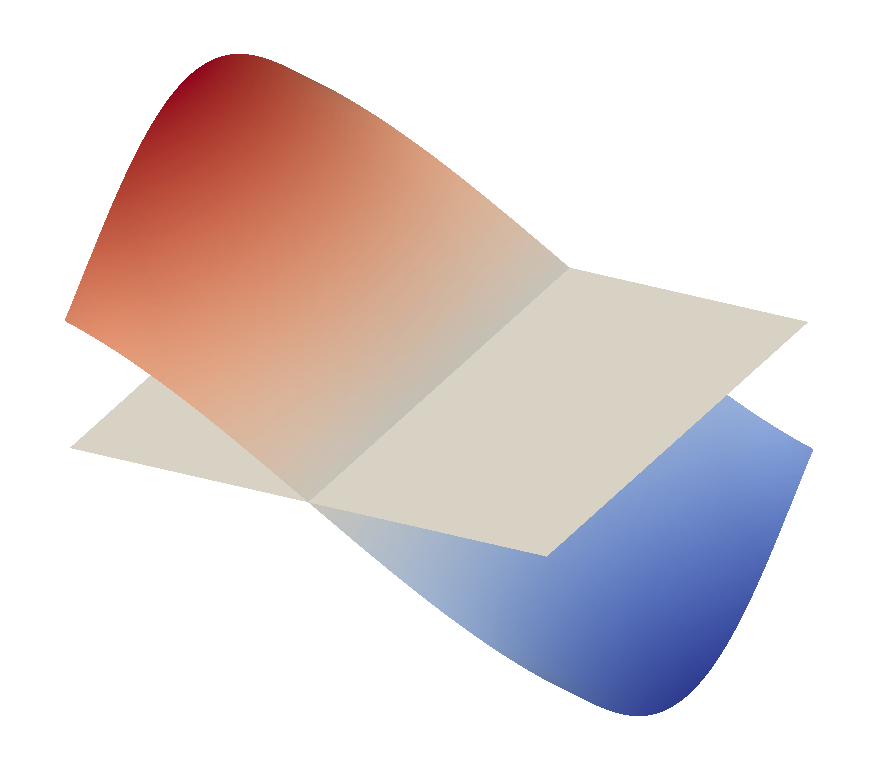}
}
    
\subfigure[$\Delta t/T = 10^{-3}$, equal-order formulation~$(k_u,k_p)=(k,k)$, refine level~$l=4-k$]{
    \includegraphics[width=0.2\textwidth]{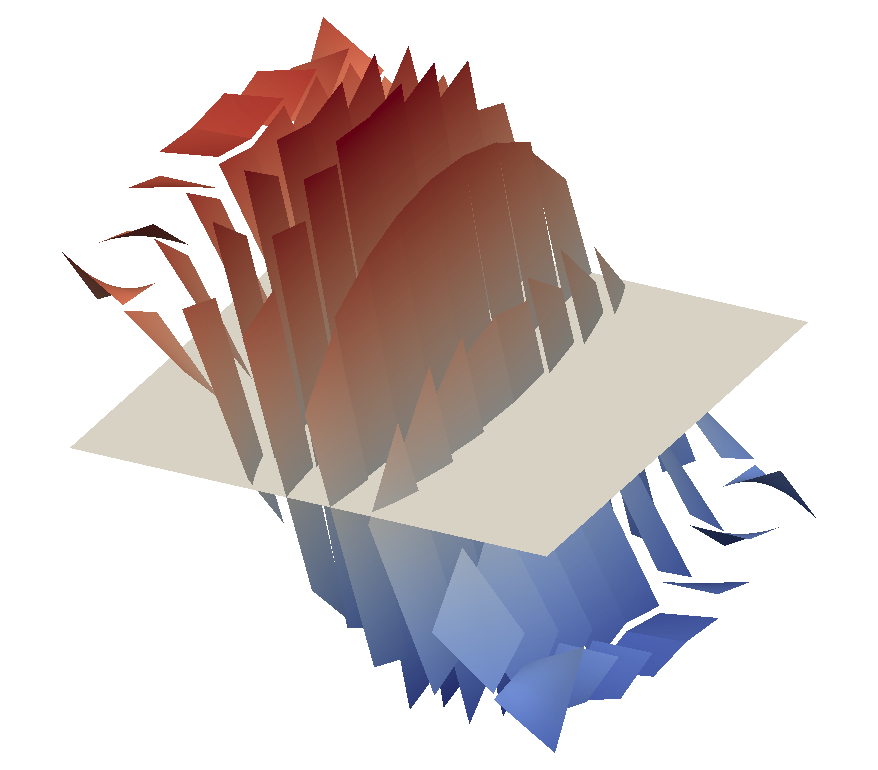}
    \includegraphics[width=0.2\textwidth]{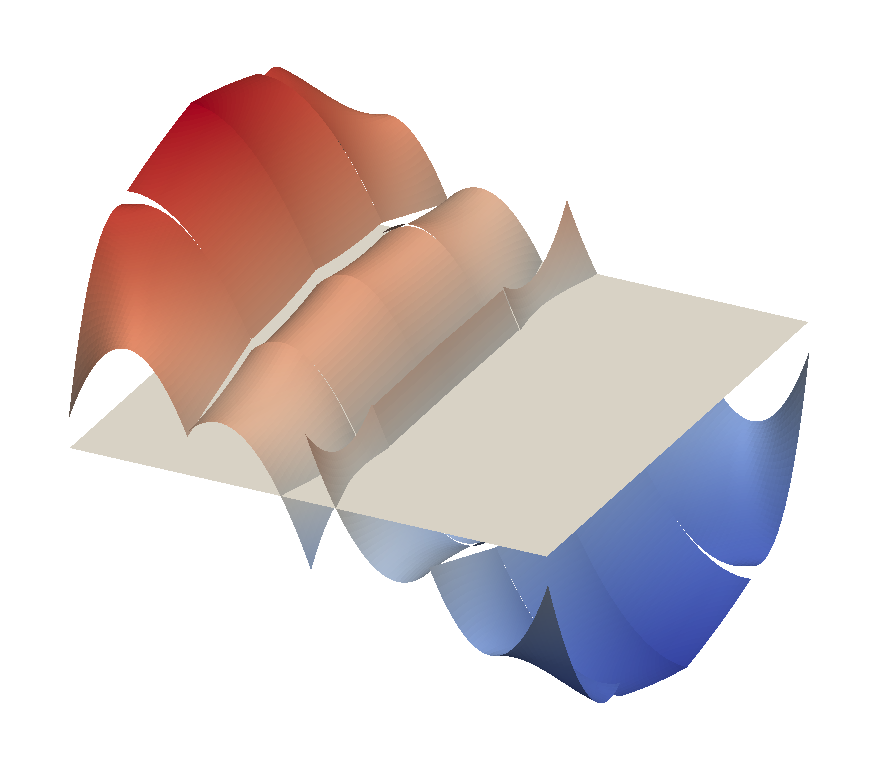}
    \includegraphics[width=0.2\textwidth]{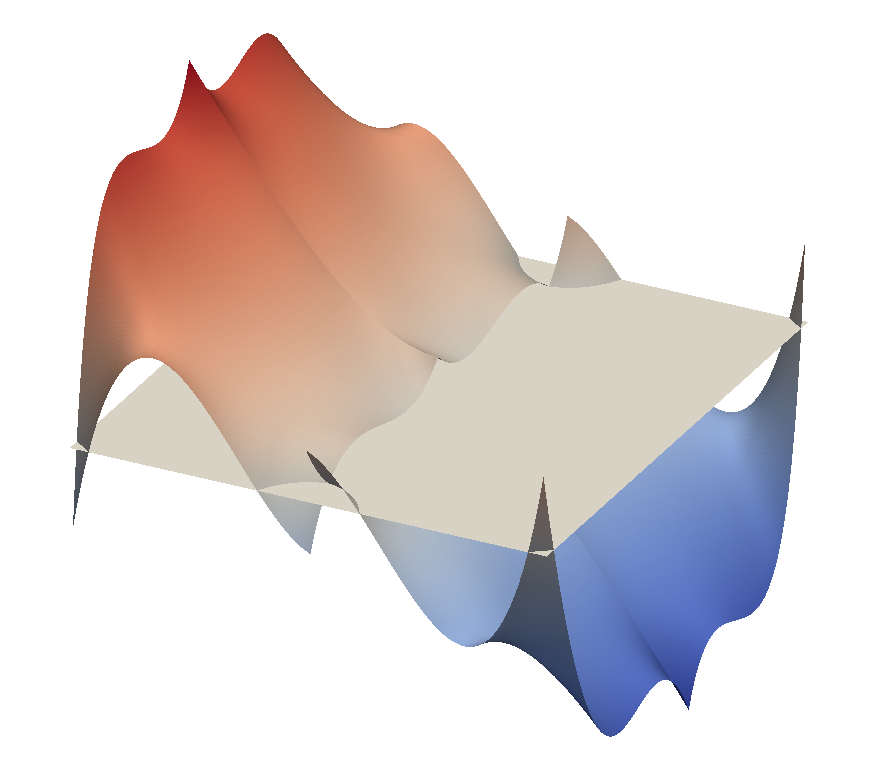}
    \includegraphics[width=0.2\textwidth]{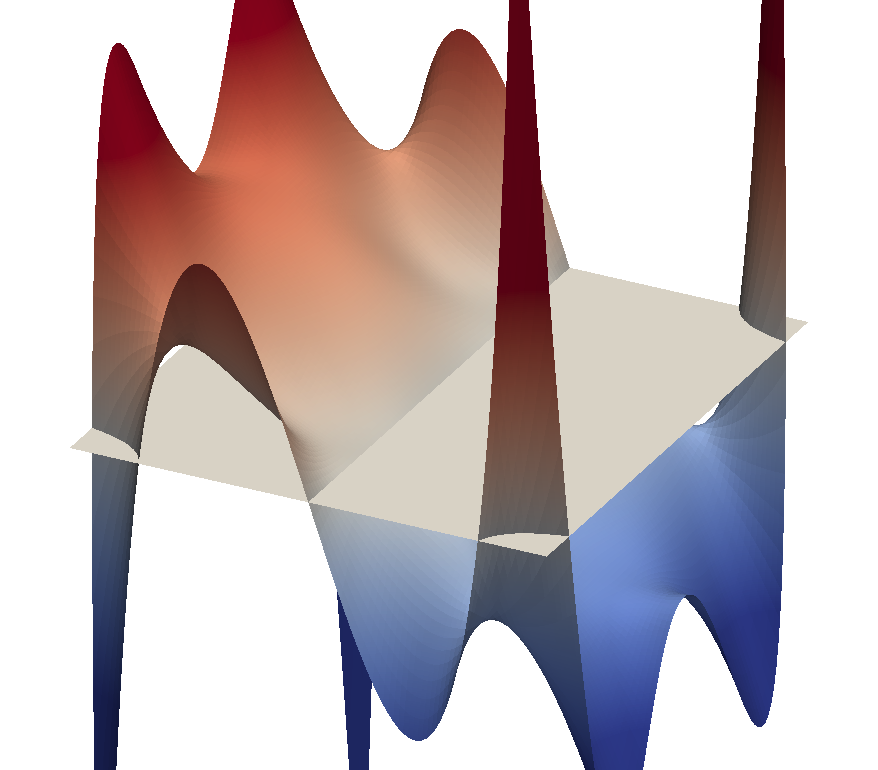}
}   
 
\subfigure[$\Delta t/T = 10^{-3}$, mixed-order formulation~$(k_u,k_p)=(k+1,k)$, refine level~$l=4-k$]{
    \includegraphics[width=0.2\textwidth]{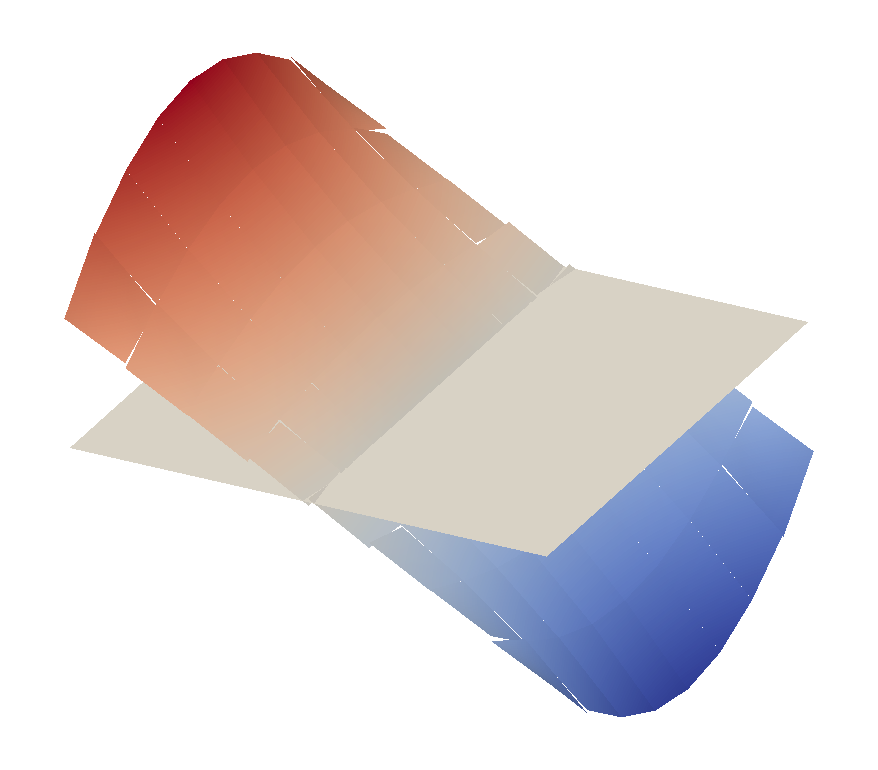}
    \includegraphics[width=0.2\textwidth]{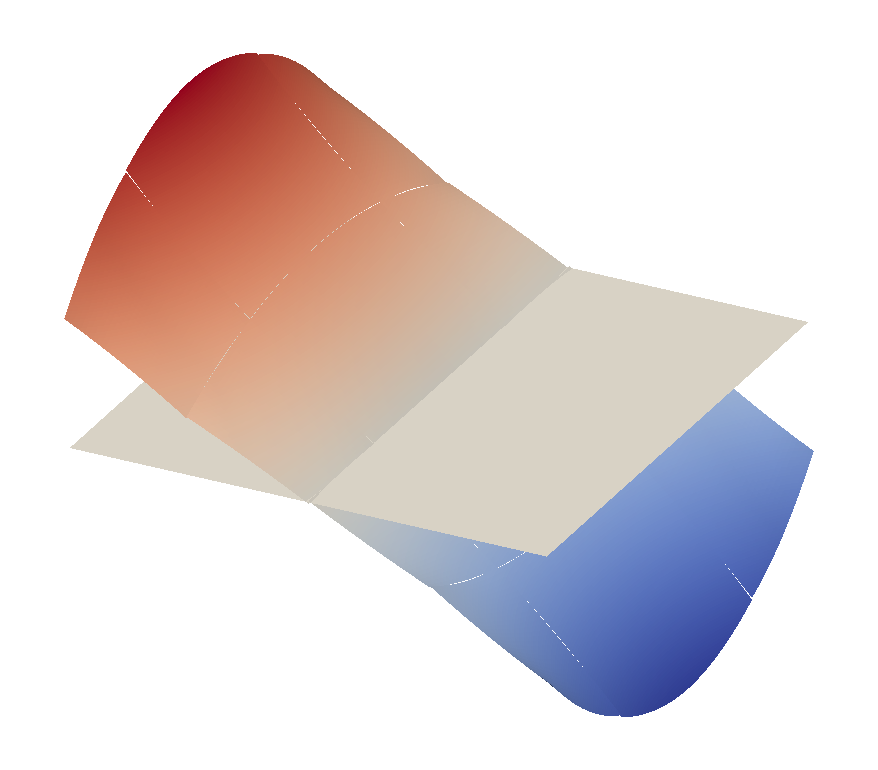}
    \includegraphics[width=0.2\textwidth]{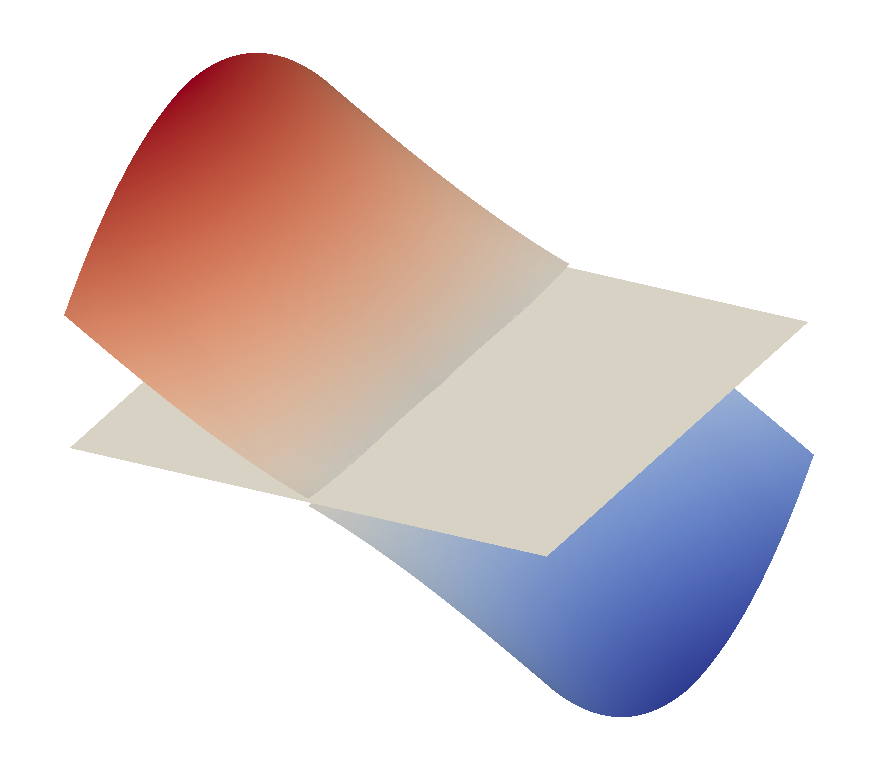}
    \includegraphics[width=0.2\textwidth]{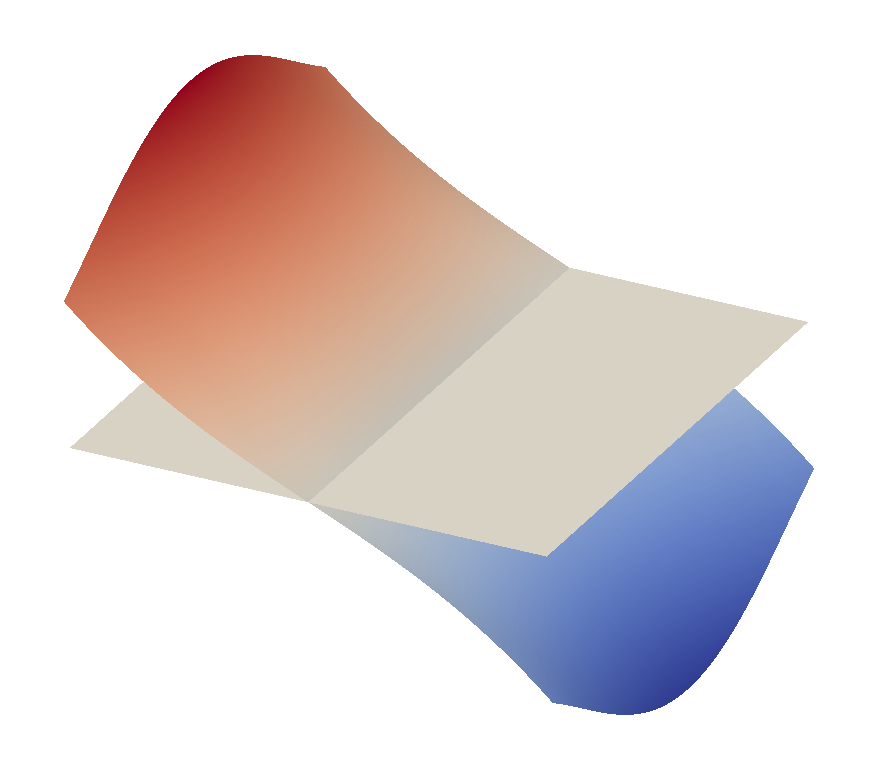}
}
\caption{Inf--sup stability of dual splitting scheme: pressure solution at time~$t=T$ for equal-order formulation and mixed-order formulation and different time step sizes. The parameter~$k=1,...,4$ varies from~$k=1$ to~$k=4$ from left to right.}
\label{fig:inf_sup_instabilities_dual_splitting}
\end{figure}

\begin{figure}[t]
 \centering 
 \subfigure[$\Delta t/T = 10^{-1}$, equal-order formulation~$(k_u,k_p)=(k,k)$, refine level~$l=4-k$]{
	\includegraphics[width=0.2\textwidth]{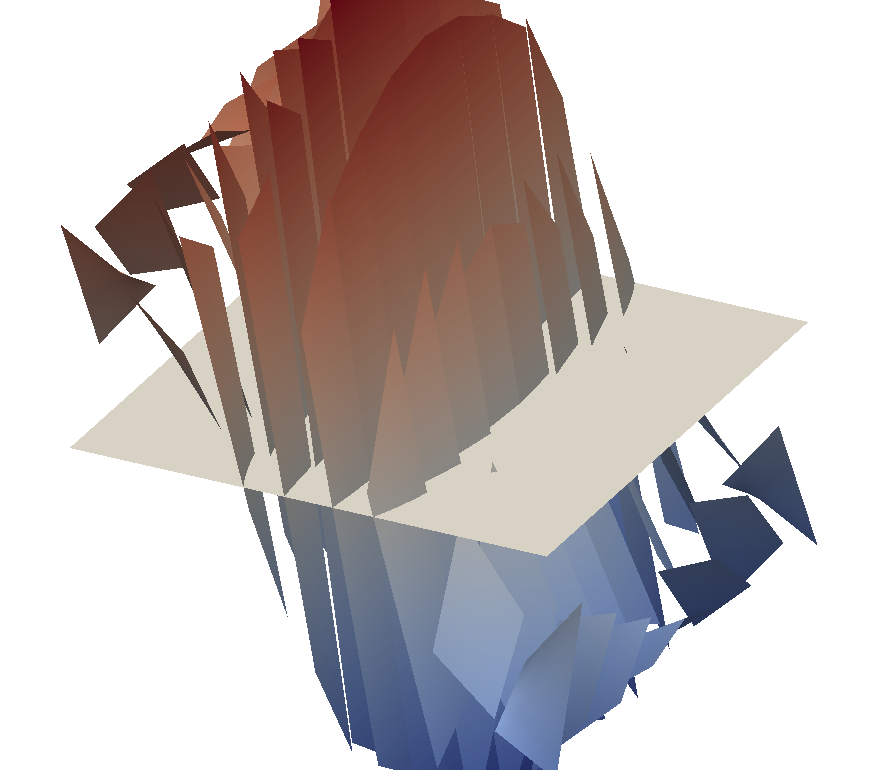}
    \includegraphics[width=0.2\textwidth]{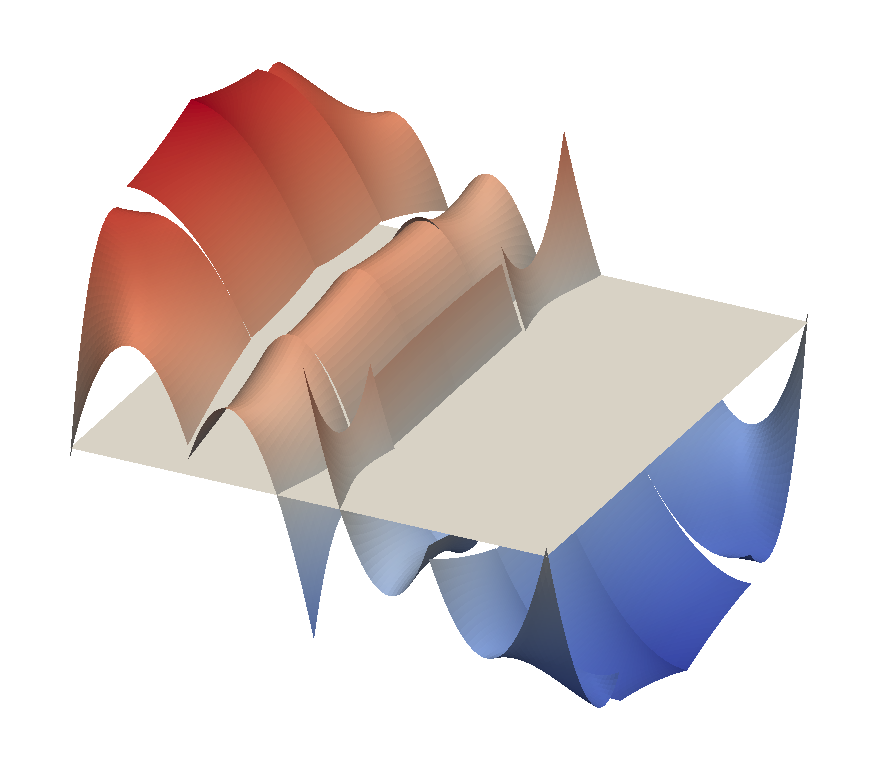}
    \includegraphics[width=0.2\textwidth]{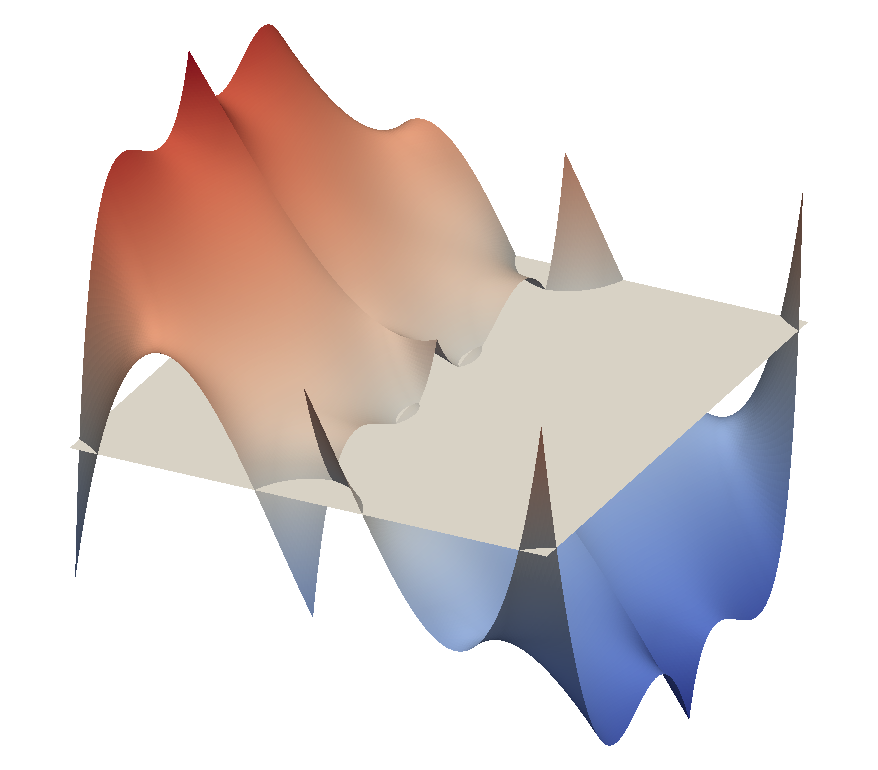}
    \includegraphics[width=0.2\textwidth]{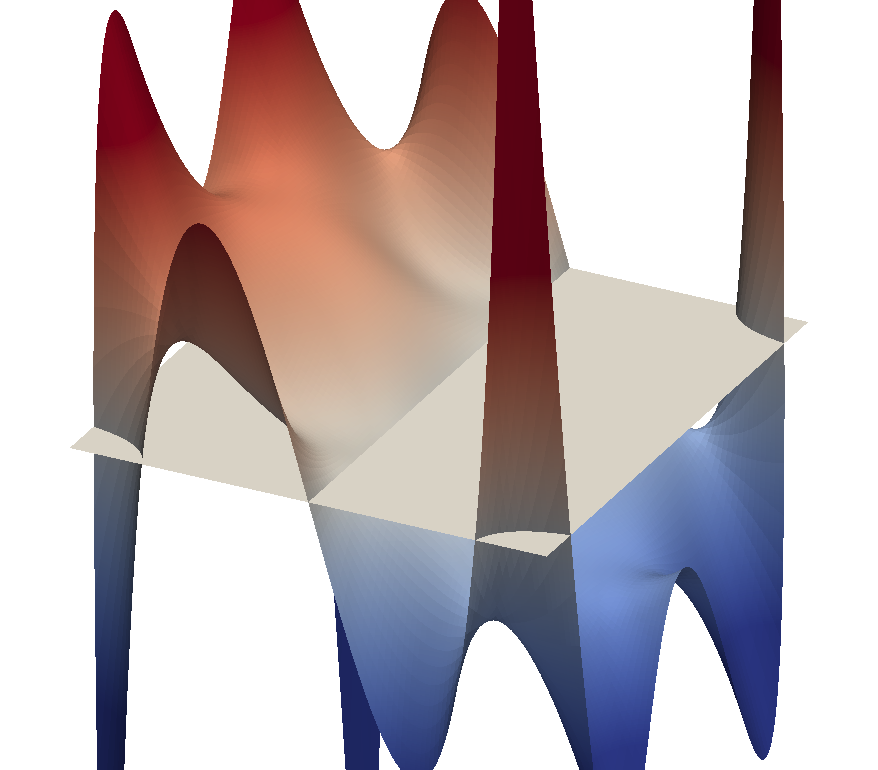}
}
    
\subfigure[$\Delta t/T = 10^{-3}$, equal-order formulation~$(k_u,k_p)=(k,k)$, refine level~$l=4-k$ \label{subfig:inf_sup_pc_dt_T_1e-3_equal_order}]{
    \includegraphics[width=0.2\textwidth]{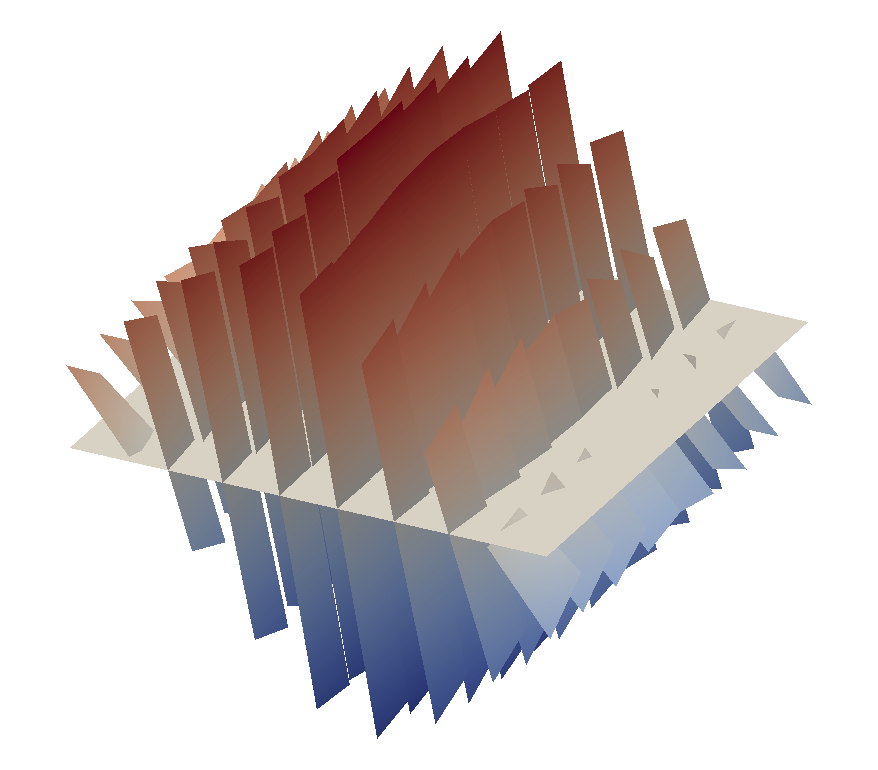}
    \includegraphics[width=0.2\textwidth]{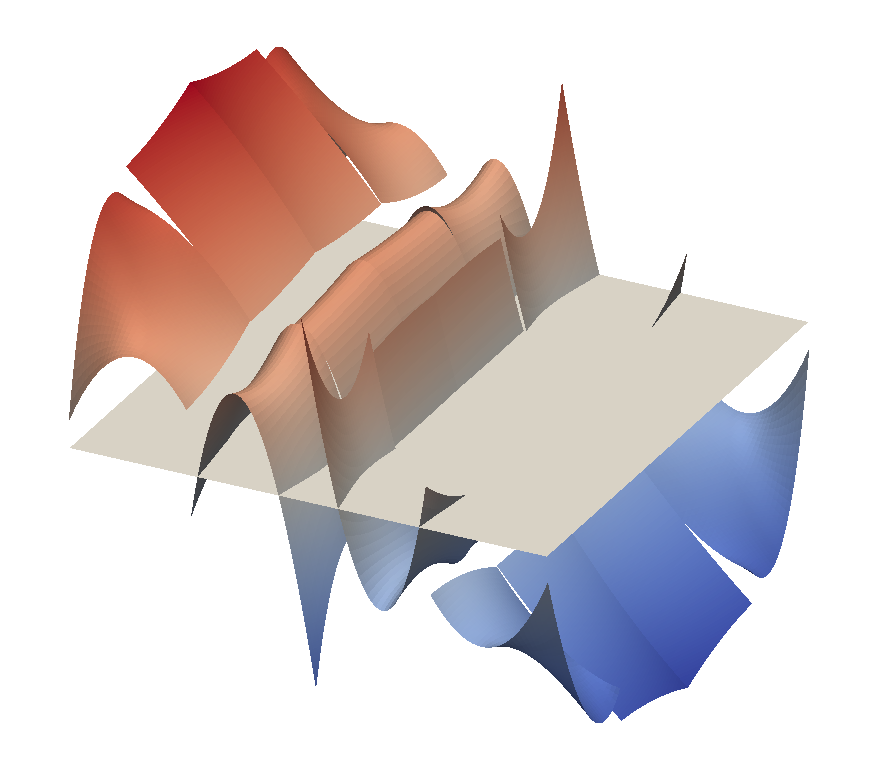}
    \includegraphics[width=0.2\textwidth]{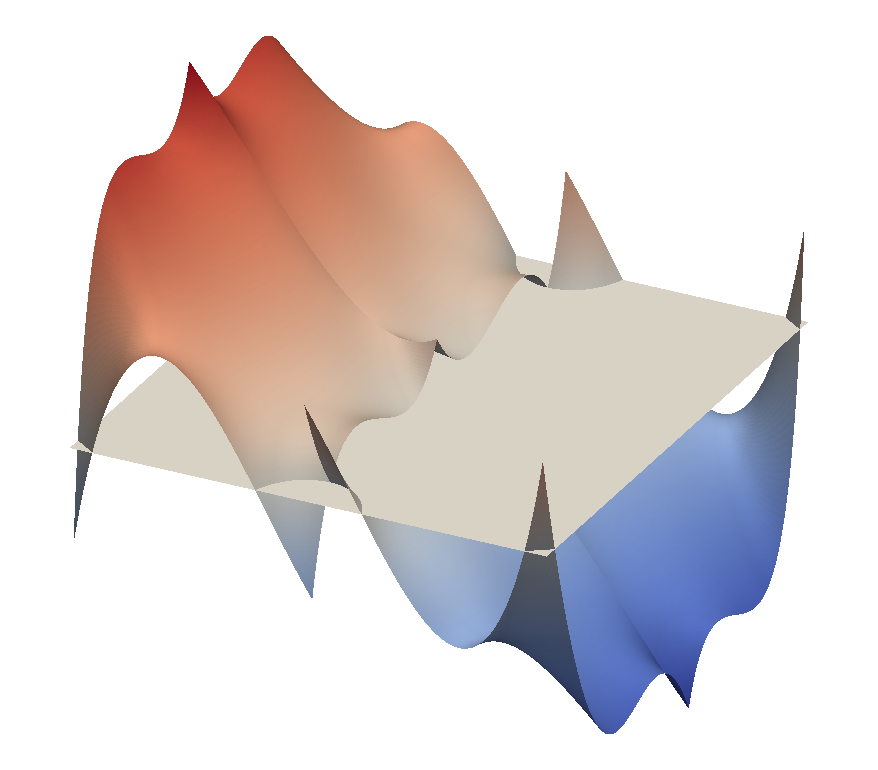}
    \includegraphics[width=0.2\textwidth]{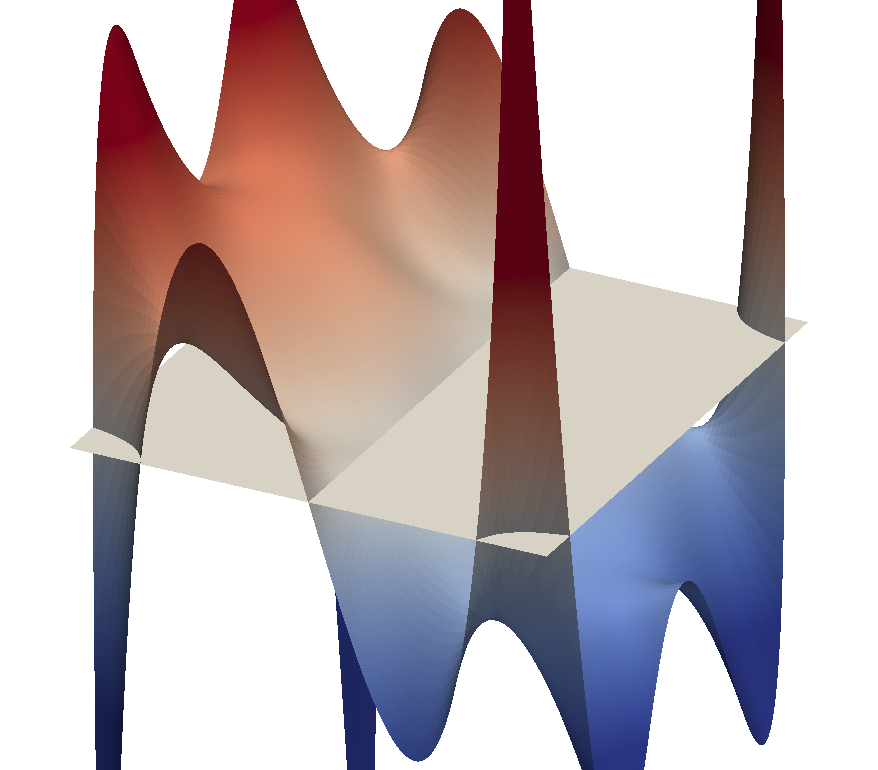}
}   
 
\subfigure[$\Delta t/T = 10^{-3}$, mixed-order formulation~$(k_u,k_p)=(k+1,k)$, refine level~$l=4-k$]{
    \includegraphics[width=0.2\textwidth]{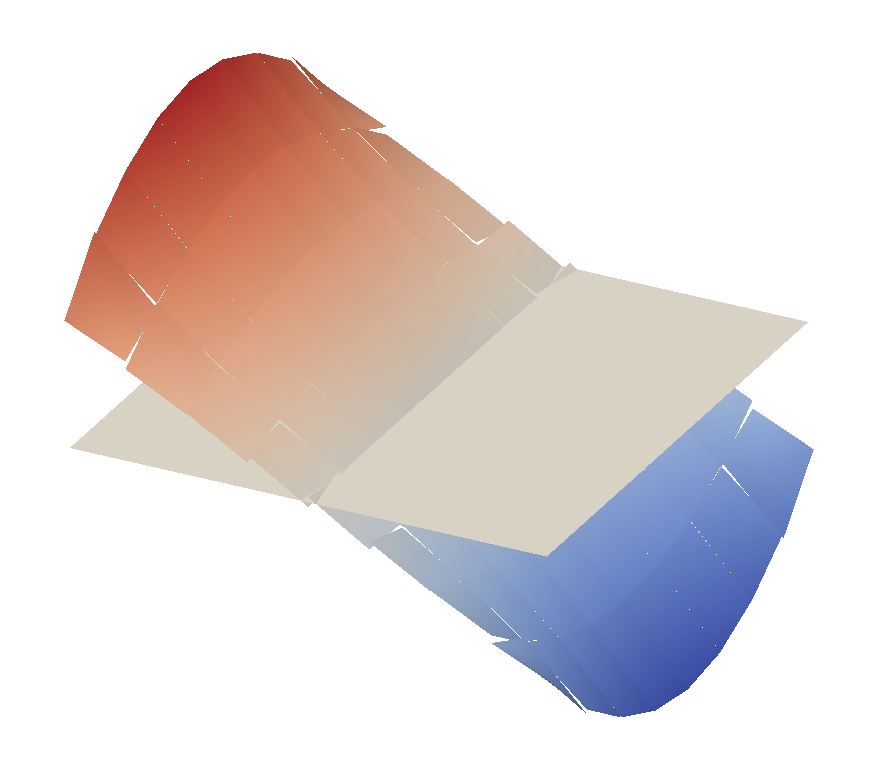}
    \includegraphics[width=0.2\textwidth]{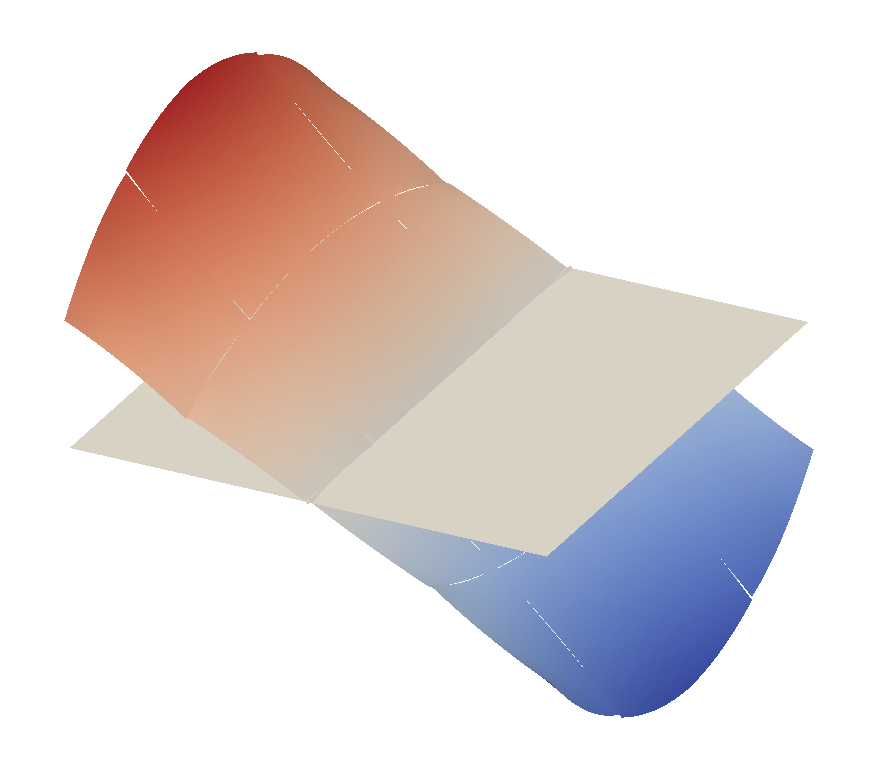}
    \includegraphics[width=0.2\textwidth]{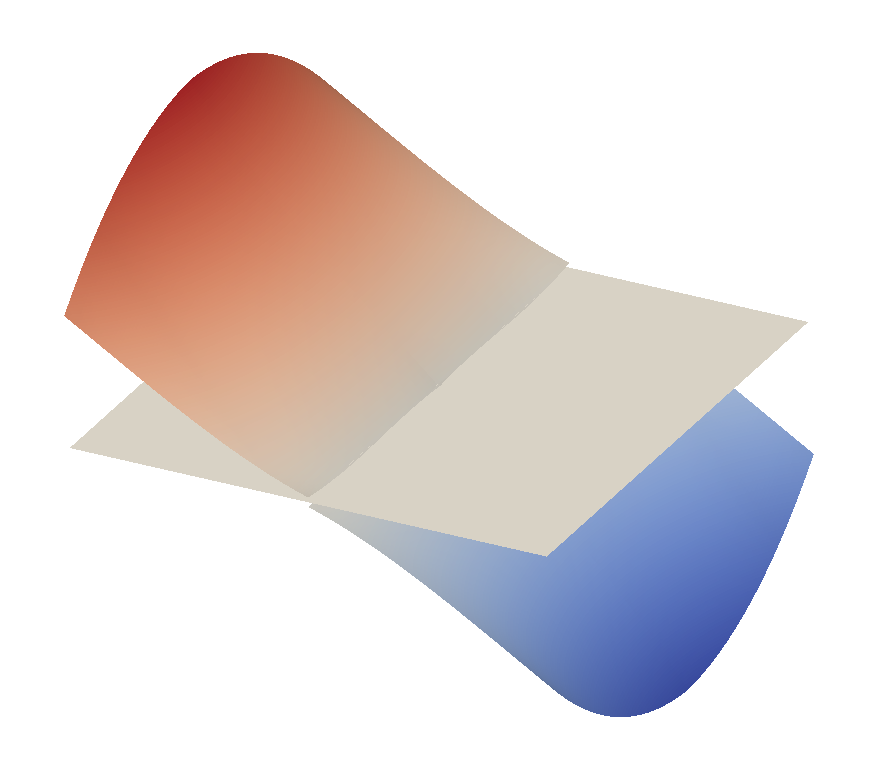}
    \includegraphics[width=0.2\textwidth]{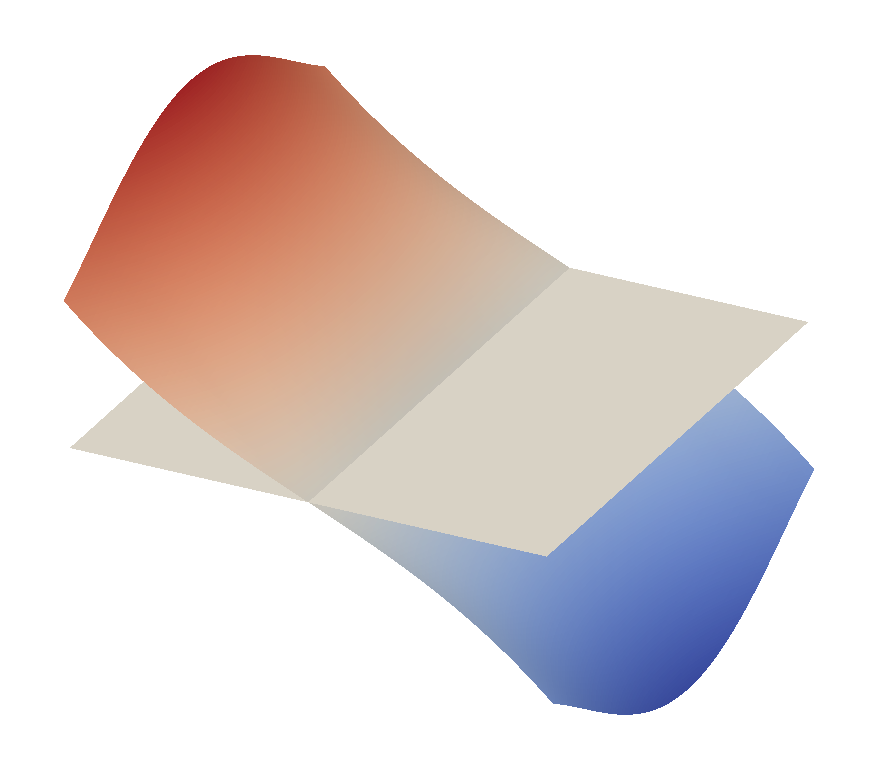}
}
\caption{Inf--sup stability of incremental pressure-correction scheme: pressure solution at time~$t=T$ for equal-order formulation and mixed-order formulation and different time step sizes. The parameter~$k=1,...,4$ varies from~$k=1$ to~$k=4$ from left to right. Due to large pressure oscillations the pressure solution for~$k=1$ (left picture) in subfigure~\ref{subfig:inf_sup_pc_dt_T_1e-3_equal_order} is scaled by a factor of~$0.2$ compared to all other pressure plots.}
\label{fig:inf_sup_instabilities_pressure_correction}
\end{figure}

As a means of verifying the above considerations on inf--sup stabilization, we perform spatial convergence tests for both equal-order and mixed-order polynomials. Moreover, the results for the dual splitting scheme and pressure-correction scheme are compared to the results obtained for the coupled solution approach. Table~\ref{SpatialConvergenceStokes} provides information on spatial convergence tests for equal-order polynomials~$(k_u,k_p)=(k,k)$ and mixed-order polynomials~$(k_u,k_p)=(k,k-1)$ for~$k=2,3,4,5$. The refine level is varied from~$l=1$ to~$l=4$ for all polynomial degrees. With respect to the temporal discretization we use BDF2 time integration and a small time step size of~$\Delta t /T = 10^{-4}$ for all simulations to ensure that the spatial discretization error is dominant. For the pressure-correction scheme the incremental formulation with~$J_p=1$ and the rotational formulation is used, resulting in the most accurate scheme as shown in Section~\ref{TemporalConvergenceStokes}.

Using equal-order polynomials for velocity and pressure, both the coupled solution approach and the pressure-correction scheme yield suboptimal rates of convergence. With respect to the approximation of the velocity solution, experimental convergence rates between~$k_u$ and~$k_u+1$ are obtained. For the pressure, convergence rates are between~$k_p$ and~$k_p+1$ for~$k_p=2,4$ but considerable lower than the theoretical values for polynomial degrees~$k_p=3,5$. The convergence behavior is significantly improved when using the dual splitting scheme: The velocity converges with rates between~$k_u+1/2$ and~$k_u+1$ and also for the pressure convergence rates between~$k_p+1/2$ and~$k_p+1$ are obtained for higher refine levels. These results are in agreement with the above theoretical considerations, equations~\eqref{InfSupDualSplitting} and~\eqref{InfSupPressureCorrection}. There is no stabilization term in case of the incremental pressure-correction scheme considered here similar to the coupled solution approach. However, the dual splitting scheme introduces an inf--sup stabilization term according to equation~\eqref{InfSupDualSplitting} where the impact of the stabilization increases for fine spatial resolutions due to the presence of second derivatives in the stabilization term. This could explain the improved convergence observed for the dual splitting scheme.

Using mixed-order polynomials, optimal rates of convergence of order~$h^{k_u+1}$ for the velocity and~$h^{k_p+1}=h^{k_u}$ for the pressure are obtained for all solution approaches. 
We assume that the occurrence of superconvergence effects of the pressure solution with rates of convergence higher than~$k_p+1$ in case of the mixed-order formulation is related to this specific flow problem. Consider also the spatial convergence results in Section~\ref{VortexProblem} where convergence rates of the pressure very close to the theoretical optimum~$h^{k_p+1}$ are obtained for the mixed-order formulation.

While velocity errors are comparable for the equal-order formulation and the mixed-order formulation, pressure errors are significantly smaller when using the mixed-order formulation instead of the equal-order formulation especially for the coupled solution approach and the pressure-correction scheme. This indicates the presence of spurious pressure oscillations in case of the equal-order formulation related to the inf--sup stability condition.

In addition to the spatial convergence tests, we consider varying time step sizes for given spatial resolutions. Focusing on the dual splitting scheme, Figure~\ref{fig:inf_sup_instabilities_dual_splitting} displays the pressure solution at final time~$t=T$ for equal-order and mixed-order polynomials of varying degree as well as for different time step sizes of~$\Delta t/T=10^{-1}$ and~$\Delta t/T=10^{-3}$. Since inf--sup instabilities are expected to be pronounced for coarse meshes, we consider comparably low spatial resolutions and simultaneously reduce the level of refinement when increasing the polynomial degree. For equal-order polynomials, the time step size has a huge influence on the pressure solution and artificial pressure modes show up for small~$\Delta t$. This is in agreement with equation~\eqref{InfSupDualSplitting} predicting that the stabilizing effect is related to~$\Delta t$ and will diminish when decreasing the time step size. For the mixed-order formulation the pressure solution is smooth and results for~$\Delta t/T=10^{-1}$ and~$\Delta t/T=10^{-3}$ are indistinguishable. Hence, we only present results for the smaller time step size which is the critical one in this respect. We performed the same simulations for the nonincremental pressure-correction scheme in standard formulation. For the equal-order formulation and the two time step sizes~$\Delta t/T=10^{-1}$ and~$\Delta t/T=10^{-3}$, very similar results are obtained as for the dual-splitting scheme in terms of spurious pressure oscillations which is in line with equations~\eqref{InfSupDualSplitting} and~\eqref{InfSupPressureCorrection}.

Figure~\ref{fig:inf_sup_instabilities_pressure_correction} shows results for the same stability experiment using the incremental pressure-correction scheme in rotational formulation. In contrast to the dual splitting scheme, inf--sup instabilities also occur for very large time step sizes when using equal-order polynomials. Again, no oscillations occur for the mixed-order formulation. We also note that the coupled solution approach yields results similar to those for the incremental pressure-correction scheme in rotational form. These results can be seen as a numerical verification of equation~\eqref{InfSupPressureCorrection} stating that the incremental pressure-correction scheme and the coupled solution approach behave similarly in terms of inf--sup stabilization.

\subsection{Unsteady Navier--Stokes equations: Vortex problem}\label{VortexProblem}
In order to verify the implementation of the different solution strategies for the incompressible Navier--Stokes equations including the convective term and to demonstrate optimal rates of convergence with respect to the temporal and the spatial discretization, we consider the vortex problem analyzed in~\cite{Hesthaven07}. This test case is an analytical solution of the unsteady incompressible Navier--Stokes equations in two dimensions for~$\bm{f}=\bm{0}$
\begin{align}
\begin{split}
\bm{u}(\bm{x},t) &=  \begin{pmatrix}
-\sin(2\pi x_2)\\
+\sin(2\pi x_1)
\end{pmatrix}
\exp\left(-4\nu\pi^2 t\right)\; ,\\
p(\bm{x},t) &=  -\cos(2\pi x_1)\cos(2\pi x_2)\exp\left(-8\nu \pi^2 t\right)\; .
\end{split}\label{AnalyticalSolutionVortex}
\end{align}
The viscosity is set to~$\nu=0.025$. The domain~$\Omega=[-L/2,L/2]^2$ is a square of length~$L=1$ and the simulations are performed for the time interval~$0\leq t\leq T=1$. On domain boundaries, Dirichlet boundary conditions are prescribed at the inflow part of the boundary and Neumann boundary conditions at the outflow part so that the coordinate axes split each of the four sides of the rectangular domain into a Dirichlet part and a Neumann part, see also~\cite{Hesthaven07}. Initial conditions as well as the solution at previous instants of time~$t_{n-J+1},...,t_{n-1}$ required by the BDF scheme and extrapolation scheme for~$J>1$ are obtained by interpolation of the analytical solution. The velocity Dirichlet boundary condition~$\bm{g}_u$, the time derivative term~$\partial \bm{g}_u/\partial t$, the velocity gradient in normal direction~$\bm{h}_u/\nu$, and the pressure Dirichlet boundary condition~$g_p$ are derived from the analytical solution. In case of the coupled solution approach the Neumann boundary condition is then given as~$\bm{h} = \bm{h}_u- g_p \bm{n}$. Since the velocity boundary conditions~$\bm{g}_u$ and~$\bm{h}_u$ and the pressure boundary condition~$g_p$ are nontrivial and time-dependent, this flow problem is an appropriate test case to verify the temporal accuracy of the different solution approaches with their respective boundary conditions. As for the Stokes flow problem, a uniform Cartesian grid is used for discretization in space where the element length is~$h=L/2^l$ with the refinement level~$l$. For this flow problem, we only consider mixed-order formulations for reasons explained above.

\subsubsection{Temporal convergence test}
We perform temporal convergence tests for BDF schemes of order 1 and 2. For the pressure-correction scheme we use an order of extrapolation of the pressure gradient term of~$J_p=J-1$ and analyze both the standard formulation and the rotational formulation. A high spatial resolution (refine level~$l=3$ and polynomial degrees~$k_u=8$,~$k_p=7$) is used so that the spatial discretization error is negligible compared to the temporal discretization error.

\begin{figure}[t]
 \centering 
	\includegraphics[width=1.0\textwidth]{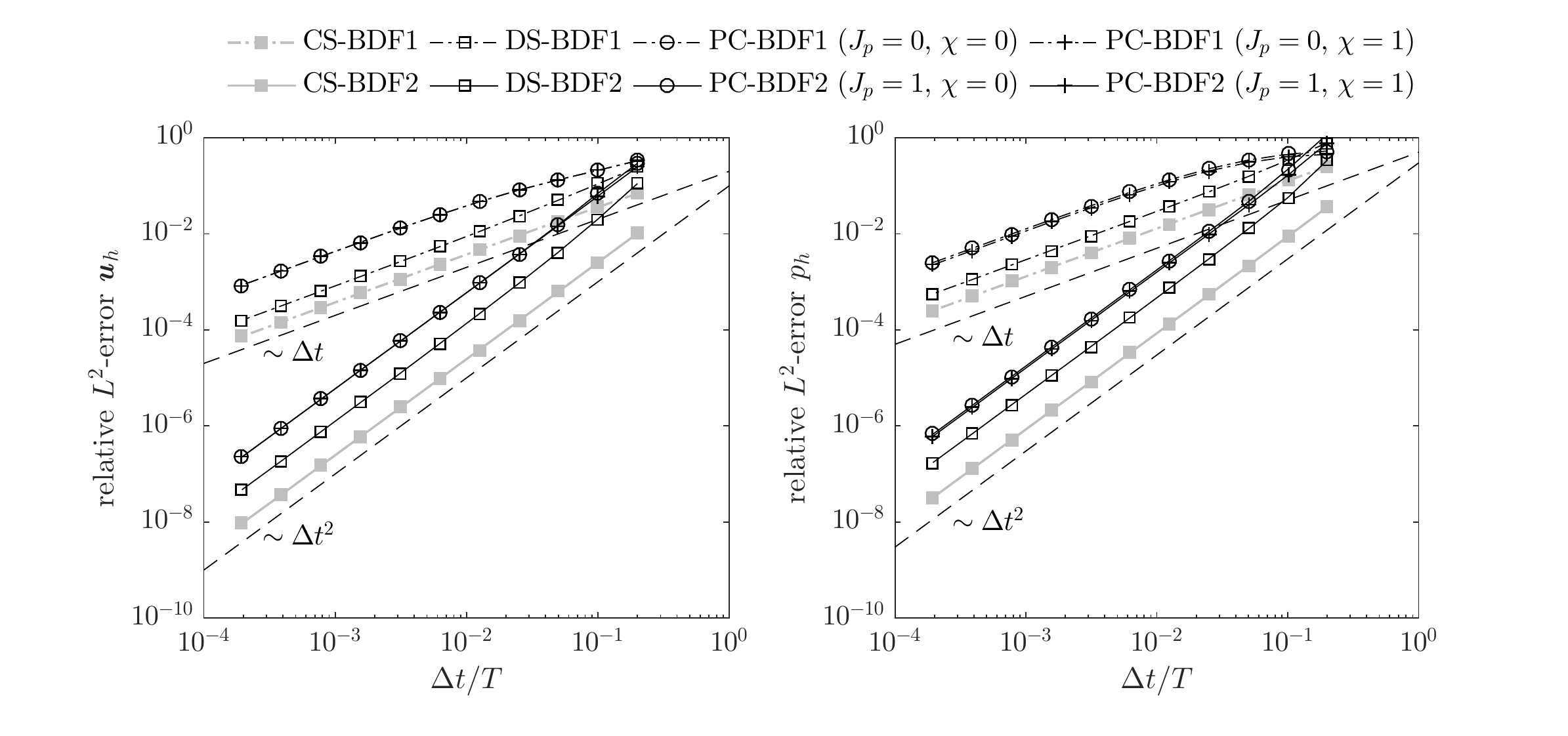}
\caption{Vortex problem: temporal convergence tests for coupled solver (CS), dual splitting scheme (DS), and pressure-correction scheme (PC) in standard form~($\chi=0$) and rotational form~($\chi=1$) for BDF schemes of order 1 and 2. The spatial resolution is~$l=3$ and~$(k_u,k_p)=(8,7)$.}
\label{fig:temporal_convergence_vortex}
\end{figure}

Results of the temporal convergence test are presented in Figure~\ref{fig:temporal_convergence_vortex}. All methods show optimal rates of convergence of order~$\Delta t^1$ for the BDF1 scheme and~$\Delta t^2$ for the BDF2 scheme. In terms of absolute errors, the coupled solution approach is the most accurate method while the pressure-correction scheme is the least accurate one for this problem. For the pressure-correction scheme, the errors are almost the same for the standard formulation and the rotational formulation which might be explained by the fact that the pressure gradient in normal direction,~$\nabla p \cdot \bm{n}$, is zero on~$\Gamma^{\rm{D}}$ according to the analytical solution~\eqref{AnalyticalSolutionVortex}.

\begin{remark}
A more detailed analysis of the absolute errors of the different solution strategies reveals that the increased error of the dual splitting scheme as compared to the coupled solver is due to the explicit treatment of the convective term, i.e., using an extrapolation of the convective term of order~$J$ in case of the coupled solution approach leads to results that are very close to those obtained for the dual splitting scheme. Interestingly, using an explicit treatment of the convective term slightly reduces the errors in case of the pressure-correction scheme. However, the error is still a factor of approximately~$2$ larger compared to the coupled solution approach or dual splitting scheme in that case.
\end{remark}

\subsubsection{Spatial convergence test}
\begin{figure}[t]
 \centering 
	\includegraphics[width=1.0\textwidth]{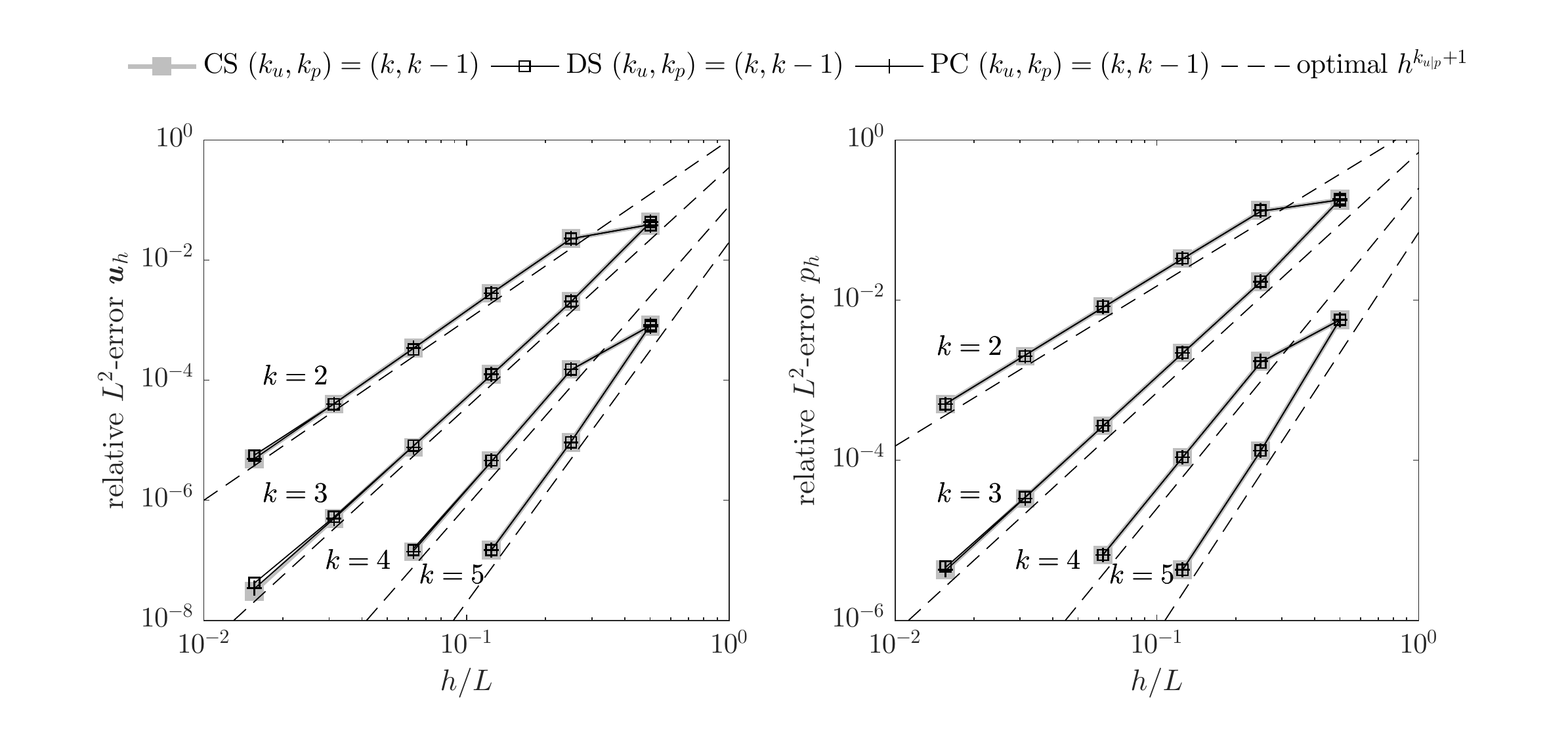}
\caption{Vortex problem: spatial convergence tests for coupled solver (CS), dual splitting scheme (DS), and incremental ($J_p=1$) pressure-correction scheme (PC) in rotational form for mixed-order formulations~$(k_u,k_p)=(k,k-1)$ where~$k=2,...,5$. For discretization in time the BDF2 scheme and a fix time step size of~$\Delta t/T = 5\cdot 10^{-5}$ is used for all simulations.}
\label{fig:spatial_convergence_vortex}
\end{figure}

Results of the spatial convergence tests are shown in Figure~\ref{fig:spatial_convergence_vortex} for mixed-order formulations~$(k_u,k_p)=(k,k-1)$ and polynomial degrees in the range~$k=2,...,5$. To ensure a small temporal discretization error the BDF2 scheme is used and a fix time step size of~$\Delta t /T = 5\cdot 10^{-5}$. For the pressure-correction scheme, the incremental formulation with~$J_p=1$ and the rotational formulation is used. Deviating from the solver tolerances specified in Section~\ref{Implementation}, the absolute tolerance of the Newton solver for the momentum equation of the pressure-correction scheme is set to~$10^{-10}$ to obtain convergence for all spatial resolutions.

For all solution approaches, experimental rates of convergence very close to the optimal rates of convergence of order~$h^{k_u+1}$ for the velocity and~$h^{k_p+1}$ for the pressure are obtained. Moreover, the errors are virtually the same for all solution techniques. 

\subsection{Unsteady Navier--Stokes equations: Laminar flow around a cylinder}
In order to demonstrate the geometric flexibility of the present Navier--Stokes solvers and to analyze the efficiency of high-order polynomial spaces for the approximation of velocity and pressure on complex domains with curved boundaries, we consider laminar flow around a cylinder with unsteady vortex shedding. This test case has been proposed in the 1990s by~\cite{Schafer1996} as a benchmark problem and has found widespread use in terms of the verification of incompressible Navier--Stokes solvers. In the present work we focus on the two-dimensional, unsteady test case named 2D-3 for which an accurate reference solution is available~\cite{John2004}.

\begin{figure}[t]
 \centering 
	\includegraphics[width=0.6\textwidth]{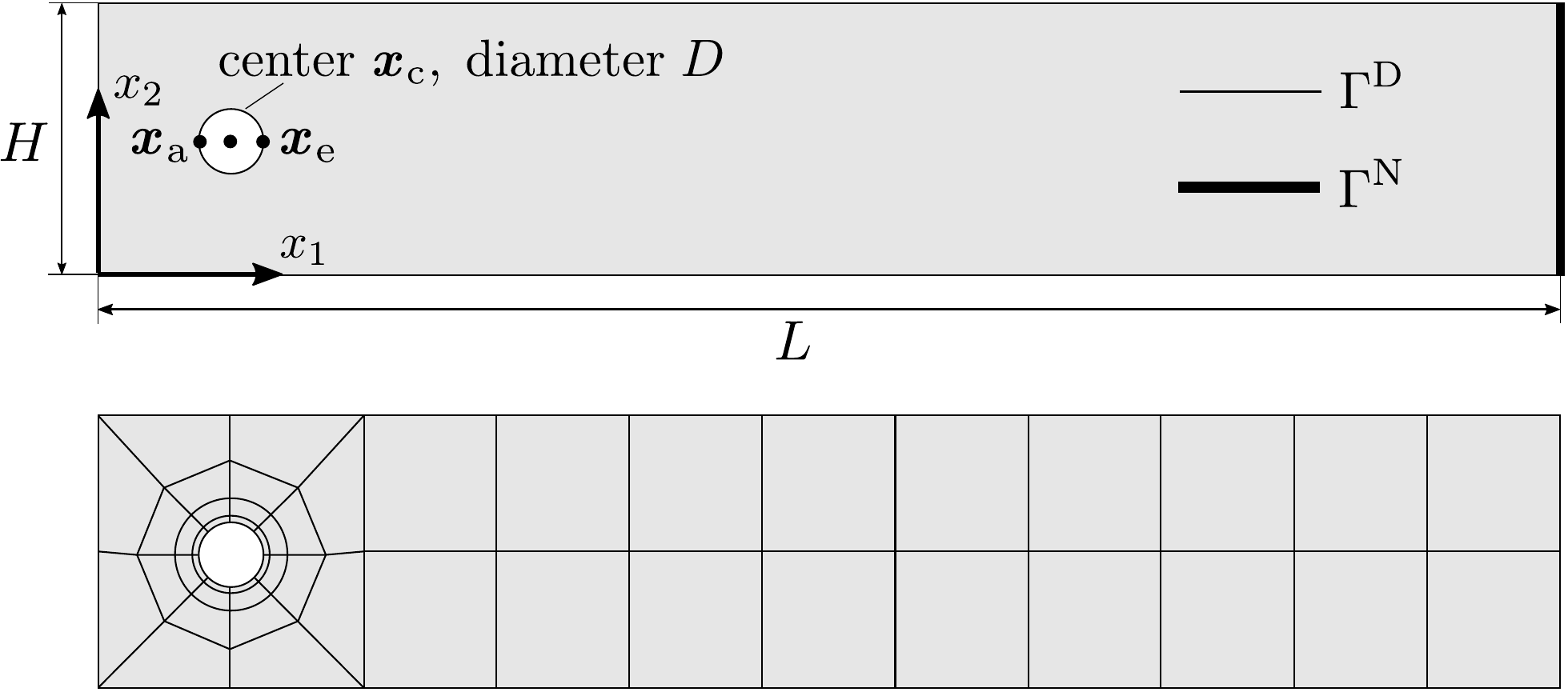}
\caption{Laminar flow around cylinder: geometry and boundary conditions according to the benchmark~\cite{Schafer1996} as well as coarsest mesh corresponding to refine level~$l=0$.}
\label{fig:FlowAroundCylinder}
\end{figure}

The geometry is displayed in Figure~\ref{fig:FlowAroundCylinder}. The cylinder with center~$(x_{1,\rm{c}},x_{2,\rm{c}})=(0.2,0.2)$ and diameter~$D=0.1$ is located slightly asymmetrically in a rectangular channel of length~$L=2.2$ and height~$H=0.41$. The inflow boundary (left boundary), the channel walls (upper and lower boundary) and the cylinder surface are treated as Dirichlet boundaries. At the inflow boundary a parabolic velocity profile is prescribed~\cite{Schafer1996}
\begin{align}
g_{u_1}(x_1=0,x_2,t) = U_{\rm{m}} \frac{4 x_2 (H-x_2)}{H^2}\sin(\pi t/T)\; , \; g_{u_2}(x_1=0,x_2,t)=0 \; ,
\end{align}
where the time interval is~$0\leq t \leq T=8$. The Reynolds number~$\mathrm{Re} = \bar{U} D/\nu$ is defined using the mean inflow velocity~$\bar{U} = 2U_{\rm{m}}/3$ and the cylinder diameter~$D$. The maximum inflow velocity is given as~$U_{\rm{m}}=1.5$ and the viscosity is~$\nu=10^{-3}$ so that the Reynolds number reaches a maximum value of~$\mathrm{Re}_{\max}=100$ at time~$t=T/2$.
On the cylinder surface and the channel walls no slip boundary conditions are prescribed for the velocity,~$\bm{g}_u=\bm{0}$. The outflow boundary (right boundary) is treated as a Neumann boundary where the benchmark itself does not define a specific boundary condition on~$\Gamma^{\rm{N}}$. For the decoupled solution approaches we prescribe~$\bm{h}_u=\bm{0}$ and~$g_p=0$ and for the coupled solution approach~$\bm{h}=\bm{0}$.

The accuracy of the numerical solution is evaluated by calculating the maximum drag coefficient~$c_{\rm{D},\max}$, the maximum lift coefficient~$c_{\rm{L},\max}$, and the pressure difference~$\Delta p(t=T)=p(\bm{x}_{\rm{a}},t=T)-p(\bm{x}_{\rm{e}},t=T)$ between the front and the back of the cylinder at final time~$t=T$, where~$\bm{x}_{\rm{a}} = (x_{1,\rm{c}}-D/2,x_{2,\rm{c}})^T$ and~$\bm{x}_{\rm{e}} = (x_{1,\rm{c}}+D/2,x_{2,\rm{c}})^T$. The drag coefficient~$c_{\rm{D}} = F_1/(\rho\bar{U}^2 D/2)$ and the lift coefficient~$c_{\rm{L}} = F_2/(\rho\bar{U}^2 D/2)$ are obtained by calculating the force vector~$\bm{F} = (F_1,F_2)^T= - \rho \int_A \left(-p\bm{I} + \nu \left(\Grad{\bm{u}}+\Grad{\bm{u}}^T\right)\right)\cdot \bm{n} \mathrm{d}A$ acting on the cylinder, where~$A$ denotes the cylinder surface and~$\bm{n}$ the outward pointing normal vector of the computational domain. Reference solutions of these quantities are listed in Table~\ref{FlowPastCylinder2D-3_Reference}.

\begin{table}[!h]
\caption{Laminar flow around cylinder: reference results for test case 2D-3}\label{FlowPastCylinder2D-3_Reference}
\renewcommand{\arraystretch}{1.1}
\begin{center}
\begin{tabular}{ccccc}
\hline  
reference & $c_{\mathrm{D},\max}$ &  $c_{\mathrm{L},\max}$ &$\Delta p (t=T)$\\ 
\hline
Sch{\"a}fer et al.~\cite{Schafer1996}  & $2.95 \pm 2\cdot 10^{-2}$ &  $0.48 \pm 1\cdot 10^{-2}$ & $-0.11 \pm 5\cdot 10^{-3}$ \\     
John~\cite{John2004}  & $2.950921575\pm 5\cdot 10^{-7}$ &  $0.47795\pm 1\cdot 10^{-4}$ &  $-0.1116\pm 1\cdot 10^{-4}$ \\
present ($l=3, (k_u,k_p)=(10,9)$) & $2.95091839$ & $0.47788776$ & $-0.11161590$\\
\hline
\end{tabular}

\end{center}
\renewcommand{\arraystretch}{1}
\end{table}

The mesh is visualized in Figure~\ref{fig:FlowAroundCylinder} for the coarsest refine level~$l=0$ which consists of~$N_{\text{el},l=0}=50$ quadrilateral elements. Finer meshes are obtained by uniform refinement so that the number of elements on level~$l$ is~$N_{\text{el},l}=N_{\text{el},l=0} (2^d)^l$. The total number of degrees of freedom is~$N_{\text{dofs}}=N_{\text{el},l}N_{\text{dofs},\text{el}}$, where the number of unknowns per element is~$N_{\text{dofs},\text{el}}=d(k_u+1)^d + (k_p+1)^d$. In order to accurately resolve the flow near the cylinder, the mesh is refined towards the cylinder and an isoparametric mapping is used for an improved approximation of the curved cylinder boundary. In this respect, the first two layers of elements around the cylinder are subject to a cylindrical manifold description to enable high-order accuracy. For the third layer of cells we implemented a volume manifold description allowing to prescribe a cylindrical manifold for one of the four faces of the quadrilateral element with straight edges on the other faces.

In Table~\ref{SpatialConvergenceFlowPastCylinder2D-3} we present spatial convergence results for mixed order polynomials~$(k_u,k_p)=(k,k-1)$ with polynomial degrees in the range~$k=2,4,6,8,10$. For this analysis, the high-order dual splitting scheme with~$J=2$ is used but very similar results are obtained for the pressure-correction scheme and the monolithic solver. Since there is no analytical solution for this test case, the time integration scheme is started using~$J=1$ in the first time step. The time step size is chosen small enough so that the overall error is dominated by the spatial discretization error.
\begin{table}[!h]
\caption{Laminar flow around cylinder: spatial convergence results for test case 2D-3 using the second order accurate dual splitting scheme.}\label{SpatialConvergenceFlowPastCylinder2D-3}
\renewcommand{\arraystretch}{1.1}
\begin{scriptsize}
\begin{center}
\begin{tabular}{ccccccccccc}
\hline
& & & \multicolumn{2}{c}{$c_{\mathrm{D},\max}$} & &\multicolumn{2}{c}{$c_{\mathrm{L},\max}$} & &\multicolumn{2}{c}{$\Delta p (t=T)$}\\ 
\cline{4-5} \cline{7-8} \cline{10-11} $(k_u,k_p)$ & $l$ & $N_{\text{dofs}}$ & value & relative error & & value & relative error & & value & relative error\\ 
\hline 
$(2,1)$ & 0 &    1100 & 3.38467316 & 1.47E--001 & & 2.39280367 & 4.01E+000  & & --0.09472515 & 1.51E--001\\
        & 1 &    4400 & 2.85701971 & 3.18E--002 & & 0.53988096 & 1.30E--001 & & --0.10305096 & 7.67E--002\\
        & 2 &   17600 & 2.94428788 & 2.25E--003 & & 0.47625908 & 3.41E--003 & & --0.09968075 & 1.07E--001\smallskip\\
        
$(4,3)$ & 0 &    3300 & 3.03777831 & 2.94E--002 & & 0.44023122 & 7.88E--002 & & --0.11266542 & 9.40E--003\\
        & 1 &   13200 & 2.94992818 & 3.36E--004 & & 0.49730575 & 4.06E--002 & & --0.10720962 & 3.95E--002\\
        & 2 &   52800 & 2.95100482 & 2.93E--005 & & 0.47746288 & 8.89E--004 & & --0.11169743 & 7.30E--004\smallskip\\
        
$(6,5)$ & 0 &    6700 & 2.94279309 & 2.75E--003 & & 0.53034748 & 1.10E--001 & & --0.10402386 & 6.80E--002\\
        & 1 &   26800 & 2.95097835 & 2.03E--005 & & 0.47861466 & 1.52E--003 & & --0.11113738 & 4.29E--003\\
        & 2 &  107200 & 2.95091829 & 3.39E--008 & & 0.47787797 & 2.05E--005 & & --0.11161756 & 1.49E--005\smallskip\\
        
$(8,7)$ & 0 &   11300 & 2.95252064 & 5.43E--004 & & 0.49339221 & 3.24E--002 & & --0.10378121 & 7.02E--002\\
        & 1 &   45200 & 2.95090792 & 3.55E--006 & & 0.47744518 & 9.26E--004 & & --0.11159388 & 1.97E--004\\
        & 2 &  180800 & 2.95091839 &            & & 0.47788804 & 5.86E--007 & & --0.11161592 & 1.79E--007\smallskip\\
        
$(10,9)$ & 0 &  17100 & 2.95153863 & 2.10E--004 & & 0.47267440 & 1.09E--002 & & --0.10902898 & 2.32E--002\\
         & 1 &  68400 & 2.95091805 & 1.15E--007 & & 0.47784187 & 9.60E--005 & & --0.11161589 & 8.96E--008\\
         & 2 & 273600 & 2.95091839 &            & & 0.47788778 & 4.19E--008 & & --0.11161590 &           \\     
\hline
\end{tabular}

\end{center}
\end{scriptsize}
\renewcommand{\arraystretch}{1}
\end{table}

\begin{figure}[t]
 \centering 
	\includegraphics[width=1.0\textwidth]{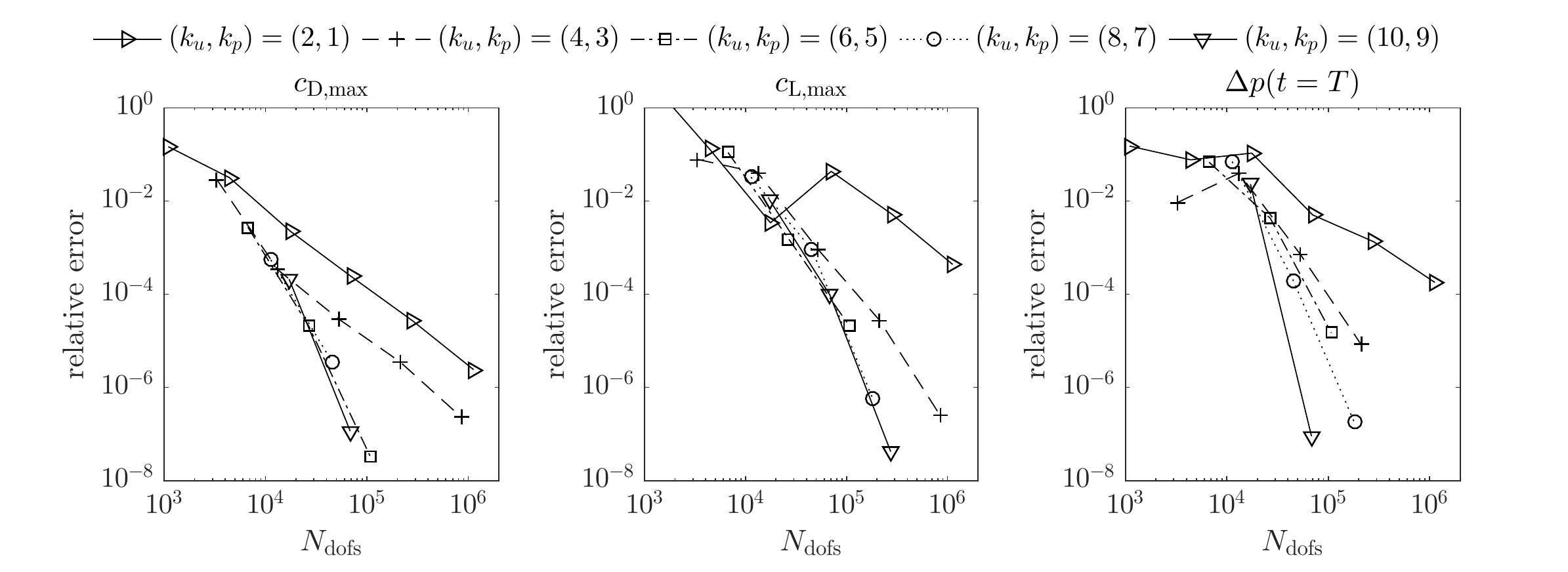}
\caption{$h$-refinement study for various polynomial degrees considering test case 2D-3 of the laminar flow around cylinder benchmark. Relative errors of~$c_{\mathrm{D},\max}$,~$c_{\mathrm{L},\max}$ and~$\Delta p(t=T)$ are shown as a function of the number of unknowns~$N_{\text{dofs}}$ in order to assess the efficiency of high polynomial degrees.}
\label{fig:flow_past_cylinder_efficiency_higher_order}
\end{figure}

The results for the maximum drag and lift coefficients and the pressure difference agree with the reference solution~\cite{Schafer1996}. While the values obtained for maximum lift coefficient and the pressure difference can be seen as exact regarding the reference solution~\cite{John2004}, the maximum drag coefficient converges to a slightly different value. In Table~\ref{FlowPastCylinder2D-3_Reference} we also list an accurate reference solution obtained for the present discretization approach using very fine spatial and temporal resolutions. Since this reference solution is obtained for two different spatial resolutions (refine level~$l=3$ and polynomial degrees~$k=8$ and~$k=10$) as well as different time step sizes, this reference solution is assumed to be accurate up to all decimal places and is used to calculate the relative errors in Table~\ref{SpatialConvergenceFlowPastCylinder2D-3}. The results in Table~\ref{SpatialConvergenceFlowPastCylinder2D-3} indicate that the use of higher order discontinuous Galerkin discretizations allows to obtain significantly more accurate results for a given number of unknowns as compared to low order methods. To further investigate this aspect quantitatively we compare in Figure~\ref{fig:flow_past_cylinder_efficiency_higher_order} the efficiency of different polynomial degrees where we use the quotient of accuracy (inverse of error) and the number of unknowns as a measure of efficiency. An alternative definition of efficiency could be based on the computational costs (in terms of the wall time) rather than the number of unknowns.

The results in Figure~\ref{fig:flow_past_cylinder_efficiency_higher_order} clearly demonstrate that high-order polynomial degrees can significantly reduce the number of unknowns required to obtain a certain level of accuracy and are essential to obtain solutions of high accuracy for a given (maximum) number of unknowns. This can be seen by comparing the curves for~$(k_u,k_p)=(6,5)$ and~$(k_u,k_p)=(2,1)$. At the same time, we observe that the increase in efficiency saturates for polynomial degrees~$(k_u,k_p)=(6,5)$ and higher. Especially for the maximum drag and lift coefficients no noticeable advantage in terms of efficiency can be observed for high polynomial degrees~$(k_u,k_p)=(10,9)$ and~$(k_u,k_p)=(8,7)$ as compared to~$(k_u,k_p)=(6,5)$.

\section{Conclusion}\label{Conclusion}
We have analyzed the stability of projection methods for incompressible flow based on discontinuous Galerkin discretizations and compared the results to a fully coupled, implicit solution approach. Our main results are twofold: Firstly, by means of numerical investigation we have shown that the DG discretization of the velocity--pressure coupling terms substantially affects the stability of projection methods in the limit of small time step sizes. Using a proper DG discretization of these terms including the definition of consistent boundary conditions is crucial to obtain a stable and higher order accurate method. Our results are in contrast to previous publications where instabilities in the limit of small time steps have been ascribed to an inaccurate projection operator resulting in a velocity field that is not exactly divergence-free or to inf--sup instabilities. We emphasize that with the modifications presented in this work we did not observe differences between projection methods and the monolithic solver in terms of stability issues related to the small time steps limit. Secondly, it has been shown that, although some projection methods include an inf--sup stabilization term, spurious pressure oscillations show up for equal-order polynomial approximations resulting in suboptimal rates of convergence in space. Using a mixed-order formulation optimal rates of convergence for velocity and pressure have been demonstrated.

\appendix

\section*{Acknowledgments}
The research presented in this paper was partly funded by the German Research Foundation (DFG) under the project ``High-order discontinuous Galerkin for the exa-scale'' (ExaDG) within the priority program ``Software for Exascale Computing'' (SPPEXA), grant agreement no. KR4661/2-1 and WA1521/18-1.


\bibliography{paper}

\begin{thebibliography}{10}
\expandafter\ifx\csname url\endcsname\relax
  \def\url#1{\texttt{#1}}\fi
\expandafter\ifx\csname urlprefix\endcsname\relax\def\urlprefix{URL }\fi
\expandafter\ifx\csname href\endcsname\relax
  \def\href#1#2{#2} \def\path#1{#1}\fi

\bibitem{Karniadakis13}
G.~E. Karniadakis, S.~J. Sherwin, Spectral/hp element methods for computational
  fluid dynamics, Oxford University Press, 2013.
\newblock \href {http://dx.doi.org/10.1093/acprof:oso/9780198528692.001.0001}
  {\path{doi:10.1093/acprof:oso/9780198528692.001.0001}}.

\bibitem{Hesthaven07}
J.~S. Hesthaven, T.~Warburton, Nodal discontinuous {G}alerkin methods:
  algorithms, analysis, and applications, Springer, 2007.
\newblock \href {http://dx.doi.org/10.1007/978-0-387-72067-8}
  {\path{doi:10.1007/978-0-387-72067-8}}.

\bibitem{Ferrer11}
E.~Ferrer, R.~H.~J. Willden, A high order discontinuous {G}alerkin finite
  element solver for the incompressible {N}avier--{S}tokes equations, Comput.
  Fluids 46~(1) (2011) 224 -- 230.
\newblock \href {http://dx.doi.org/10.1016/j.compfluid.2010.10.018}
  {\path{doi:10.1016/j.compfluid.2010.10.018}}.

\bibitem{Ferrer14}
E.~Ferrer, D.~Moxey, R.~H.~J. Willden, S.~J. Sherwin, Stability of projection
  methods for incompressible flows using high order pressure-velocity pairs of
  same degree: Continuous and discontinuous {G}alerkin formulations, Commun.
  Comput. Phys. 16 (2014) 817--840.
\newblock \href {http://dx.doi.org/10.4208/cicp.290114.170414a}
  {\path{doi:10.4208/cicp.290114.170414a}}.

\bibitem{Karniadakis1991}
G.~E. Karniadakis, M.~Israeli, S.~A. Orszag, High-order splitting methods for
  the incompressible {N}avier--{S}tokes equations, J. Comput. Phys. 97~(2)
  (1991) 414 -- 443.
\newblock \href {http://dx.doi.org/10.1016/0021-9991(91)90007-8}
  {\path{doi:10.1016/0021-9991(91)90007-8}}.

\bibitem{Guermond06}
J.~L. Guermond, P.~Minev, J.~Shen, An overview of projection methods for
  incompressible flows, Comput. Methods in Appl. Mech. Eng. 195~(44–47)
  (2006) 6011 -- 6045.
\newblock \href {http://dx.doi.org/10.1016/j.cma.2005.10.010}
  {\path{doi:10.1016/j.cma.2005.10.010}}.

\bibitem{Orszag1986}
S.~A. Orszag, M.~Israeli, M.~O. Deville, Boundary conditions for incompressible
  flows, J. Sci. Comput. 1~(1) (1986) 75--111.
\newblock \href {http://dx.doi.org/10.1007/BF01061454}
  {\path{doi:10.1007/BF01061454}}.

\bibitem{Guermond2003}
J.-L. Guermond, J.~Shen, Velocity-correction projection methods for
  incompressible flows, SIAM J. Numer. Anal. 41~(1) (2003) 112--134.
\newblock \href {http://dx.doi.org/10.1137/S0036142901395400}
  {\path{doi:10.1137/S0036142901395400}}.

\bibitem{Chorin68}
A.~J. Chorin, Numerical solution of the {N}avier--{S}tokes equations, Math.
  Comp. 22~(104) (1968) 745--762.
\newblock \href {http://dx.doi.org/10.1090/S0025-5718-1968-0242392-2}
  {\path{doi:10.1090/S0025-5718-1968-0242392-2}}.

\bibitem{hirt1972}
C.~Hirt, J.~Cook, Calculating three-dimensional flows around structures and
  over rough terrain, Journal of Computational Physics 10~(2) (1972) 324--340.

\bibitem{Goda1979}
K.~Goda, A multistep technique with implicit difference schemes for calculating
  two-or three-dimensional cavity flows, Journal of Computational Physics
  30~(1) (1979) 76--95.

\bibitem{VanKan1986}
J.~Van~Kan, A second-order accurate pressure-correction scheme for viscous
  incompressible flow, SIAM Journal on Scientific and Statistical Computing
  7~(3) (1986) 870--891.

\bibitem{Timmermans1996}
L.~Timmermans, P.~Minev, F.~Van De~Vosse, An approximate projection scheme for
  incompressible flow using spectral elements, International Journal for
  Numerical Methods in Fluids 22~(7) (1996) 673--688.

\bibitem{Guermond2004}
J.~Guermond, J.~Shen, On the error estimates for the rotational
  pressure-correction projection methods, Mathematics of Computation 73~(248)
  (2004) 1719--1737.

\bibitem{cockburn2002local}
B.~Cockburn, G.~Kanschat, D.~Sch{\"o}tzau, C.~Schwab, Local discontinuous
  {G}alerkin methods for the {S}tokes system, SIAM Journal on Numerical
  Analysis 40~(1) (2002) 319--343.

\bibitem{cockburn2004local}
B.~Cockburn, G.~Kanschat, D.~Sch{\"o}tzau, The local discontinuous {G}alerkin
  method for the {O}seen equations, Mathematics of Computation 73~(246) (2004)
  569--593.

\bibitem{cockburn2005locally}
B.~Cockburn, G.~Kanschat, D.~Sch{\"o}tzau, A locally conservative {L}{D}{G}
  method for the incompressible {N}avier--{S}tokes equations, Mathematics of
  Computation 74~(251) (2005) 1067--1095.

\bibitem{cockburn2009equal}
B.~Cockburn, G.~Kanschat, D.~Sch{\"o}tzau, An equal-order {D}{G} method for the
  incompressible {N}avier--{S}tokes equations, Journal of Scientific Computing
  40~(1-3) (2009) 188--210.

\bibitem{girault2005discontinuous}
V.~Girault, B.~Rivi{\`e}re, M.~Wheeler, A discontinuous {G}alerkin method with
  nonoverlapping domain decomposition for the {S}tokes and {N}avier--{S}tokes
  problems, Mathematics of Computation 74~(249) (2005) 53--84.

\bibitem{girault2005splitting}
V.~Girault, B.~Rivi{\`e}re, M.~F. Wheeler, A splitting method using
  discontinuous {G}alerkin for the transient incompressible {N}avier--{S}tokes
  equations, ESAIM: Mathematical Modelling and Numerical
  Analysis-Mod{\'e}lisation Math{\'e}matique et Analyse Num{\'e}rique 39~(6)
  (2005) 1115--1147.

\bibitem{bassi2006artificial}
F.~Bassi, A.~Crivellini, D.~A. Di~Pietro, S.~Rebay, An artificial
  compressibility flux for the discontinuous {G}alerkin solution of the
  incompressible {N}avier--{S}tokes equations, Journal of Computational Physics
  218~(2) (2006) 794--815.

\bibitem{bassi2007implicit}
F.~Bassi, A.~Crivellini, D.~A. Di~Pietro, S.~Rebay, An implicit high-order
  discontinuous {G}alerkin method for steady and unsteady incompressible flows,
  Computers \& Fluids 36~(10) (2007) 1529--1546.

\bibitem{Shahbazi07}
K.~Shahbazi, P.~F. Fischer, C.~R. Ethier, A high-order discontinuous {G}alerkin
  method for the unsteady incompressible {N}avier--{S}tokes equations, J.
  Comput. Phys. 222~(1) (2007) 391 -- 407.
\newblock \href {http://dx.doi.org/10.1016/j.jcp.2006.07.029}
  {\path{doi:10.1016/j.jcp.2006.07.029}}.

\bibitem{Klein13}
B.~Klein, F.~Kummer, M.~Oberlack, A {SIMPLE} based discontinuous {G}alerkin
  solver for steady incompressible flows, J. Comput. Phys. 237 (2013) 235 --
  250.
\newblock \href {http://dx.doi.org/10.1016/j.jcp.2012.11.051}
  {\path{doi:10.1016/j.jcp.2012.11.051}}.

\bibitem{Klein15}
B.~Klein, F.~Kummer, M.~Keil, M.~Oberlack, An extension of the {SIMPLE} based
  discontinuous {G}alerkin solver to unsteady incompressible flows, Int. J.
  Numer. Meth. Fluids 77~(10) (2015) 571--589, fld.3994.
\newblock \href {http://dx.doi.org/10.1002/fld.3994}
  {\path{doi:10.1002/fld.3994}}.

\bibitem{Botti11}
L.~Botti, D.~A.~D. Pietro, A pressure-correction scheme for
  convection-dominated incompressible flows with discontinuous velocity and
  continuous pressure, J. Comput. Phys. 230~(3) (2011) 572 -- 585.
\newblock \href {http://dx.doi.org/10.1016/j.jcp.2010.10.004}
  {\path{doi:10.1016/j.jcp.2010.10.004}}.

\bibitem{Piatkowski16}
M.~Piatkowski, S.~M\"uthing, P.~Bastian, A {S}table and {H}igh-{O}rder
  {A}ccurate {D}iscontinuous {G}alerkin {B}ased {S}plitting {M}ethod for the
  {I}ncompressible {N}avier--{S}tokes {E}quations, arXiv preprint
  arXiv:1612.00657.

\bibitem{nguyen2011implicit}
N.~C. Nguyen, J.~Peraire, B.~Cockburn, An implicit high-order hybridizable
  discontinuous {G}alerkin method for the incompressible {N}avier--{S}tokes
  equations, Journal of Computational Physics 230~(4) (2011) 1147--1170.

\bibitem{Lehrenfeld16}
C.~Lehrenfeld, J.~Sch{\"o}berl, High order exactly divergence-free hybrid
  discontinuous {G}alerkin methods for unsteady incompressible flows, Comput.
  Methods in Appl. Mech. Eng. 307 (2016) 339 -- 361.
\newblock \href {http://dx.doi.org/10.1016/j.cma.2016.04.025}
  {\path{doi:10.1016/j.cma.2016.04.025}}.

\bibitem{Steinmoeller13}
D.~T. Steinmoeller, M.~Stastna, K.~G. Lamb, A short note on the discontinuous
  {G}alerkin discretization of the pressure projection operator in
  incompressible flow, J. Comput. Phys. 251 (2013) 480 -- 486.
\newblock \href {http://dx.doi.org/10.1016/j.jcp.2013.05.036}
  {\path{doi:10.1016/j.jcp.2013.05.036}}.

\bibitem{Krank16b}
B.~Krank, N.~Fehn, W.~A. Wall, M.~Kronbichler, A high-order semi-explicit
  discontinuous {G}alerkin solver for 3{D} incompressible flow with application
  to {DNS} and {LES} of turbulent channel flow, arXiv preprint
  arXiv:1607.01323.

\bibitem{Joshi16}
S.~M. Joshi, P.~J. Diamessis, D.~T. Steinmoeller, M.~Stastna, G.~N. Thomsen, A
  post-processing technique for stabilizing the discontinuous pressure
  projection operator in marginally-resolved incompressible inviscid flow,
  Comput. Fluids\href {http://dx.doi.org/10.1016/j.compfluid.2016.04.021}
  {\path{doi:10.1016/j.compfluid.2016.04.021}}.

\bibitem{Emamy14}
N.~Emamy, Numerical simulation of deformation of a droplet in a stationary
  electric field using dg, Ph.D. thesis, Technische Universit{\"a}t Darmstadt
  (2014).

\bibitem{Emamy17}
N.~Emamy, F.~Kummer, M.~Mrosek, M.~Karcher, M.~Oberlack, Implicit-explicit and
  explicit projection schemes for the unsteady incompressible
  {N}avier--{S}tokes equations using a high-order d{G} method, Computers \&
  Fluids Accepted.

\bibitem{Leriche2000}
E.~Leriche, G.~Labrosse, High-order direct {S}tokes solvers with or without
  temporal splitting: numerical investigations of their comparative properties,
  SIAM J. Sci. Comput. 22~(4) (2000) 1386--1410.
\newblock \href {http://dx.doi.org/10.1137/S1064827598349641}
  {\path{doi:10.1137/S1064827598349641}}.

\bibitem{Leriche2006}
E.~Leriche, E.~Perchat, G.~Labrosse, M.~O. Deville, Numerical evaluation of the
  accuracy and stability properties of high-order direct {S}tokes solvers with
  or without temporal splitting, J. Sci. Comput. 26~(1) (2006) 25--43.
\newblock \href {http://dx.doi.org/10.1007/s10915-004-4798-0}
  {\path{doi:10.1007/s10915-004-4798-0}}.

\bibitem{arnold2000discontinuous}
D.~N. Arnold, F.~Brezzi, B.~Cockburn, D.~Marini, Discontinuous {G}alerkin
  methods for elliptic problems, in: B.~Cockburn, G.~Karniadakis, C.-W. Shu
  (Eds.), Discontinuous Galerkin Methods, Vol.~11 of Lecture Notes in
  Computational Science and Engineering, Springer Berlin Heidelberg, 2000, pp.
  89--101.

\bibitem{arnold2002unified}
D.~N. Arnold, F.~Brezzi, B.~Cockburn, L.~D. Marini, Unified analysis of
  discontinuous {G}alerkin methods for elliptic problems, SIAM Journal on
  Numerical Analysis 39~(5) (2002) 1749--1779.

\bibitem{Shahbazi05}
K.~Shahbazi, An explicit expression for the penalty parameter of the interior
  penalty method, J. Comput. Phys. 205~(2) (2005) 401 -- 407.
\newblock \href {http://dx.doi.org/10.1016/j.jcp.2004.11.017}
  {\path{doi:10.1016/j.jcp.2004.11.017}}.

\bibitem{Hillewaert13}
K.~Hillewaert, Development of the discontinuous {G}alerkin method for
  high-resolution, large scale {CFD} and acoustics in industrial geometries,
  Ph.D. thesis, Univ. de Louvain (2013).

\bibitem{dealII85}
D.~Arndt, W.~Bangerth, D.~Davydov, T.~Heister, L.~Heltai, M.~Kronbichler,
  M.~Maier, J.-P. Pelteret, B.~Turcksin, D.~Wells, The \texttt{deal.II}
  library, version 8.5, Journal of Numerical Mathematics\href
  {http://dx.doi.org/10.1515/jnma-2017-0058}
  {\path{doi:10.1515/jnma-2017-0058}}.

\bibitem{Kronbichler12}
M.~Kronbichler, K.~Kormann, A generic interface for parallel cell-based finite
  element operator application, Comput. Fluids 63 (2012) 135--147.
\newblock \href {http://dx.doi.org/10.1016/j.compfluid.2012.04.012}
  {\path{doi:10.1016/j.compfluid.2012.04.012}}.

\bibitem{Kormann16}
K.~Kormann, M.~Kronbichler, Efficient matrix-free implementations for
  discontinuous {G}alerkin methods, In preparation.

\bibitem{Ethier1994}
C.~R. Ethier, D.~Steinman, Exact fully 3{D} {N}avier--{S}tokes solution for
  benchmarking, International Journal for Numerical Methods in Fluids 19~(5)
  (1994) 369--376.

\bibitem{Schafer1996}
M.~Sch{\"a}fer, S.~Turek, F.~Durst, E.~Krause, R.~Rannacher, Benchmark
  computations of laminar flow around a cylinder, Springer, 1996.

\bibitem{John2004}
V.~John, Reference values for drag and lift of a two-dimensional time-dependent
  flow around a cylinder, International Journal for Numerical Methods in Fluids
  44~(7) (2004) 777--788.

\end{thebibliography}

\end{document}